\documentclass[a4paper, 11pt, reqno]{article}

\pdfoutput=1

\usepackage{amsmath}
\usepackage{amsfonts}
\usepackage{amssymb}
\usepackage{amsthm}
\usepackage{authblk}
\usepackage{bm}
\usepackage{indentfirst}
\usepackage{graphicx}
\usepackage{subfigure}
\usepackage{array}
\usepackage{fullpage}
\usepackage{color}
\usepackage{algorithm}
\usepackage{algorithmic}
\usepackage{url}

\DeclareMathAlphabet{\mathsfsl}{OT1}{cmss}{m}{sl}
\newcommand{\tensor}[1]{\mathbf{#1}}

\newcommand{\dif}{\mathrm{d}}

\newcommand{\ba}{\bm\alpha}
\newcommand{\bx}{\bm\xi}
\newcommand{\ve}{\varepsilon}
 
\newcommand{\mexp}[1]{\mathbb{E}\left\{{#1}\right\}}

\newcommand{\mr}{\mathbb{R}}
\newcommand{\mn}{\mathbb{N}}
\newcommand{\hc}{\hat{\bm c}} 
\newcommand{\tc}{\tilde{\bm c}} 
\newcommand{\trans}{^\top} 
\DeclareMathOperator*{\argmin}{arg\,min}
\newcommand{\PreserveBackslash}[1]{\let\temp=\\#1\let\\=\temp}
\newcolumntype{C}[1]{>{\PreserveBackslash\centering}p{#1}}
\newcolumntype{R}[1]{>{\PreserveBackslash\raggedleft}p{#1}}
\newcolumntype{L}[1]{>{\PreserveBackslash\raggedright}p{#1}}

\newtheorem{thm}{Theorem}[section]
\newtheorem{lem}[thm]{Lemma}

\newcommand\Def{\stackrel{\textrm{def}}{=}}
\DeclareMathOperator{\sech}{sech}

\title{A General Framework for Enhancing Sparsity of Generalized Polynomial Chaos Expansions}

\author{Xiu Yang\footnote{(Corresponding author) xiu.yang@pnnl.gov}}
\affil{Advanced Computing, Mathematics and Data Division, Pacific Northwest
National Laboratory, Richland, WA, 99352.}
\author{Xiaoliang Wan\footnote{xlwan@math.lsu.edu}}
\affil{Department of Mathematics and Center of Computation and Technology,
Louisiana State University, Baton Rouge, LA, 70803.}
\author{Lin Lin\footnote{linlin@math.berkeley.edu}}
\affil{Department of Mathematics, University of California, Berkeley and Computational
Research Division, Lawrence Berkeley National Laboratory, Berkeley, CA 94720.}
\author{Huan Lei\footnote{huan.lei@pnnl.gov}}
\affil{Advanced Computing, Mathematics and Data Division, Pacific Northwest
National Laboratory, Richland, WA, 99352.}

\begin{document}

\maketitle

\abstract{
Compressive sensing has become a powerful addition to uncertainty quantification
when only limited data are available. In this paper, we provide a general 
framework to enhance the sparsity of the representation of uncertainty in the 
form of generalized polynomial chaos expansion. We use an alternating direction
method to identify new sets of random variables through iterative rotations so 
the new representation of the uncertainty is sparser. Consequently, we increase
both the efficiency and accuracy of the compressive-sensing-based uncertainty
quantification method. We demonstrate that the previously developed
rotation-based methods \cite{LeiYZLB15,YangLBL16} to enhance the sparsity of 
Hermite polynomial expansion is a special case of this general framework. 
Moreover, we use Legendre and Chebyshev polynomial expansions to demonstrate the
effectiveness of this method with applications in solving stochastic partial 
differential equations and high-dimensional ($\mathcal{O}(100)$) problems.

\textbf{keywords}: Uncertainty quantification, generalized polynomial chaos, 
compressive sensing, iterative rotation, alternating direction.
}

\section{Introduction}
\label{sec:intro}
Surrogate-model-based uncertainty quantification (UQ) plays an important role in constructing
computational models as it helps to understand the influence of uncertainties on
the quantity of interest (QoI). In this paper, we focus on parametric uncertainty,
i.e., some of the parameters in the system are random variables.
We assume that these random variables are defined on a probability space
$(\Omega,\mathcal{F},P)$, where $\Omega$ is the event space and $P$ is a
probability measure on the $\sigma$-field $\mathcal{F}$. We consider a system 
depending on a $d$-dimensional random vector
$\bx(\omega)=(\xi_1(\omega),\xi_2(\omega),\cdots,\xi_d(\omega))\trans$, where
$\omega$ is an event in $\Omega$. For simplicity, we denote $\xi_i(\omega)$ as
$\xi_i$. A useful surrogate model of the QoI $u(\bx)$ is the 
generalized polynomial chaos (gPC) expansion \cite{GhanemS91, XiuK02}:
\begin{equation}\label{eq:gpc}
u(\bx) = \sum_{n=1}^Nc_n\psi_n(\bx) + \ve(\bx),
\end{equation}
where $\ve$ is the truncation error, $N$ is a positive integer, $c_n$ are
coefficients, and $\psi_n$ are multivariate polynomials which are orthonormal 
with respect to the measure of $\bx$:
\begin{equation}
  \int_{\mr^d} \psi_i(\bm x)\psi_j(\bm x)\rho_{\bx}(\bm x)\dif\bm x = \delta_{ij},
\end{equation}
where $\rho_{\bx}(\bm x)$ is the probability density function (PDF) of $\bx$ and 
$\delta_{ij}$ is the Kronecker delta function. This approximation converges in the $L_2$ 
sense as $N$ increases if $u$ is in the Hilbert space associated with the 
measure of $\bx$ (i.e., the weight of the inner product is the PDF of $\bx$) 
\cite{XiuK02, CameronM1947, Ogura1972, ernst2012convergence}. 
Both \emph{intrusive} methods (e.g., stochastic Galerkin) and 
\emph{non-intrusive} methods (e.g., probabilistic collocation method)
are developed \cite{GhanemS91,XiuK02,TatangPPM97, XiuH05, BabuskaNT10} to compute the 
gPC coefficients $\bm c=(c_1,c_2,\cdots,c_N)\trans$. Specifically, the non-intrusive
method is more suited to study complex system as it does not require
modifying the computational models or simulation, codes while the intrusive 
method requires rewriting these models or codes. 

Non-intrusive methods use the samples of input $\{\bx^q\}_{q=1}^M$ and 
corresponding output of the computational model $\{u^q=u(\bx^q)\}_{q=1}^M$ to 
compute the coefficients $\bm c$. Subsequently, we can write Eq.~\eqref{eq:gpc}
as the following linear system:
\begin{equation}\label{eq:cs_eq}
\tensor\Psi \bm c = \bm u - \bm\ve,
\end{equation}
where $\bm u=(u^1,u^2,\cdots,u^M)\trans$ is the vector of output samples, 
$\bm c=(c_1, c_2, \cdots, c_N)\trans$ is the vector of gPC coefficients,
$\tensor\Psi$ is an $M\times N$ matrix with $\Psi_{ij}=\psi_j(\bx^i)$ and 
$\bm\ve=(\ve^1,\ve^2,\cdots,\ve^M)\trans$ is a vector of error samples with 
$\ve^q=\ve(\bx^q)$. In many practical problems, it is costly to obtain $u^q$
because of the model complexity and the limited computational
{\color{red}resources}. As such,
we often must consider the scenario when $M<N$ or even $M\ll N$. Namely, the
number of samples is smaller than the number of basis functions, which implies
that Eq.~\eqref{eq:cs_eq} is underdetermined. The compressive sensing method is 
effective at solving this type of underdetermined problem when $\bm c$ is sparse 
\cite{CandesRT06,DonohoET06,Candes08,BrucksteinDE09} and recent studies have 
applied this approach to UQ problems 
\cite{DoostanO11,YanGX12,YangK13,LeiYZLB15,XuZ14,SargsyanSNDRT14,PengHD15,LeiYLK17}. 
Several useful approaches have been developed to improve the efficiency of 
solving Eq.~\eqref{eq:cs_eq} in UQ applications. For example, 
weighted/re-weighted $\ell_1$ minimization assigns a weight to each $c_n$ and 
solves a weighted $\ell_1$ minimization problem to enhance the 
sparsity \cite{CandesWB08,YangK13,PengHD14,RauhutW15}; better sampling 
strategies can improve the property of 
$\tensor\Psi$ \cite{RauhutW12,HamptonD15,AlemazkoorM17,jakeman2017generalized};
adaptive basis selection reduces the number of unknown \cite{JakemanES14,dai2008subspace,blatman2011adaptive,alemazkoor2017divide,hampton2018basis}; dimension reduction techniques
also reduce the number of unknown \cite{yang2017sliced,tsilifis2018compressive}.

We proposed to enhance the sparsity $\bm c$ for an arbitrarily 
distributed random variable using a non-linear mapping via optimization 
in \cite{YangWK13}. Later, motivated by quantifying the QoIs of complex 
biomolecule systems, we proposed an approach to enhance the sparsity of $\bm c$
by unitary transform of the random vector $\bm\xi$ \cite{LeiYZLB15} when $\xi_i$
are independent and identically distribution (i.i.d.) Gaussian random variables.
Subsequently, we developed an iterative-rotation algorithm \cite{YangLBL16} to 
enhance the sparsity of $\bm c$ more effectively in a successive manner. In this 
work, we provide a general framework to enhance the sparsity of the 
representation of uncertainty in the form of gPC expansion, and $\bx$ are not 
restricted to Gaussian random variables. We use an alternating direction method
to iteratively identify a rotation matrix $\tensor A$, which transforms $\bx$ to
a new set of random variables 
$\bm\eta=\tensor A\bx$, where $\bm\eta=(\eta_1, \eta_2, \cdots, \eta_d)\trans$,
such that the ``gPC expansion" of $u$ with respect to $\bm\eta$ is sparser. In
other words,
\begin{equation}\label{eq:gpc2}
u(\bx) \approx \sum_{n=1}^N c_n\psi_n(\bx)=\sum_{n=1}^N \tilde c_n \psi_n(\bm\eta(\bx))\approx u(\bm\eta(\bx)),
\end{equation}
and we intend to make 
$\tilde{\bm c}=(\tilde c_1, \tilde c_2, \cdots, \tilde c_N)\trans$ 
sparser than $\bm c$. By increasing the sparsity, we can improve both the 
efficiency and accuracy of the compressive-sensing-based UQ
method because we can use fewer samples of $u^q$ to obtain a more 
accurate representation of $u$. Particularly, making $\tilde{\bm c}$ sparser
than $\bm c$ has the potential to yield
{\color{red}
\[\left\Vert u(\bx)-\sum_{n=1}^N \hat{\tilde c}_n\psi_n(\tensor A\bx)\right\Vert_2\leq
\left\Vert u(\bx)-\sum_{n=1}^N \hat c_n\psi_n(\bx)\right\Vert_2,\]
when $\hat{\bm c}$ and $\hat{\tilde{\bm c}}$ are approximations of $\bm c$ and
$\tilde{\bm c}$ by compressive sensing, respectively. This is because the 
compressive sensing method computes sparser coefficients more accurately with
a given samples size (see Eq.~\eqref{eq:l1_thm} for the error estimate).}
Of note, $\psi_n$ may not be orthonormal with respect to the measure of 
$\bm\eta$. Thus, the new representation of $u$ with respect to $\bm\eta$, 
$\displaystyle\sum_{n=1}^N \tilde c_n \psi_n(\bm\eta)$, is not necessarily a 
standard gPC expansion. Instead, it is a polynomial-based surrogate model for 
$u$. We use Legendre and Chebyshev polynomial expansions to demonstrate the 
effectiveness of our proposed method.

\section{Brief review of the compressive-sensing-based gPC method}
\label{sec:review}

\subsection{Generalized polynomial chaos expansions}
In this paper, we study systems relying on $d$-dimensional i.i.d. random
variables $\bx$. Hence, the gPC basis
functions are constructed by tensor products of univariate orthonormal 
polynomials. For a multi-index 
$\ba=(\alpha_1,\alpha_2,\cdots,\alpha_d), \alpha_i\in\mn\cup\{0\}$, we set
\begin{equation}\label{eq:tensor}
\psi_{\ba}(\bx) =
\psi_{\alpha_1}(\xi_1)\psi_{\alpha_2}(\xi_2)\cdots\psi_{\alpha_d}(\xi_d).
\end{equation}
For two different multi-indices
$\ba_i=((\alpha_i)_{_1}, \cdots, (\alpha_i)_{_d})$ and 
$\ba_j=((\alpha_j)_{_1}, \cdots, (\alpha_j)_{_d})$, we have the
property
\begin{equation}
  \int_{\mr^d} \psi_{\ba_i}(\bm x)\psi_{\ba_j}(\bm x) \rho_{\bx}(\bm x) \dif\bm x=
\delta_{\ba_i\ba_j} = \delta_{(\alpha_i)_{_1}(\alpha_j)_{_1}}
\delta_{(\alpha_i)_{_2}(\alpha_j)_{_2}}\cdots
\delta_{(\alpha_i)_{_d}(\alpha_j)_{_d}},
\end{equation}
where 
\begin{equation}
  \rho_{\bx}(\bm x) = \rho_{\xi_1}(x_1)\rho_{\xi_2}(x_2)\cdots\rho_{\xi_d}(x_d).
\end{equation}
For simplicity, we denote $\psi_{\ba_i}(\bx)$ as $\psi_i(\bx)$.


\subsection{Compressive sensing}

We first introduce the concept of \emph{sparsity} as it is critical in the error
estimates for solving the under-determined system Eq.~\eqref{eq:cs_eq} with the
compressive sensing method. The number of non-zero entries of a vector 
$\bm x=(x_1,x_2,\cdots,x_N)$ is denoted as 
$\Vert\bm x\Vert_0\Def \#\{i:x_i\neq 0\}$ \cite{Donoho06,CandesRT06,BrucksteinDE09},
and the $\ell_1$ norm of $\bm x$ is defined as
$\Vert\bm x\Vert_1\Def \sum_{n=1}^N |x_n|$.
Of note, $\Vert\cdot\Vert_0$ is named ``$\ell_0$ norm" in \cite{Donoho06},
although it is not a norm nor a semi-norm. The vector $\bm x$ is called 
\emph{$s$-sparse} if $\Vert \bm x\Vert_0\leq s$, and $\bm x$ 
is considered a sparse vector if $s\ll N$. Few practical systems have truly 
sparse gPC coefficients $\bm c$. However, in many cases, the $\bm c$ are 
compressible, i.e., only a few entries make significant contribution to its
$\ell_1$ norm. Consequently, a vector $\bm x_s$ is defined as the 
\textbf{best $s$-sparse approximation} that can be obtained knowing the exact
locations and amplitudes of the $s$-largest entries of $\bm x$, i.e., $\bm x_s$
is the vector $\bm x$ with all but the $s$-largest entries set to 
zero \cite{Candes08}. Subsequently, $\bm x$ is considered sparse if 
$\Vert \bm x - \bm x_s\Vert_1$ is small for $s\ll N$. 

Under some conditions, the sparse vector $\bm c$ in Eq.~\eqref{eq:cs_eq} can be 
approximated by solving the following $\ell_1$ minimization problem:
\begin{equation}\label{eq:lh}
(P_{1,\epsilon}):~\arg \min_{\hat{\bm c}}\Vert\hat{\bm c}\Vert_1, 
\text{~~subject to~~} \Vert\tensor\Psi\hat{\bm c}-\bm u\Vert_2\leq\epsilon,
\end{equation}
where $\epsilon=\Vert\bm\ve\Vert_2$. 
The error bound for solving Eq.~\eqref{eq:cs_eq} with $\ell_1$ minimization 
requires definiting the \textit{restricted isometry property} (RIP) constant
\cite{CandesT05}. For each integer $s=1,\cdots,N$, the isometry constant 
$\delta_s$ of a matrix $\tensor\Phi$ is defined as the smallest number such that
\begin{equation}
(1-\delta_s)\Vert\bm x\Vert_2^2\leq \Vert\tensor\Phi\bm x\Vert_2^2\leq
(1+\delta_s)\Vert\bm x\Vert_2^2
\end{equation}
holds for all $s$-sparse vectors $\bm x$.
\noindent With some restrictions, Cand\`{e}s et al.\ showed $\bm x$ can be stably
reconstructed \cite{Candes08}.
Assume that the matrix $\tensor\Psi$ satisfies $\delta_{2s}<\sqrt{2}-1$, and
$\Vert\bm\ve\Vert_2\leq\epsilon$, then solution 
$\hat{\bm c}$ to $(P_{1,\epsilon})$ obeys
\begin{equation}\label{eq:l1_thm}
\Vert \bm c - \hat{\bm c}\Vert_2 \leq C_1\epsilon + C_2 
\dfrac{\Vert\bm c-\bm c_s\Vert_1}{\sqrt{s}},
\end{equation}
where $C_1$ and $C_2$ are constants, $\bm c$ is the exact vector we aim to 
approximate and $\hat{\bm c}$ is the solution of $(P_{1,\epsilon})$. 
This result implies that the upper bound of the error is related to the 
truncation error and the sparsity of $\bm c$, which is indicated in the first
and second terms on the right-hand side of Eq. \eqref{eq:l1_thm}, respectively.
We use  $\Vert\bm c-\bm c_s\Vert_1/\sqrt{s}$ to examine the sparsity in our 
numerical examples.

In practice, the error term $\epsilon$ is not known \emph{a priori}. Hence,
in the present work, we use cross-validation to estimate it. One such algorithm
is \cite{DoostanO11} summarized in Algorithm \ref{algo:cross}.
\begin{algorithm}
\caption{Cross-validation to estimate the error $\epsilon$}
\label{algo:cross}
\begin{algorithmic}[1]
\STATE Divide the $M$ output samples to $M_r$ reconstruction ($\bm u_r$) and
$M_v$ validation ($\bm u_v$) samples and divide the measurement matrix 
$\tensor\Psi$ correspondingly into $\tensor\Psi_r$ and $\tensor\Psi_v$.
\STATE Choose multiple values for $\epsilon_r$ such that the exact error
$\Vert\tensor\Psi_r\bm c-\bm u_r\Vert_2$ of the reconstruction samples is
within the range of $\epsilon_r$ values.
\STATE For each $\epsilon_r$, {\color{red}solve $(P_{1,\epsilon})$} with $\bm u_r$ and
$\tensor\Psi_r$ to obtain $\hat{\bm c}$, then compute 
$\epsilon_v=\Vert\tensor\Psi_v\hat{\bm c}-\bm u_v\Vert_2$.
\STATE Find the minimum value of $\epsilon_v$ and its corresponding 
$\epsilon_r$. Set $\epsilon=\sqrt{M/M_r}\epsilon_r$.
\end{algorithmic}
\end{algorithm}
We note that some techniques may be applied to avoid the cross-validation step
and we refer interested readers to \cite{Adcock17}. 


\subsection{Compressive-sensing-based gPC methods}
Given $M$ samples of $\bx$, the QoI $u$ is approximated by a gPC expansion as in
Eq.~\eqref{eq:gpc}:
\begin{equation}
u(\bx^q) = \sum_{n=1}^N c_n\psi_n(\bx^q) + \ve(\bx^q), \quad q=1,2,\cdots,M,
\end{equation}
which can be rewritten as Eq.~\eqref{eq:cs_eq}. A typical approach to 
compressive-sensing-based-gPC is summarized in Algorithm \ref{algo:cs1}.
\begin{algorithm}
\caption{Compressive-sensing-based gPC}
\label{algo:cs1}
\begin{algorithmic}[1]
\STATE Generate input samples $\bx^q, q=1,2,\cdots, M$ based on the distribution
of $\bx$. 
\STATE Generate output samples $u^q=u(\bx^q)$ by solving the complete model, 
e.g., running simulations, solvers, etc.
\STATE Select gPC basis functions $\{\psi_n\}_{n=1}^N$ associated with $\bx$ and
then generate the measurement matrix $\tensor\Psi$ by setting
$\Psi_{ij}=\psi_j(\bx^i)$.
\STATE Solve the optimization problem $(P_{h,\epsilon})$:
\[\arg \min_{\hat{\bm c}}\Vert\hat{\bm c}\Vert_h, ~
\text{subject to}~ \Vert\tensor\Psi\hat{\bm c}-\bm u\Vert_2\leq\epsilon,\]
where $h=0$ or $1$, $\bm u=(u^1,u^2,\cdots,u^M)^T$, and $\epsilon$ is obtained by
cross-validation. 
\STATE Set $\bm c=\hat{\bm c}$, and construct gPC expansion as 
$u(\bx)\approx \sum_{n=1}^N c_n\psi_n(\bx)$.
\end{algorithmic}
\end{algorithm}
Moreover, we use the \emph{re-weighted} $\ell_1$ minimization approach 
\cite{CandesWB08} in the numerical examples to improve the accuracy  of solving
Eq.~\eqref{eq:lh}. This approach solves the following optimization problem:
\begin{equation}\label{eq:wl1}
(P_{1,\epsilon}^W):~\arg\min_{\hat{\bm c}}\Vert\tensor W\hat{\bm c}\Vert_1, 
~\text{subject to}~\Vert\tensor\Psi\hat{\bm c}-\bm u\Vert_2\leq\epsilon,
\end{equation}
where $\tensor W$ is a diagonal matrix: 
$\tensor W=\text{diag}(w_1,w_2,\cdots,w_N)$. Clearly, $(P_{1,\epsilon})$ can be
considered as a special case of $(P_{1,\epsilon}^W)$ by setting 
$\tensor W=\tensor I$. The elements $w_i$ of the diagonal matrix can be 
estimated iteratively \cite{CandesWB08, YangK13}. More precisely, in the $l$-th
iteration, $(P_{1,\epsilon}^W)$ is solved to obtain $\hat{\bm c}^{(l)}$.
Then we set $w_i^{(l+1)}=1/(|\hat c_i^{(l)}|+\delta)$ for the next iteration. 
The parameter $\delta>0$ is introduced to provide stability and to ensure that a
zero-valued component in $\hat{\bm c}^{(l)}$ does not prohibit a nonzero 
estimate at the next step. In Cand\`{e}s et al.~\cite{CandesWB08}, the authors 
suggest two to three iterations of this procedure. Subsequent analytical 
work \cite{Needell09} provides an error bound for each iteration, as well as
the limit of computing $\hat{\bm c}$ with re-weighted $\ell_1$ minimization. 
The form is similar to Eq.~\eqref{eq:l1_thm} with different constants.


\section{Iterative rotations}
\label{sec:method}

\subsection{Basic idea}
\label{subsec:basic}
We aim to find a linear map $g:\mr^d\mapsto\mr^d$ such that we have a new 
set of random variables:
\begin{equation}
  \bm\eta=g(\bx)=\tensor A\tensor\bx, \quad \bm\eta =
  (\eta_1,\eta_2,\cdots,\eta_d)\trans,
\end{equation}
where $\tensor A$ is an orthogonal matrix satisfying $\tensor A\tensor A\trans=\tensor I$
and the PDF of $\bm\eta$ is denoted as $\rho_{\bm\eta}$. In other words, the map
from $\bx$ to $\bm\eta$ can be considered as a rotation in $\mathbb{R}^d$.
As such, the new polynomial expansion for $u$ is
\begin{equation}\label{eq:new_poly}
u(\bx) \approx u_g(\bx)=\sum_{n=1}^Nc_n\psi_n(\bx)=
\sum_{n=1}^N\tilde c_n\psi_n(\tensor A\bx)=
\sum_{n=1}^N\tilde c_n\psi_n(\bm\eta)=v_g(\bm\eta),
\end{equation}
namely, 
\[u(\bx)=v(\bm\eta(\bx))=v(\bm\eta)\approx v_g(\bm\eta).\]
Here, $u_g(\bx)$ is understood as a polynomial $u_g(\bm x)$ evaluated at the random 
variables $\bx$ and the same for $v_g$. Ideally, $\tilde{\bm c}$ is sparser than 
$\bm c$. In the previous work \cite{LeiYZLB15,YangLBL16}, we assume that 
$\bm\xi\sim\mathcal{N}(\bm 0, \tensor I)$. Hence, 
$\bm\eta\sim\mathcal{N}(\bm 0, \tensor I)$. For general cases where
$\{\xi_i\}_{i=1}^d$ are not i.i.d. Gaussian, $\{\eta_i\}_{i=1}^d$ are not 
necessarily independent. Moreover, $\{\psi_n\}_{n=1}^N$ are not necessarily 
orthogonal to each other with respect to $\rho_{\bm\eta}$. Therefore,
$v_g(\bm\eta)$ may not be a standard gPC expansion of $v(\bm\eta)$. It is
a polynomial equivalent to $u_g(\bm\xi)$ with potentially sparser coefficients.
The idea of using a linear map is also used in sliced inverse 
regression (SIR) \cite{Li91}, active subspace \cite{Russi10, ConstantineDW14},
and basis adaptation \cite{TipG14}, while these methods compute the matrix in
different manners. More importantly, in contrast to these methods, our method does
not truncate the dimension, and $\tensor A$ is a square matrix. We use an 
iterative algorithm to identify $\tensor A$, and the initial guess may not be 
sufficiently accurate. Thus, reducing dimension before the iterations terminate
may lead to less accurate results. The dimension reduction can be integrated 
with the iterative method proposed in this work, e.g., an algorithm for i.i.d.
Gaussian random variables with dimension reduction after iterative rotations was
proposed in \cite{yang2017pdf}. An iterative rotation method preceded with 
SIR-based dimension reduction was proposed in \cite{yang2017sliced}. We refer
interested readers to the respective literatures.

In this work, we use the gradient information of $u$ to identify the
rotation matrix $\tensor A$ based on the framework of active subspace
\cite{Russi10}. We define a ``stiff matrix" $\tensor G$ and compute 
its eigendecomposition as
\begin{equation}\label{eq:grad_mat}
\tensor G \Def
\mathbb{E}\left\{\nabla u(\bx)\cdot\nabla u(\bx)\trans\right\} = \tensor U\tensor \Lambda\tensor
U\trans,
  \quad \tensor U\tensor U\trans = \tensor I.
\end{equation}
Here, $\tensor\Lambda=\text{diag}(\lambda_1,\cdots,\lambda_d)$,  
$\lambda_1\geq\cdots\geq \lambda_d\geq 0$, 
$\nabla u(\bm x) = (\partial u/\partial x_1, \cdots, \partial u/\partial x_d)\trans$
is a column vector, and $\nabla u(\bm\xi)=\nabla u(\bm x)|_{\bm x=\bm\xi}$.
We then set $\tensor A=\tensor U\trans$. The magnitude of $\lambda_i$ 
indicates the importance of the corresponding eigenspace. If only a few
$\lambda_i$ are large and others are very small, we can expect to obtain a
sparse $\tilde{\bm c}$ because only a few $\eta_i$ make major contributions in
expansion $v_g(\bm\eta)$.

In practice, it is difficult to compute the expectation in 
Eq.~\eqref{eq:grad_mat}, which is a high-dimensional integral and $u$ is 
unknown. Thus, we used $u_g$ to approximate $u$ in 
Eq.~\eqref{eq:grad_mat} \cite{LeiYZLB15}, and provided an analytic form to 
approximate $\tensor G$ \cite{YangLBL16} (see Section~\ref{subsec:comparison}). 
Alternatively, in this work, we use another approach to compute the rotation 
matrix as in \cite{ConstantineDW14}. We define
\begin{equation}\label{eq:grad_u}
  \tensor W=\dfrac{1}{\sqrt{M}}[\nabla u(\bm\xi^1),\nabla u(\bm\xi^2),\cdots,\nabla u(\bm\xi^{M})],
\end{equation}
where $M$ is the number of available samples, and we set $M\geq d$ in this work.
Thus, $\tensor W$ is a $d\times M$ matrix. The singular value decomposition
(SVD) or principle component analysis (PCA) of $\tensor W$ yields
\begin{equation}\label{eq:pca}
\tensor W = \tensor U_{W} \tensor \Sigma_{W}\tensor V\trans_{W},
\end{equation}
where $\tensor U_{W}$ is a $d\times d$ orthogonal matrix and 
$\tensor \Sigma_{W}$ is a $d\times M$ matrix, whose diagonal consists of singular
values $\sigma_1\geq\cdots\geq\sigma_d\geq 0$. We set the rotation matrix as
$\tensor A = \tensor U_W\trans$. As such, the rotation projects $\bm\xi$ to the
directions of principle components of $\nabla u$. {\color{red}Of note, we do not use the
information of the orthogonal matrix $\tensor V_W$.}

The connections between Eq.~\eqref{eq:grad_mat} and \eqref{eq:pca} are as
follows: when samples $\bx^q$ are generated based on the distribution of $\bx$
and $M\rightarrow\infty$, we have
$\tensor\Sigma_W\trans\tensor\Sigma_W\rightarrow\tensor \Lambda$, i.e.,
$\sigma_i^2\rightarrow\lambda_i$. Also, the eigenspaces corresponding to
$\sigma_i$ converge to eigenspaces corresponding to $\lambda_i$. In other words,
the rotation identified by Eq.~\eqref{eq:pca} is the same as that from
Eq.~\eqref{eq:grad_mat} asymptotically. These connections are summarized in
the following two lemmas.
\begin{lem}
  \label{lem:1}
Assume that $\nabla u(\bm x)$ exists for any $\bm x\in\mathbb{R}^d$, and 
$\nabla u$ is square-integrable with respect to the measure of $\bx$. If the
matrix $\tensor W$ in Eq.~\eqref{eq:grad_u} is constructed with the samples
$\bx^q$ that are generated based on the distribution of $\bx$, then 
$\Vert\tensor\Sigma_W\tensor\Sigma_W\trans-\tensor\Lambda\Vert_F\rightarrow 0$
as $M\rightarrow\infty$, where $\Vert\cdot\Vert_F$ is the Frobenius norm. 
\end{lem}
\begin{proof}
  Because $\bx^q$ are generated based on the distribution of $\bx$, 
  $\Vert\tensor W\tensor W\trans-\tensor G\Vert_F\rightarrow 0$ as
  $M\rightarrow\infty$. This is straightforward because it uses the Monte
  Carlo method to approximate integrals. Then, because of the following
  eigendecomposition
  \begin{equation}\label{eq:eigen_w}
    \tensor W\tensor W\trans = 
     \tensor U_W\left(\tensor\Sigma_W\tensor\Sigma_W\trans\right)\tensor U_W\trans,
   \end{equation}
  according to Hoffman and Wielandt \cite{hoffman2003variation}, we have
  \[\Vert\tensor\Sigma_W\tensor\Sigma_W\trans-\tensor\Lambda\Vert_F
    =\sqrt{\sum_{i=1}^d (\sigma_i^2-\lambda_i)^2}\leq \Vert\tensor W\tensor
  W\trans-\tensor G\Vert_F \rightarrow 0.\]
\end{proof}

To demonstrate the convergence of eigenspaces, we consider the case when 
$\tensor G$ has at least two distinct eigenvalues. Otherwise, all eigenvector 
directions are equally important, and the rotation can't enhance the sparsity.
Assume without losing generality that
$\lambda_1=\lambda_2=\cdots=\lambda_k>\lambda_{k+1}$, where $0<k<d$. We rewrite
Eq.~\eqref{eq:grad_mat} as
\begin{equation}\label{eq:eigen_g2}
  \tensor G=\tensor U\tensor\Lambda\tensor U\trans=
  \begin{pmatrix}
    \tensor U_1, \tensor U_2
  \end{pmatrix}
  \begin{pmatrix}
    \tensor\Lambda_1 & 0 \\ 0 & \tensor\Lambda_2
  \end{pmatrix}
  \begin{pmatrix}
    \tensor U\trans_1  \\  \tensor U\trans_2
  \end{pmatrix},
\end{equation}
where $\tensor\Lambda_1=\text{diag}(\lambda_1,\lambda_2,\cdots,\lambda_k)$, and
$\tensor U_1$ consists of basis of eigenspace corresponding to $\lambda_1$,
i.e., the first $k$ columns of $\tensor U$.
Accordingly, Eq.~\eqref{eq:eigen_w} can be written as
\begin{equation}\label{eq:eigen_w2}
  \tensor W\tensor W\trans=\tensor U_W(\tensor\Sigma_W\tensor\Sigma_W\trans)\tensor U_W\trans
  =\begin{pmatrix}
    \widetilde{\tensor U}_1, \widetilde{\tensor U}_2
  \end{pmatrix}
  \begin{pmatrix}
    \widetilde{\tensor\Lambda}_1 & 0 \\ 0 & \widetilde{\tensor\Lambda}_2
  \end{pmatrix}
  \begin{pmatrix}
    \widetilde{\tensor U}\trans_1  \\  \widetilde{\tensor U}\trans_2
  \end{pmatrix},
\end{equation}
where $\widetilde{\tensor\Lambda}_1=\text{diag}(\sigma_1^2,\sigma_2^2,\cdots,\sigma_k^2)$,
and $\widetilde{\tensor U}_1$ consists of the first $k$ columns of $\tensor U_W$.
We denote the subspace spanned by columns of $\tensor U_1$ as $\mathbb{H}$, and the
subspace spanned by columns of $\widetilde{\tensor U}_1$ as
$\widetilde{\mathbb{H}}$. We demonstrate that the \emph{principle angles} (or
\emph{canonical angles}) \cite{stewart1990matrix} between $\mathbb{H}$ and
$\widetilde{\mathbb{H}}$ convergence to $0$ as $M\rightarrow\infty$, which 
indicates that $\widetilde{\mathbb{H}}$ converges to $\mathbb{H}$.
\begin{lem}\label{lem:2}
 Let $\tensor U_1$  and $\widetilde{\tensor U}_1$ be defined in
 Eq.~\eqref{eq:eigen_g2} and Eq.~\eqref{eq:eigen_w2}. The 
 \emph{angle matrix} between $\tensor U_1$ and $\widetilde{\tensor U}_1$ is
defined as
\[\tensor\Theta(\tensor U_1,\widetilde{\tensor U}_1)\Def
\arccos\left(\tensor U_1\trans\widetilde{\tensor U}_1\widetilde{\tensor
U}\trans_1\tensor U_1\right)^{-1/2}.\]
 We have 
 $\Vert\sin\tensor\Theta(\tensor U_1,\widetilde{\tensor U}_1)\Vert_F\rightarrow 0$
 as $M\rightarrow \infty$.
\end{lem}
\begin{proof}
  Let $\beta=\displaystyle\min_{1\leq j\leq d-k}|\lambda_k-\lambda_{k+j}|/2>0$. According to Lemma~\ref{lem:1}, there
  exists an integer $M_0>0$, such that for any $M>M_0$, 
  $\displaystyle\max_{1\leq i\leq d}|\sigma_i^2-\lambda_i|<\beta$. Therefore,
  $\displaystyle\min_{1\leq i\leq k, 1\leq j\leq
  d-k}|\lambda_i-\sigma_{k+j}^2|>\beta$ for any $M>M_0$. According to \cite{davis1970rotation}, 
  we have
  \[\Vert\sin\tensor\Theta(\tensor U_1,\widetilde{\tensor U}_1)\Vert_F\leq
  \dfrac{\Vert\left(\tensor G-\tensor W\tensor W\trans\right)\tensor
U_1\Vert_F}{\beta}
\leq \dfrac{\Vert\tensor G-\tensor W\tensor W\trans\Vert_F\Vert\tensor
U_1\Vert_F}{\beta}\rightarrow 0.\]
\end{proof}
The principle angles are singular values of 
$\tensor\Theta(\tensor U_1,\widetilde{\tensor U}_1)$,
and Lemma~\ref{lem:2} indicates that they converge to $0$. Thus,
$\widetilde{\mathbb{H}}$ converges to $\mathbb{H}$. Similarly, the eigenspace
associated with other singular values of $\tensor W$ converges to eigenspace
associated with corresponding eigenvalues of $\tensor G$.

Because $u$ is not known and we assume that samples of $\nabla u$ are not
available, we replace $u$ with $u_g$ in Eq.~\eqref{eq:grad_u} for approximation: 
\begin{equation}
  \tensor W\approx\tensor W_g=\dfrac{1}{\sqrt{M}}[\nabla u_g(\bm\xi^1),\nabla
u_g(\bm\xi^2),\cdots,\nabla u_g(\bm\xi^{M})],
\end{equation}
and the rotation matrix is constructed based on the SVD of $\tensor W_g$:
\begin{equation}\label{eq:pcag}
\tensor W_g = \tensor U_{W_g} \tensor \Sigma_{W_g}\tensor V\trans_{W_g}, \quad
\tensor A = \tensor U_{W_g}\trans.
\end{equation}
Here, $u_g$ can be computed with different methods, and we use standard $\ell_1$
(or re-weighted $\ell_1$) minimization to compute it. On obtaining $\tensor A$,
we define $\bm\eta=\tensor A\bm\xi$ and compute the corresponding input samples
as $\bm\eta^q=\tensor A\bm\xi^q$. After, we construct a new measurement 
matrix $\tensor\Psi(\bm\eta)$ as $(\tensor\Psi(\bm\eta))_{ij}=\psi_j(\bm\eta^i)$.
Next, we solve the $\ell_1$ minimization problem $(P_{1,\varepsilon})$ to obtain
$\tilde{\bm c}$. If a few singular values of $\tensor W$ are much larger than
others, we can expect that the representation of $u$ with respect to
$\bm\eta$ mainly replies on the eigenspace associated with these singular
values, i.e., $\eta_i$ corresponding to large $\sigma_i$. Therefore, we can
obtain a sparser representation. On the other hand, if the differences
between $\sigma_i$ are small, the rotation will not enhance the sparsity, and
this method may not result in a more accurate representation of uncertainty.


\subsection{Iterative method}
\label{subsec:iteration}

We generalize the ideas in Section \ref{subsec:basic} and formulate it as the
following modified $\ell_1$ minimization problem: 
\begin{equation}\label{eq:rot_l1}
  (P_{1,\varepsilon}^R)\qquad\qquad  \argmin_{\tc,~ \tensor A} \Vert \tc\Vert_1,
  \quad \text{subject to}\quad \Vert\tensor \Psi\tc-\bm u\Vert_2\leq\varepsilon,
  \quad \tensor A\tensor A\trans=\tensor I,
\end{equation}
where $\tensor\Psi_{ij}=\Psi_j(\tensor A\bm\xi^i)$. Here, the orthogonality
requirement on $\tensor A$ can be released to obtain a better result, i.e.,
making $\Vert\epsilon(\bx)\Vert_2$ smaller. However, even solving optimization 
problem $(P_{1,\varepsilon}^R)$ can be difficult when dimension is high, i.e., 
when $d$ is large, because the degree of freedom is $d+d(d-1)=d^2$. Similarly,
the degree of freedom is $d+d^2=d(d+1)$ when the constraint 
$\tensor A\tensor A\trans=\tensor I$ is not imposed.
Therefore, we propose to use an alternating direction method that combines
$\ell_1$ minimization solver and rotation matrix based on the gradients to
approximate the solution of $(P_{1,\varepsilon}^R)$.

We start with an initial guess $u_g(\bx)=\sum_{n=1}^N\tilde c^{(0)}\psi_n(\bx)$,
which is obtained by standard $\ell_1$ minimization or re-weighted $\ell_1$
minimization (Algorithm \ref{algo:cs1}). Here, the superscript
$\cdot^{(l)}$ stands for the $l$-th iteration, e.g., $\tc^{(l)}$ are the gPC
coefficients of the $l$-th iteration. We also set $\tensor A^{(0)}=\tensor I,
\bm\eta^{(0)}=\bm\xi, v_g^{(0)}(\bm\eta^{(0)})=u_g(\bm\xi)$.
In the $l$-th iteration ($l\geq 1$), given $v_g^{(l-1)}$ and input samples
$\{(\bm\eta^{(l-1)})^q\}_{q=1}^M$, we first collect the gradient of $v_g^{(l-1)}$:
\begin{equation}\label{eq:grad_iter}
  \tensor W_g^{(l-1)}=\dfrac{1}{\sqrt{M}}\left [\nabla_{_\xi} v_g^{(l-1)}\left((\bm\eta^{(l-1)})^1\right),
    \cdots, \nabla_{_\xi} v_g^{(l-1)}\left((\bm\eta^{(l-1)})^{M}\right)\right ],
\end{equation}
where $\nabla_{_\xi}\cdot=(\partial\cdot/\partial \xi_1, \partial\cdot/\partial \xi_2,
\cdots, \partial\cdot/\partial \xi_d)\trans$. Next, we compute the SVD of 
$\tensor W_g^{(l-1)}$:
\begin{equation}\label{eq:svd_iter}
\tensor W_g^{(l-1)} = \tensor U_{W_g}^{(l-1)} \tensor
\Sigma_{W_g}^{(l-1)}\left(\tensor V_{W_g}^{(l-1)}\right)\trans, 
\end{equation}
and set $\tensor A^{(l)}=\left(\tensor U_{W_g}^{(l-1)}\right)\trans$. Now, we can 
define a new set of random variables as $\bm\eta^{(l)}=\tensor A^{(l)}\bm\xi$ and
compute their samples accordingly: $(\bm\eta^{(l)})^q=\tensor A^{(l)}\bm\xi^q$. 
We then construct a new measurement matrix $\tensor\Psi^{(l)}$ as 
$\Psi^{(l)}_{ij}=\psi_j((\bm\eta^{(l)})^i)$, and solve the $\ell_1$ minimization
problem $(P_{1,\varepsilon})$ to obtain $\bm c^{(l)}$. Of note, the gradient in
Eq.~\eqref{eq:grad_iter} is computed with respect to $\bm\xi$. Specifically, we
use the chain rule in the computing:
\begin{equation}\label{eq:chain}
  \begin{split}
 &  \nabla_{_\xi} v_g^{(l)}\left((\bm\eta^{(l)})^q\right)=
  \nabla_{_\xi} v_g^{(l)}\left(\tensor A^{(l)}\bm\xi^q\right) 
  = (\tensor A^{(l)})\trans \nabla v_g^{(l)}(\bm x) \bigg |_{\bm x=\tensor A^{(l)}\xi^q} \\
  = &  (\tensor A^{(l)})\trans\nabla \sum_{n=1}^N \tilde c_n^{(l)}\psi_n(\bm x)\bigg |_{\bm x=\tensor A^{(l)}\xi^q}
  = (\tensor A^{(l)})\trans \sum_{n=1}^N \tilde c_n^{(l)}\nabla \psi_n(\bm x)\bigg |_{\bm x=\tensor A^{(l)}\xi^q} .
  \end{split}
\end{equation}
Therefore, we only need to evaluate $\nabla\psi_n$ at $(\bm\eta^{(l)})^q$. This
is straightforward because we construct $\psi_n$ using the tensor product of
univariate polynomials (Eq.~\eqref{eq:tensor}) and derivatives of widely used
orthogonal polynomials in UQ study, e.g., Hermite, Laguerre, Legendre and
Chebyshev, are known. We summarize the entire procedure in
Algorithm \ref{algo:cs_rot1}.
\begin{algorithm}[t]
 \caption{Alternating direction method of solving $(P^R_{1,\epsilon})$}
\label{algo:cs_rot1}
\begin{algorithmic}[1]
  \STATE Generate input samples $\{\bx^q\}_{q=1}^M$ based on the distribution
     of $\bx$. 
\STATE Generate corresponding output samples $\{u^q=u(\bx^q)\}_{q=1}^M$ by
    solving the complete model, e.g., running simulations, solvers, etc.
\STATE Select gPC basis functions $\{\psi_n\}_{n=1}^N$ associated with $\bx$ and
then generate the measurement matrix $\tensor\Psi$ by setting
$\Psi_{ij}=\psi_j(\bx^i)$.
\STATE Solve the optimization problem $(P_{1,\epsilon})$:
\[\arg \min_{\hc}\Vert\hc\Vert_1, ~\text{subject to}~ \Vert\tensor\Psi\hat{\bm c}-\bm u\Vert_2\leq\epsilon.\]
\STATE Set counter $l=0$, $\bm\eta^{(0)}=\bm\xi, \tilde{\bm c}^{(0)}=\hat{\bm
c}$, $\tensor A^{(0)}=\tensor I$, and $v_g^{(0)}=\sum_{n=1}^N\tilde
c^{(0)}\psi_n(\bx)$.
\STATE $l=l+1$. Construct $\tensor W_g^{(l-1)}$ in Eq.~\eqref{eq:grad_iter} by 
computing $\nabla_{_\xi}v^{(l-1)}_g\left((\bm\eta^{(l-1)})^q\right)$ according 
to Eq.~\eqref{eq:chain}. Then compute SVD of $\tensor W_g^{(l)}$:
$\tensor W_g^{(l)}= \tensor U^{(l)}_{W_g}\tensor\Sigma^{(l)}_{W_g}\left(\tensor
  V^{(l)}_{W_g}\right)\trans.$
\STATE Set $\tensor A^{(l)}=\left(\tensor U^{(l)}_{W_g}\right)\trans$ and 
  $\bm\eta^{(l)}=\tensor A^{(l)}\bm\xi$. Then compute samples 
  $(\bm\eta^{(l)})^q=\tensor A^{(l)}\bx^q, q=1,2,\cdots,M$.
Construct the new measurement matrix $\tensor\Psi^{(l)}$ with
$\Psi^{(l)}_{ij}=\psi_j\left((\bm\eta^{(l)})^i\right)$. 
\STATE Solve the optimization problem $(P_{1,\epsilon^{(l)}})$:
\[\arg \min_{\hat{\bm c}}\Vert\hat{\bm c}\Vert_1, \quad\text{subject to}~
\Vert\tensor\Psi^{(l)}\hat{\bm c}-\bm u\Vert_2\leq\epsilon^{(l)},\]
then set $\tilde{\bm c}^{(l)}=\hat{\bm c}$ and 
$v_g^{(l)}(\bm\eta^{(l)})=\sum_{n=1}^N\tilde c^{(l)}\psi_n(\bm\eta^{(l)})$. 
\STATE If $\Vert\tensor A^{(l)}-\tensor A^{(l-1)}\Vert_2<\theta$, where the 
threshold $\theta$ is a positive real number, then stop. Otherwise, go to Step 6.
\STATE Construct gPC expansion as 
$u(\bx)\approx u_g(\bm\xi)=v_g(\bm\eta^{(l)})
=\sum_{n=1}^N \tilde{c}^{(l)}_n\psi_n(\tensor A^{(l)}\bx)$.
\end{algorithmic}
\end{algorithm}
The stopping criterion we use is the difference between to successive
rotation matrix $\tensor A^{(l)}$ and $\tensor A^{(l+1)}$. The threshold 
$\theta$ can be set as a fraction of the dimension, e.g., $0.2d$. Alternatively,
it can be set as $l\leq l_{\max}$, where $l_{\max}$ is can taken as $2$ or $3$
empirically for practical problems. 
Because it is not necessary to use the optimal rotation to improve the sparsity. 
In addition, the sparsity structure is problem-dependent (see examples in
Section~\ref{sec:numeric}), and, for many practical problems, the improvement is
not made significant by using more iterations. Notably, we update $\bm c^{(l)}$
and $\tensor A^{(l)}$ separately in each iteration, which is the spirit of the
alternating direction method. 


\subsection{A comparison with the Gaussian case}
\label{subsec:comparison}

An alternating direction algorithm approximating the solution of
$(P_{1,\epsilon}^R)$ when $\xi_i$ are i.i.d. Gaussian random variables,
i.e., $\psi_n$ are Hermite polynomials, was proposed in \cite{YangLBL16},
and we present it in Algorithm~\ref{algo:cs_rot3}.
\begin{algorithm}[t]
  \caption{Alternating direction algorithm for $(P_{1,\varepsilon}^R)$ when $\xi_i$ are i.i.d. Gaussian.}
\label{algo:cs_rot3}
\begin{algorithmic}[1]
  \STATE Generate input samples $\{\bx^q\}_{q=1}^M$ based on the distribution
of $\bx$. 
\STATE Generate output samples $\{u^q=u(\bx^q)\}_{q=1}^M$ by solving the complete model.
\STATE Set $\{\psi_n\}_{n=1}^N$ as multivariate Hermite polynomials, 
then generate the measurement matrix $\tensor\Psi$ as $\Psi_{ij}=\psi_j(\bx^i)$.
\STATE Solve the optimization problem $(P_{1,\epsilon})$:
\[\arg \min_{\hat{\bm c}}\Vert\hat{\bm c}\Vert_1, ~
\text{subject to}~ \Vert\tensor\Psi\hat{\bm c}-\bm u\Vert_2\leq\epsilon.\]
\STATE Set counter $l=0$, $\eta^{(0)}=\bx$, $\tilde{\bm c}^{(0)}=\hat{\bm c}$,
and compute $\tensor K_{ij}, i,j=1,2,\cdots,N$ as
\begin{equation}\label{eq:kernel}
  (K_{ij})_{km} = \mathbb{E}\left\{\dfrac{\partial\psi_k}{\partial x_i}(\bm\xi)
  \cdot \dfrac{\partial\psi_m}{\partial x_j}(\bm \xi)\right\}. 
\end{equation}                
\STATE $l=l+1$. Construct ${\tensor G}^{(l-1)}$ as 
$G^{(l)}_{ij}=(\tilde{\bm c}^{(l-1)})\trans \tensor K_{ij} \tilde{\bm c}^{(l-1)},
i,j=1,2,\cdots, d$.
Then, compute eigendecomposition of $\tensor G^{(l-1)}$:
$\tensor G^{(l-1)}= 
\tensor U^{(l-1)}\tensor\Lambda^{(l-1)}\left(\tensor U^{(l-1)}\right)\trans.$
\STATE Set $\bm\eta^{(l)}=\left(\tensor U^{(l)}\right)\trans\bm\eta^{(l-1)}$,
then compute samples $(\bm\eta^{(l)})^q=\left(\tensor
U^{(l-1)}\right)\trans(\bm\eta^{(l-1)})^q, q=1,2,\cdots,M$.
Construct the new measurement matrix $\tensor\Psi^{(l)}$ with
$\Psi^{(l)}_{ij}=\psi_j\left((\bm\eta^{(l)})^i\right)$. 
\STATE Solve the optimization problem $(P_{1,\epsilon^{(l)}})$:
\[\arg \min_{\hat{\bm c}}\Vert\hat{\bm c}\Vert_1, \quad\text{subject to}~
\Vert\tensor\Psi^{(l)}\hat{\bm c}-\bm u\Vert_2\leq\epsilon^{(l)},\]
and set $\tilde{\bm c}^{(l)}=\hat{\bm c}$. 
\STATE If $|\Vert\widetilde{\tensor U}^{(l)}\Vert_1-d|<\theta$, where the 
threshold $\theta$ is a positive real number, then stop. Otherwise, go to Step 6.
\STATE Set 
\[\tensor A^{(l)}=\left(\tensor U^{(0)}\tensor U^{(2)}\cdots\tensor
U^{(l-1)}\right)\trans,\]
and construct gPC expansion as $u(\bx)\approx u_g(\bm\xi)=v_g(\bm\eta^{(l)})=\sum_{n=1}^N \tilde{c}^{(l)}_n\psi_n(\tensor A^{(l)}\bx)$.
\end{algorithmic}
\end{algorithm}
In this section, we illustrate the connections and differences between 
Algorithms~\ref{algo:cs_rot1} and \ref{algo:cs_rot3}. 

These algorithms use different approaches to obtain rotation matrix that maps
$\bx$ to $\bm\eta$. Both of them use gradient information. However, 
Algorithm \ref{algo:cs_rot1} stems from SVD of $\tensor W$ in 
Eq.~\eqref{eq:pca}, while Algorithm~\ref{algo:cs_rot3} roots on 
eigendecomposition of the variance of the gradients in Eq.~\eqref{eq:grad_mat}
as in \cite{LeiYZLB15}.
In particular, Steps 5-8 of Algorithm~\ref{algo:cs_rot3} is a more sophisticated
version of Algorithm 1\ in \cite{LeiYZLB15} (because $\tensor G$ is approximated
with an analytic form), which can be considered as a special case of 
Algrothm~\ref{algo:cs_rot3} without iteration. When $u$ is known, they are 
equivalent asymptotically as proven in Lemmas~\eqref{lem:1} and ~\eqref{lem:2}.
Although, in practice, $u$ is unknown and $M$ is limited, these lemmas still 
provide intuitive understanding of the connection between 
Algorithm~\ref{algo:cs_rot1} and Algorithm~\ref{algo:cs_rot3}: they both use the
distinct importance of subspaces identified by $\nabla u$ to enhance the 
sparsity via projecting $\bx$ to these subspaces.

For the Gaussian case, i.e., Algorithm~\ref{algo:cs_rot3}, we use
Eq.~\eqref{eq:grad_mat} and approximate $\tensor G$ as
\begin{equation}\label{eq:stiff}
  \tensor G\approx \mathbb{E}\left\{\nabla u_g(\bm\xi)\cdot\nabla
  u_g(\bm\xi)\trans\right\} =
  \mathbb{E}\left\{\nabla\left(\sum_{n=1}^Nc_n\psi_n(\bx)\right)
    \cdot \nabla\left(\sum_{n'=1}^Nc_{n'}\psi_{n'}(\bx)\right)\trans\right\}.
\end{equation}
For simplicity, we denote $\frac{\partial h(\bm x)}{\partial x_i}$ as
$\partial_i h(\bm x)$ for any differentiable function $h$. Then,
Eq.~\eqref{eq:stiff} implies that
\begin{equation}\label{eq:grad}
\begin{aligned}
G_{ij} & \approx\mexp{
\partial_i\left(\sum_{n=1}^N c_n\psi_n(\bm\xi)\right)\cdot
\partial_j\left(\sum_{n'=1}^Nc_{n'}\psi_{n'}(\bm\xi)\right)}  
 = \mexp{\left(\sum_{n=1}^Nc_n \partial_i\psi_n(\bm\xi)\right) \cdot
    \left(\sum_{n'=1}^Nc_{n'}\partial_j\psi_{n'}(\bm\xi)\right)} \\ 
& = \sum_{n=1}^N\sum_{n'=1}^Nc_nc_{n'}
    \mexp{\partial_i\psi_n(\bm\xi)\cdot \partial_j\psi_{n'}(\bm\xi)} 
 = \bm c^T \tensor K_{ij} \bm c.
\end{aligned}
\end{equation}
In each iteration of Algorithm~\ref{algo:cs_rot3}, 
\begin{equation}
    G^{(l)}_{ij}  =
\mathbb{E}\left\{
\partial_i\left(\sum_{n=1}^N \tilde c^{(l-1)}_n\psi_n(\bm\eta^{(l-1)})\right)\cdot
\partial_j\left(\sum_{n'=1}^N\tilde
c^{(l-1)}_{n'}\psi_{n'}(\bm\eta^{(l-1)})\right)\right\} 
  = (\tilde{\bm c}^{(l-1)})\trans\tensor K_{ij}^{(l-1)}\tilde{\bm c}^{(l-1)},
\end{equation}
and
\begin{equation}
 (K_{ij}^{(l-1)})_{km}  =
 \mathbb{E}\left\{\partial_i\psi_k(\bm\eta^{(l-1)})\cdot\partial_j\psi_m(\bm\eta^{(l-1)})\right\}.
\end{equation}
In this instance, we do not need to update
$\tensor K_{ij}$ because $\eta^{(l)}_i$ are i.i.d. Gaussian for each $l$, and 
\[\mathbb{E}\left\{\partial_i\psi_k(\bm\eta^{(l)})\cdot\partial_j\psi_m(\bm\eta^{(l)})\right\}\equiv
  \dfrac{1}{(2\pi)^{d/2}}\int_{\mathbb{R}^d} \partial_i\psi_k(\bm x)\cdot\partial_j\psi_m(\bm
  x)\exp\left(-\dfrac{\Vert \bm x\Vert_2^2}{2}\right)\dif\bm x\]
are fixed so it can be precomputed (see \cite{YangLBL16}). The 
Algorithms~\ref{algo:cs_rot3}'s design employs this advantage. However, 
for general cases, this good property may not exist, and we need to update
$\tensor K_{ij}$ in each step if we use Eq.~\eqref{eq:grad}.
Specifically, in the $l$-th iteration of Algorithm~\ref{algo:cs_rot1}, 
\begin{equation}
  \left | \det\left(\dfrac{D\bm\eta^{(l-1)}}{D\bm\xi}\right)\right| = \left |
  \det\left((\tensor A^{(l-1)})^{-1}\right)\right| = 1,
\end{equation}
and the PDF of $\bm\eta^{(l-1)}$ is
\begin{equation}
  \rho_{\bm\eta^{(l-1)}}(\bm x) = \rho_{\bx}((\tensor A^{(l-1)})^{-1} \bm x)
  \left |\det\left(\dfrac{D\bm\bx}{D\bm\eta^{(l-1)}}\right)\right| =
  \rho_{\bx}\left((\tensor A^{(l-1)})^{-1} \bm x\right ).
\end{equation}
Thus, if we want to approximate $\tensor G$ as in Algorithm~\ref{algo:cs_rot3},
we need to compute
\begin{equation}\label{eq:kernel_general}
\begin{aligned}
  \mexp{\partial_i\psi_n(\bm\eta^{(l-1)})\cdot\partial_j\psi_{n'}(\bm\eta^{(l-1)})}
  & = \int_{\Omega_{\bm\eta^{(l-1)}}}\partial_i\psi_n(\bm x)\cdot
  \partial_j\psi_{n'}(\bm x)\rho_{\bm\eta^{(l-1)}}(\bm x)\dif\bm x \\
& = \int_{\Omega_{\bm\eta^{(l-1)}}}\partial_i\psi_n(\bm x)\cdot
  \partial_j\psi_{n'}(\bm x)\rho_{\bm\xi}((\tensor A^{(l-1)})^{-1}\bm x)\dif\bm x,
\end{aligned}
\end{equation}
where $\Omega_{\bm\eta^{(l-1)}}$ is the domain of multivariate random variable
$\bm\eta^{(l-1)}$. If $\bm\xi\sim\mathcal{N}(\bm 0, \tensor I)$, then
$\Omega_{\bm\eta^{(l-1)}}=\Omega_{\bx}$, and
\begin{equation}\label{eq:gauss_pdf_rot}
  \rho_{\bm\eta^{(l-1)}}(\bm x)=\rho_{\bx}((\tensor A^{(l-1)})^{-1}\bm
  x)=\rho_{\bx}(\bm x).
\end{equation}
Hence, we only need to compute $\tensor K_{ij}$ once.
In other cases, updating $\tensor K_{ij}$ in each iteration can be costly
because of the high-dimensional integral in Eq.~\eqref{eq:kernel_general}.

A possible solution is to compute $\tensor K_{ij}$ based on $\rho_{\bm\xi}$
first. Then, in each iteration, after obtaining 
$v^{(l)}_g(\bm\eta^{(l)})=\sum_{i=1}^N\tilde c_n^{(l)}\psi_n(\bm\eta^{(l)})$, we
compute the corresponding $u^{(l)}_g(\bm\xi)=\sum_{i=1}^N
c_n^{(l)}\psi_n(\bm\xi)$ through algebraic computing or by accurate numerical
integral:
\[c^{(l)}_n = \int_{\Omega_{\bx}} v_g^{(l)}(\tensor A^{(l)}\bm x)\psi_n(\bm x)\rho_{\bm\xi}(\bm x)\dif \bm x 
    =\sum_{q=1}^{N_q} v_g^{(l)}(\tensor A^{(l)}\bm x^q)\psi_n(\bm x^q)w^q,\]
where $\bm x^q$ and $w^q$ are quadrature points and weights with respect to
$\rho_{\bx}(\bm x)$. As such, $G_{ij}^{(l)}$ can be approximated by $(\bm
c^{l})^T\tensor K_{ij}\bm c^{(l)}$. However, the additional computation cost for
converting $v_g^{(l)}$ to $u_g^{(l)}$ makes this algorithm less efficient.
Therefore, we choose to use the SVD of $\tensor W_g$ in
Algorithm~\ref{algo:cs_rot1}.

In addition to the different approaches for computing the rotation matrix,
another difference is that in each iteration, Algorithm~\ref{algo:cs_rot1}
directly identifies rotation matrix $\tensor A^{(l)}$ that maps $\bx$ to
$\bm\eta^{(l)}$, while Algorithm~\ref{algo:cs_rot3} seeks for $\tensor U^{(l)}$,
which is a ``correction" of the existing rotation, that maps $\bm\eta^{(l-1)}$
to $\bm\eta^{(l)}$. In each iteration of Algorithm~\ref{algo:cs_rot3}, 
we do not need to compute $\tensor A^{(l)}$, and we compute it after the
iterations terminate. Thus, Algorithm~\ref{algo:cs_rot3} does not rely on the
chain rule. We also can design Algorithm~\ref{algo:cs_rot1} in this manner,
i.e., each iteration maps $\bm\eta^{(l-1)}$ to $\bm\eta^{(l)}$. Instead, we
use the current design to explicitly fit the description of the alternating 
direction method, i.e., in each iteration, we explicitly identify 
$\tilde{\bm c}^{(l)}$ and $\tensor A^{(l)}$ separately.

Finally, according to our numerical tests (not shown), there is no
significant difference between the accuracy of Algorithms~\ref{algo:cs_rot1} and
\ref{algo:cs_rot3} for Hermite polynomial expansions. This indicates that the
general framework in Algorithm \ref{algo:cs_rot1} also is efficient for $u$
relying on i.i.d. Gaussian random variables.


\subsection{Compromising the property of $\tensor\Psi$}
\label{subsec:compromise}

We already have discussed the possible enhancement of the sparsity in 
$\tilde{\bm c}$ by introducing the rotation. However, as we indicated at the 
beginning of this section, $\{\psi_n\}_{n=1}^N$ are not necessarily orthonormal
to each other with respect to $\rho_{\bm\eta}$. The property of matrix
$\tensor\Psi^{(l)}$ may become less favorable for the $\ell_1$ minimization. The
most straightforward conclusion we can obtain is that the \emph{mutual 
coherence} of $\tensor\Psi^{(l)}$ ($l\geq 1$) can be larger than that of
$\tensor\Psi^{(0)}$. Here, the mutual coherence \cite{BrucksteinDE09}, defined as
\begin{equation}\label{eq:mutual}
\mu(\tensor\Psi) \Def \max_{1\leq j,k\leq N,j\neq k}
\dfrac{|\bm\Psi_j^T\bm\Psi_k|}{\Vert\bm\Psi_j\Vert_2\cdot\Vert\bm\Psi_k\Vert_2},
\end{equation}
where $\bm\Psi_j$ and $\bm\Psi_k$ are columns of $\tensor\Psi$, is a more
tractable property of the measurement matrix than the RIP. Generally, a
measurement matrix with smaller mutual coherence is better able to recover a 
sparse solution with the compressive sensing method. When $\psi_n$ are
Hermite polynomials ($\{\xi_i\}_{i=1}^d$ are i.i.d. Gaussian), 
$\mu(\tensor\Psi)$ is conserved (statistically) because $\bm\eta^{(l)}$
are still i.i.d. Gaussian. This adds to the specialty of the Hermite polynomial
expansions. For other polynomials, the aforementioned rotational method can be
less efficient in some cases due to the increase of $\mu(\tensor\Psi)$ (detailed
in Section~\ref{sec:numeric}). We note that in the compressive sensing theory, 
the number of samples needed for an accurate computing of $\bm c$ is related to
both the sparsity of $\bm c$ and the property of $\tensor\Psi$. Therefore, for
general cases, although we can increase the sparsity of $\bm c$, if $\tensor\Psi$
becomes worse (i.e., it loses RIP or $\mu$ increases), out method's efficiency 
can be affected.

\section{Numerical Examples}
\label{sec:numeric}
In this section, we revisit the five numerical examples in \cite{YangLBL16} with
different types of random variables in the systems and various polynomial 
expansions to approximate the solution. By testing the same examples
(with different types of random variables), we demonstrate the efficiency of 
the proposed general framework and we can compare the performance with the
special case (Hermite polynomial expansion) in our previous study. Specifically, 
the random variables considered in this section are uniform random variables 
$\mathcal{U}[-1,1]^d$ (associated with Legendre polynomials) and Chebyshev random 
variables with PDF $\rho(\bm x)=\left(\dfrac{1}{\pi\sqrt{1-x^2}}\right)^d$ 
(associated with Chebyshev polynomials of the first kind, which is denoted as
Chebyshev polynomial for simplicity). 
The accuracies of different methods
are measured by the relative $L_2$ error: $(\Vert u - u_g\Vert_2)/\Vert u\Vert_2$, 
where $u_g$ is the Legendre polynomial expansion or Chebyshev polynomial 
expansion of $u$. The integral
\begin{equation}
\Vert u(\bx)\Vert_2 = \left(\int_{\mr^d} u(\bx)^2\rho(\bx)\dif\bx \right)^{1/2}
\end{equation}
(and $\Vert u-u_g\Vert_2$) is approximated with a high-level sparse grids method,
based on one-dimensional Gaussian quadrature and the Smolyak 
structure \cite{Smolyak63}. The term ``level'' $p$ means that the algebraic 
accuracy of the sparse grids method is $2p-1$. We use $P$ to denote the 
truncation order, which implies that polynomials up to order $P$ are
included in expansion $u_g$. Hence, the number of unknowns can be computed as
$N=\bigl(\begin{smallmatrix}P+d\\d\end{smallmatrix}\bigr)$.

The relative errors we present in this section are obtained from $50$ 
independent replicates for each sample size $M$. For example, we generate $50$
independent sets of input samples $\bx^q,q=1,2,\cdots,M$, compute $50$ 
different relative errors, and report the average of these error samples.
To investigate the effectiveness of increasing the output samples, we set 
the $x$-axis in our figures as the ratio $M/N$, which is the fraction of 
available data with respect to number of unknowns. We use MATLAB package 
\texttt{SPGL1} \cite{BergF08,spgl1} to solve $(P_{1,\epsilon})$. If not otherwise
indicated, we use Algorithm \ref{algo:cs_rot1}, and results are obtained with 
$l_{\max}=3$ iterations. We use three iterations in the re-weighted $\ell_1$
minimization to solve $(P_{1,\epsilon})$.

\subsection{Ridge function}
\label{subsec:ex1}
Consider the following ridge function:
\begin{equation}\label{eq:ex1}
u(\bx) = \sum_{i=1}^d \xi_i + 0.25\left(\sum_{i=1}^d \xi_i\right)^2
       + 0.025\left(\sum_{i=1}^d \xi_i\right)^3,
\end{equation}
where all $\xi_i$ are equally important. In this case, adaptive methods that 
build the surrogate model hierarchically based on the importance of $\xi_i$ 
(e.g., \cite{MaZ10,YangCLK12, ZhangYOKD15}) may not be efficient. 
A simple rotation matrix for this example has the form
\begin{equation}\label{eq:ex_rot}
\tensor A = 
\begin{pmatrix}
d^{-1/2} & d^{-1/2} & \cdots & d^{-1/2} \\
 &   &  & \\
 & \multicolumn{2}{c}{\tilde{\tensor{A}}}  & \\
 &   &  & 
\end{pmatrix},
\end{equation}
where $\tilde{\tensor A}$ is a $(d-1)\times d$ matrix chosen to ensure that
$\tensor A$ is orthonormal. Given this choice for $\tensor A$, 
$\eta_1=(\sum_{i=1}^d\xi_i)/d^{1/2}$ and $u$ has a very simple representation:
\[u(\bx)=v(\bm\eta)=d^{1/2}\eta_1+0.25d\eta_1^2+0.025d^{3/2}\eta_1^3.\]
Therefore, as we keep the set of the basis functions unchanged, all of the
polynomials not related to $\eta_1$ make no contribution to the expansion, which
implies that we obtain a very sparse representation of $u$. Becasue the 
optimal structure is not known \textit{a priori}, the standard 
compressive sensing cannot take advantage of it.

\begin{figure}[h]
\centering
\includegraphics[width=0.45\textwidth]{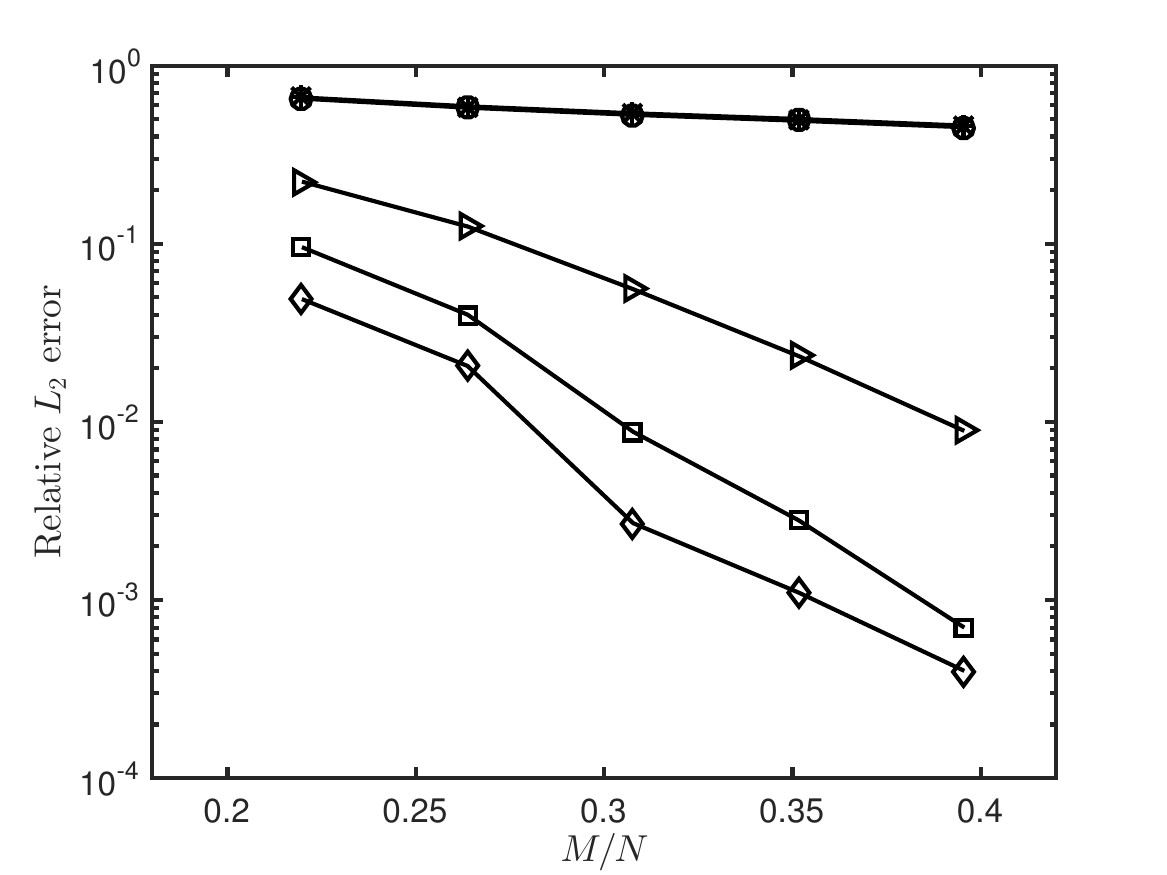}\quad
\includegraphics[width=0.45\textwidth]{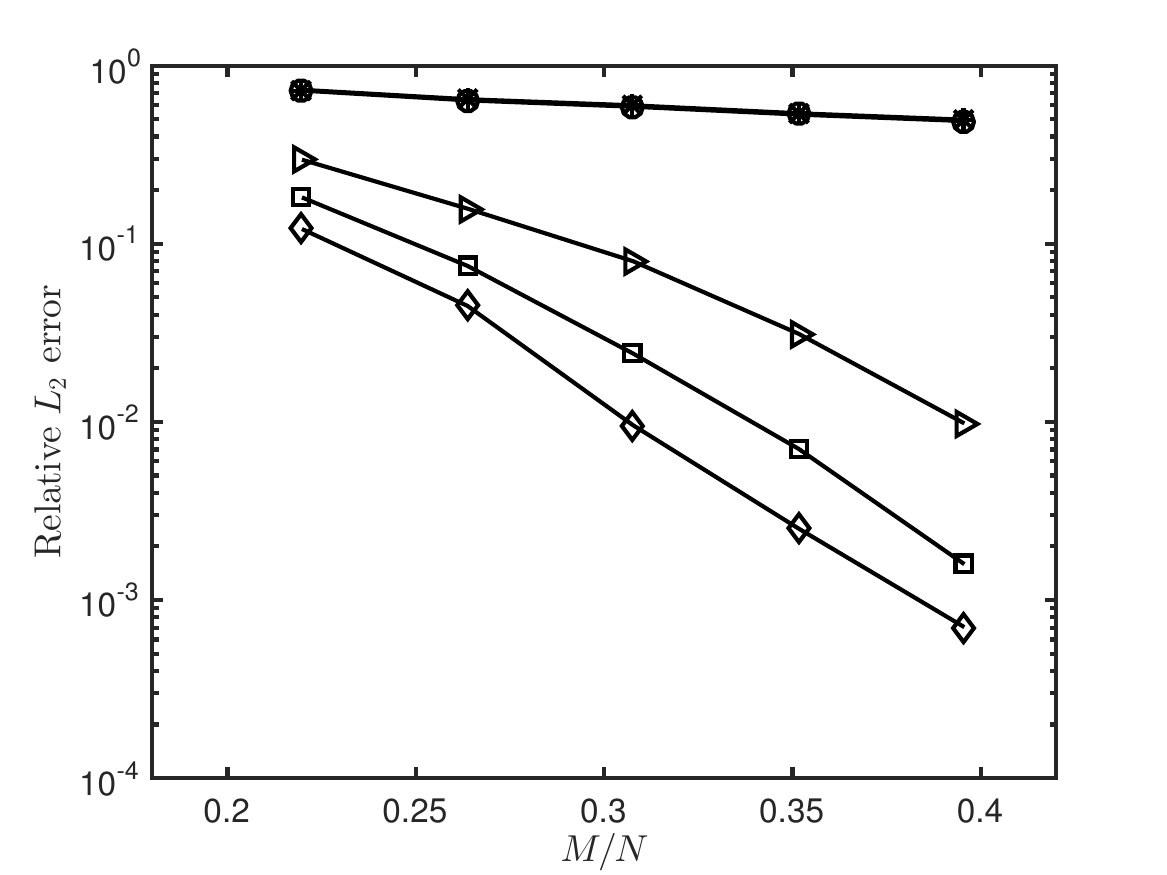}
\caption{Results for the ridge function. Left: Legendre polynomial expansion
  (when $\xi_i$ are i.i.d. uniform random variables); Right: Chebyshev polynomial expansion
  (when $\xi_i$ are i.i.d. Chebyshev random variables). ``$\circ$": standard $\ell_1$,
  ``$\ast$": re-weighted $\ell_1$, ``$\triangleright$": $\ell_1$ with three
  rotations, ``$\square$": $\ell_1$ with six rotations, ``$\diamond$":
  $\ell_1$ with nine rotations.}
\label{fig:ex1_rmse}
\end{figure}

In this test, we set $d=12$ (hence, $N=455$ for $P=3$) and demonstrate the 
effectiveness of our new method. The integrals for calculating the $L_2$ error 
are computed by a level $4$ sparse grids method. Therefore, they are exact. The 
relative errors are presented in Figure~\ref{fig:ex1_rmse} for Legendre polynomial expansion (assuming $\xi_i$
are i.i.d. uniform random variables) and Chebyshev polynomial expansion (assuming
$\xi_i$ are i.i.d. Chebyshev random variables). Clearly, the standard $\ell_1$
minimization is not effective as the relative error is close to $50\%$ even 
when $M/N$ approaches $0.4$. Also, the re-weighted $\ell_1$ does not help in this
case. However, our new iterative rotation demonstrates much better accuracy, 
especially when $M$ is large. We notice that the accuracy increases as more 
iterations are included. Moreover, the improvement from six iterations to nine
iterations is less significant as that from three iterations to six iterations
because the rotation detected by the algorithm approaches the optimal one.

\begin{figure}[h]
\centering
\subfigure[Legendre $|c_n|$]
{\includegraphics[width=0.45\textwidth]{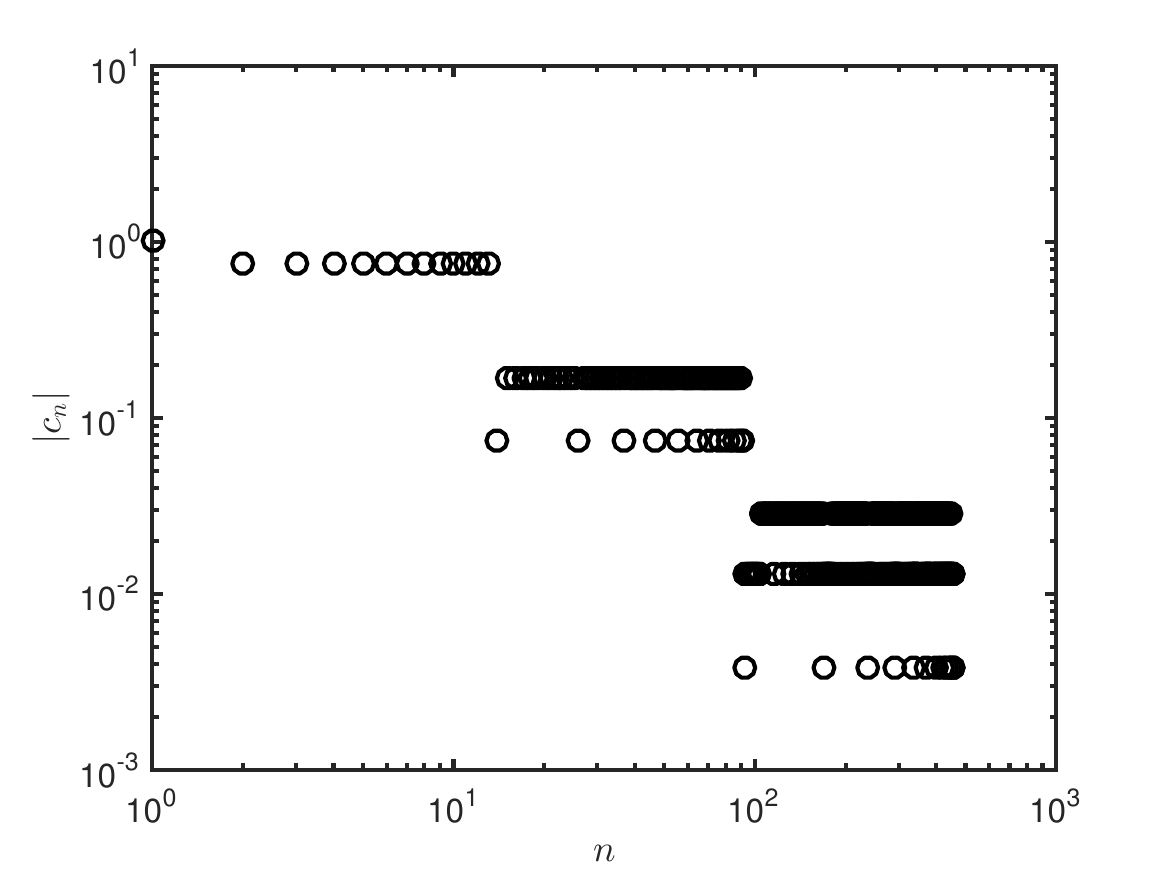}}\quad
\subfigure[Legendre $|\tilde c_n|$]
{\includegraphics[width=0.45\textwidth]{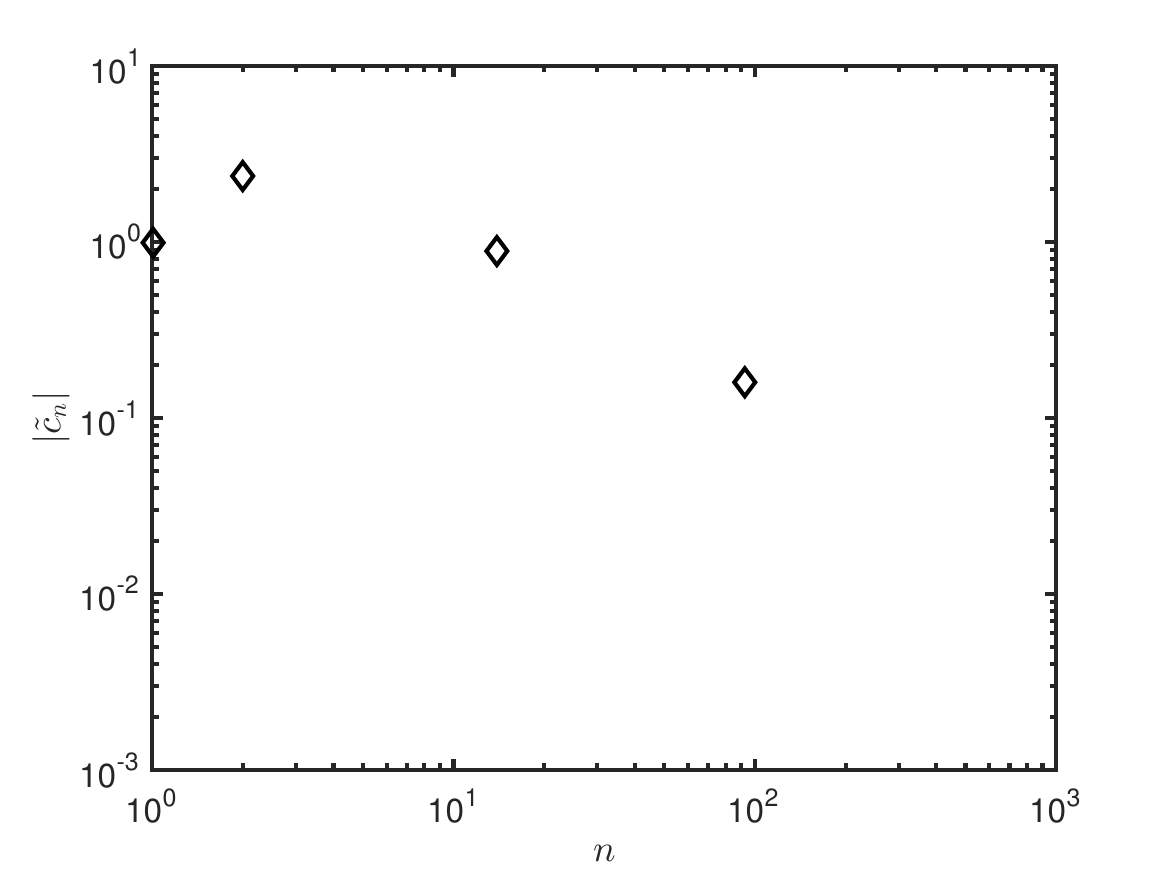}} \\
\subfigure[Chebyshev $|c_n|$]
{\includegraphics[width=0.45\textwidth]{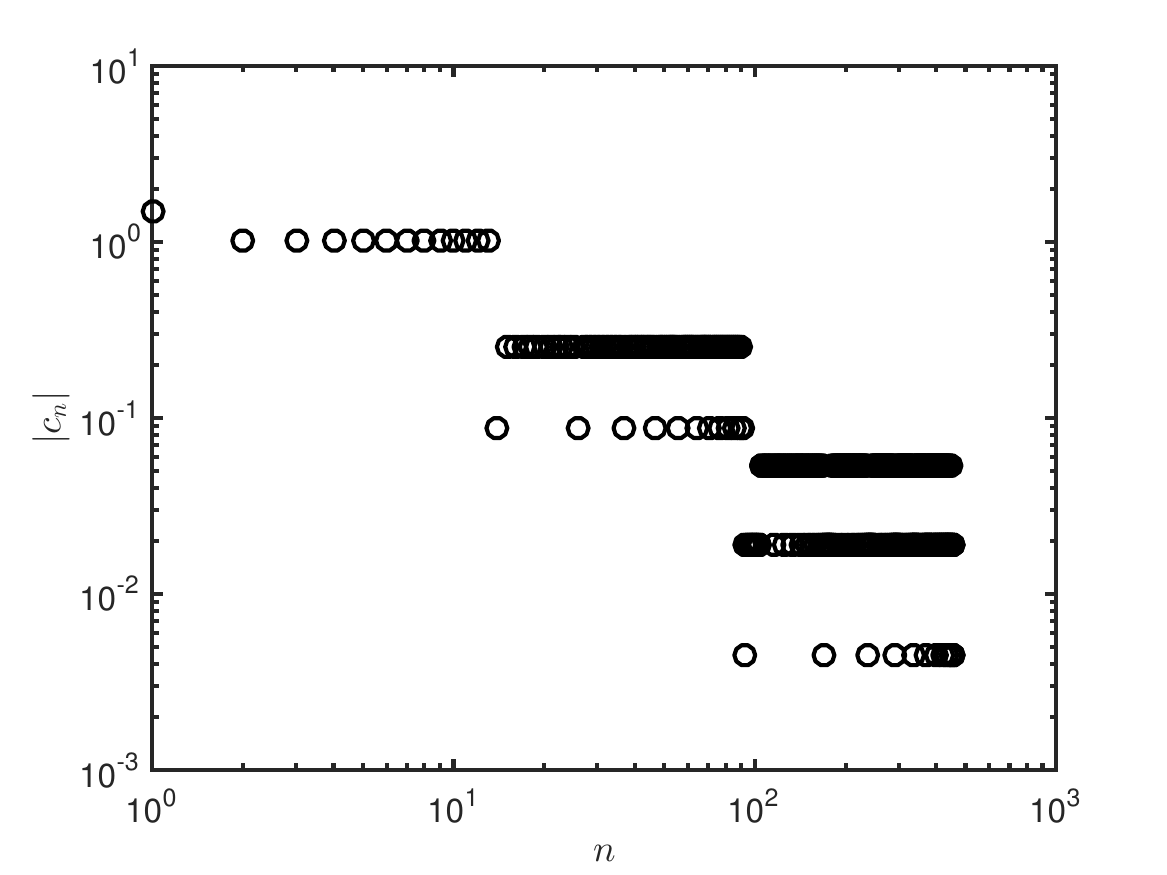}}\quad
\subfigure[Chebyshev $|\tilde c_n|$]
{\includegraphics[width=0.45\textwidth]{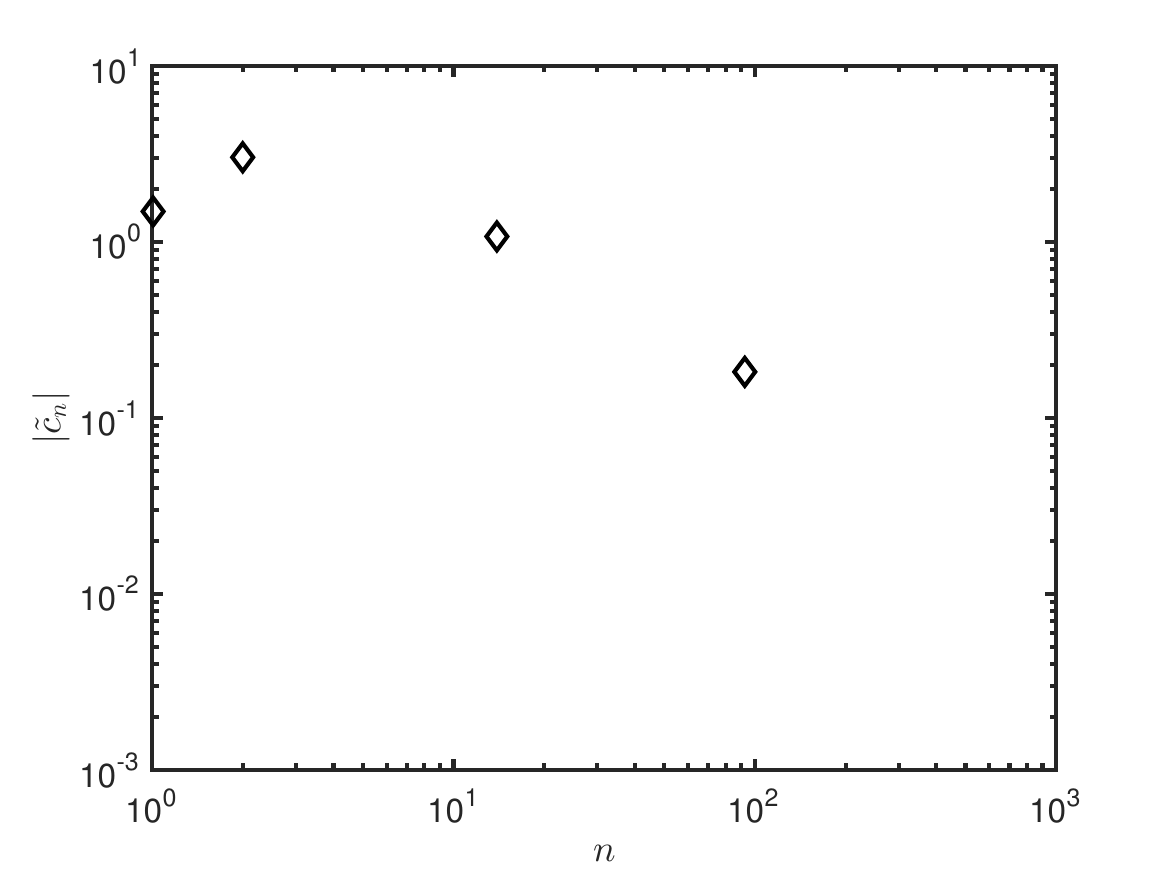}}
\caption{Results for the ridge function. Absolute values of exact coefficients 
 $c_n$ and coefficients $\tilde c_n$ after rotations using $180$ samples.}
\label{fig:ex1_coef}
\end{figure}

Figure~\ref{fig:ex1_coef} compares the absolute values of exact coefficients
$c_n$ and the coefficients $\tilde c_n$ after nine iterations using $180$
samples. In this figure we, do not present $\tilde c_n$ with absolute value 
smaller than $10^{-3}$ because they are more than two magnitudes smaller than 
the dominating ones. As demonstrated in Figure~\ref{fig:ex1_coef}, the iterative
rotation creates a much sparser representation of $u$. Thus, the efficiency of
compressive sensing method is substantially enhanced. Notably, this is a special
example in that ridge function has very good low-dimensional structure (it is a
one-dimensional function after an appropriate linear transform). In general,
many systems does not have this ideal structure, and the improvement afforded
by iterative rotations usually diminishes after two to three iterations.


\subsection{Function with high compressibility}
\label{subsec:ex2}
\begin{figure}[h]
\centering
\includegraphics[width=0.45\textwidth]{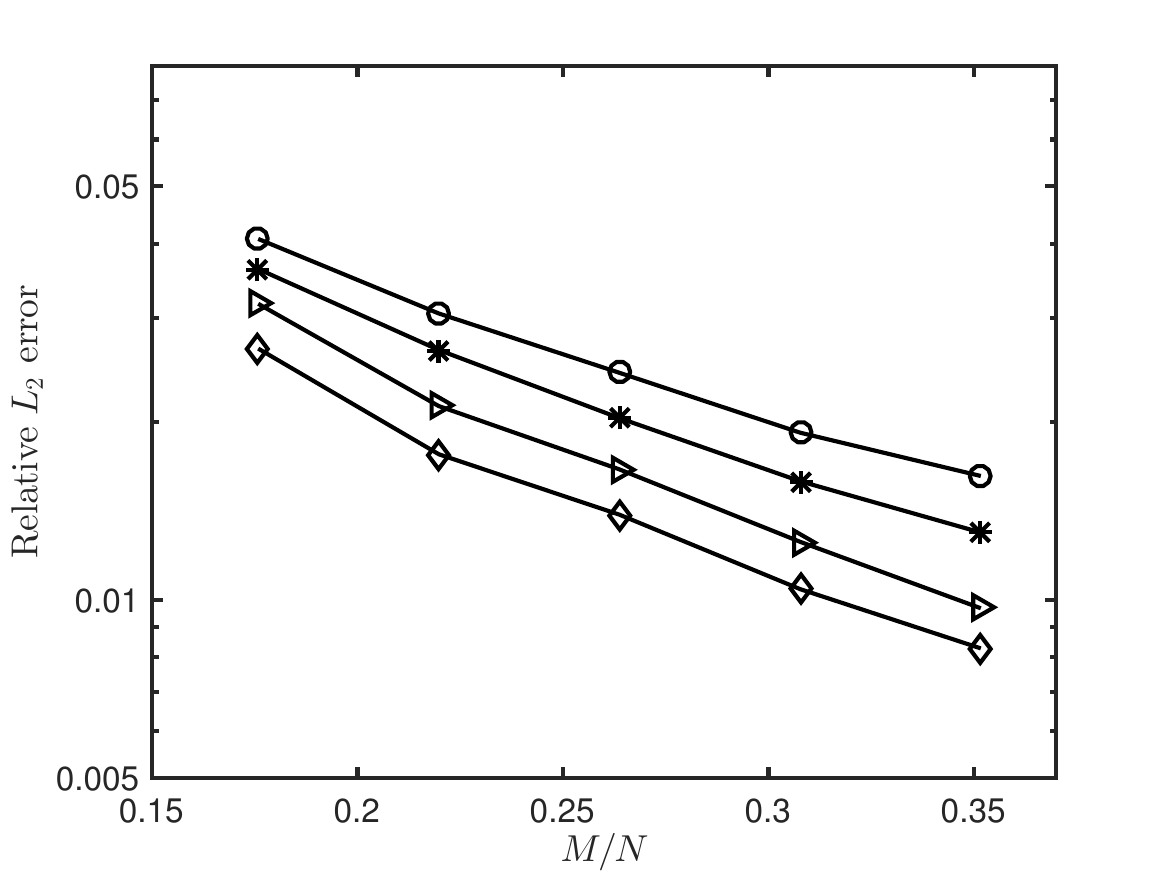}\quad
\includegraphics[width=0.45\textwidth]{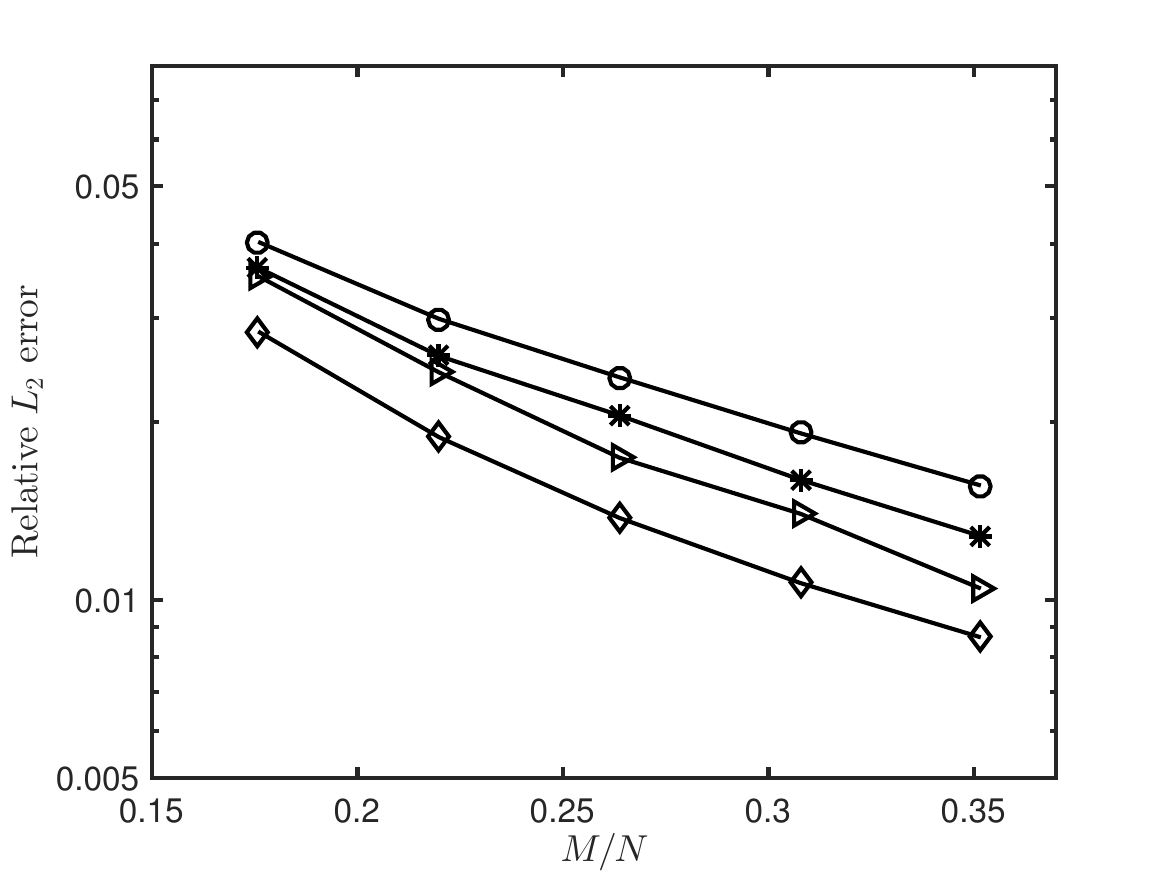}
\caption{Results for the highly compressible function. Left: Legendre polynomial expansion
  (when $\xi_i$ are i.i.d. uniform random variables). Right: Chebyshev polynomial expansion
  (when $\xi_i$ are i.i.d. Chebyshev random variables). ``$\circ$": standard $\ell_1$,
  ``$\ast$": re-weighted $\ell_1$, ``$\triangleright$": rotated $\ell_1$, ``$\diamond$":
  re-weighted+rotated $\ell_1$.}
\label{fig:ex2_rmse}
\end{figure}

\begin{figure}[h]
\centering
\subfigure[Legendre $|c_n|$]
{\includegraphics[width=0.32\textwidth]{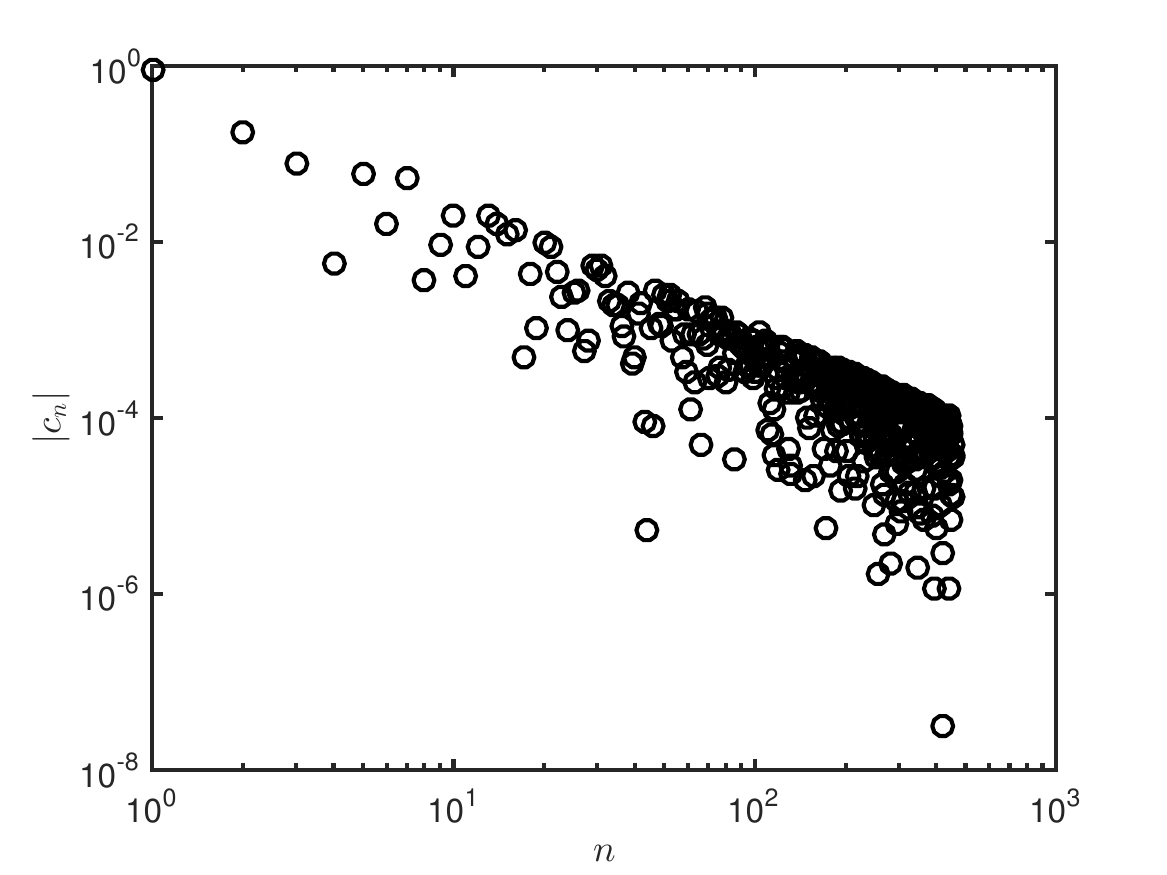}}~
\subfigure[Legendre $|\tilde c_n|$]
{\includegraphics[width=0.32\textwidth]{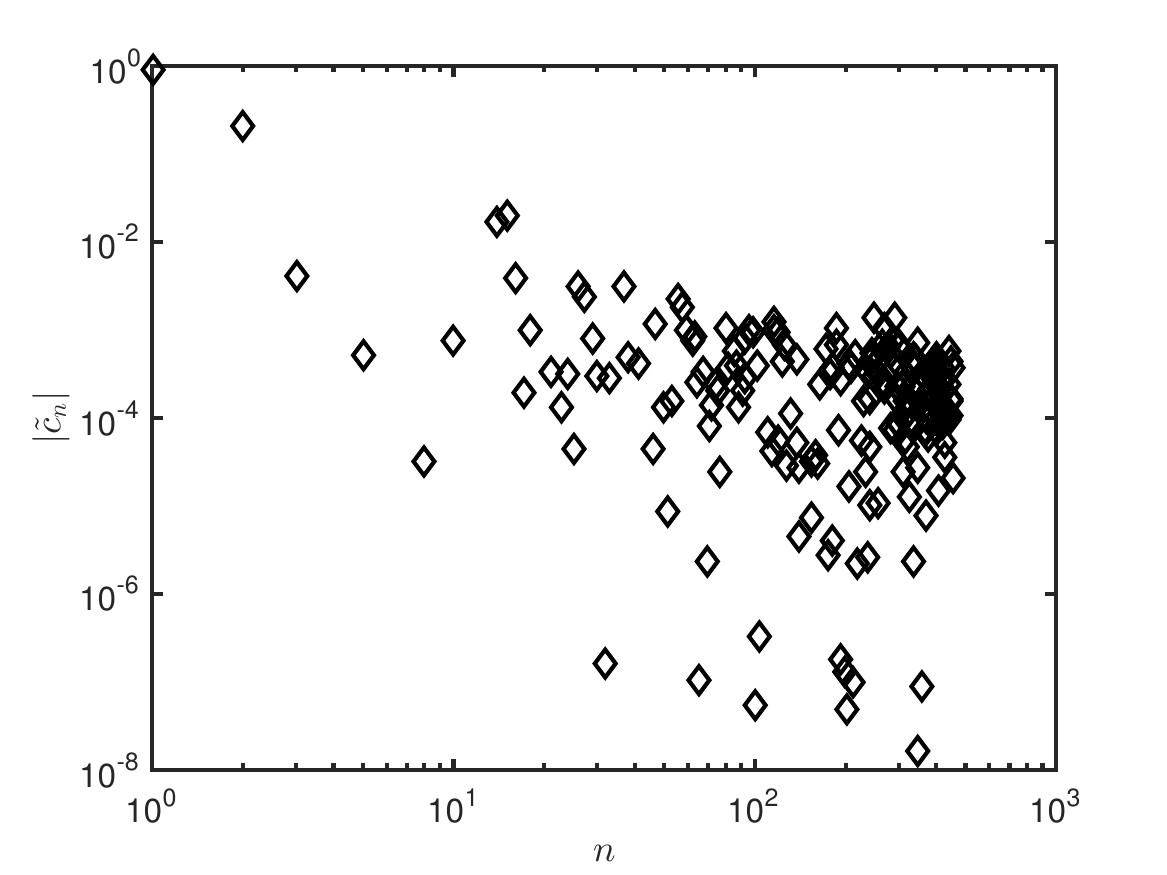}}~
\subfigure[Legendre comparison of sparsity]
{\includegraphics[width=0.32\textwidth]{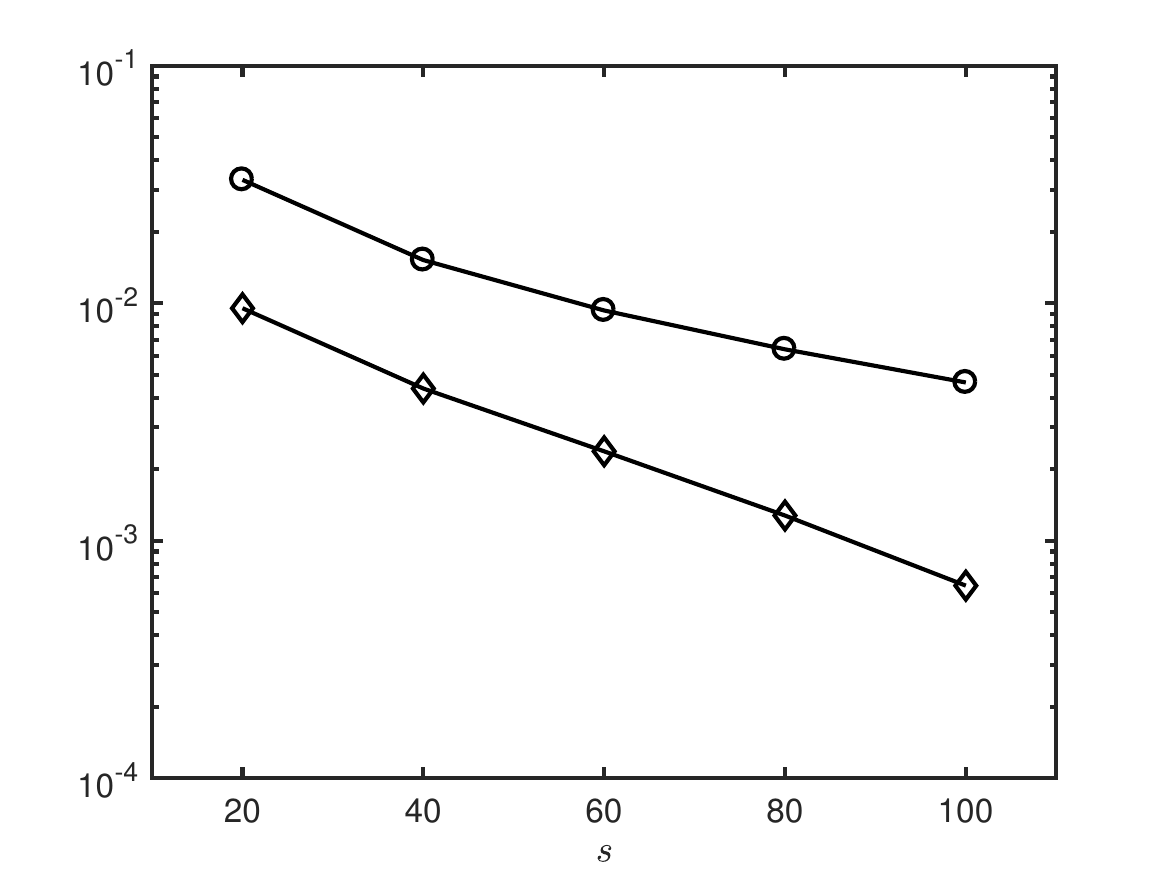}} \\
\subfigure[Chebyshev $|c_n|$]
{\includegraphics[width=0.32\textwidth]{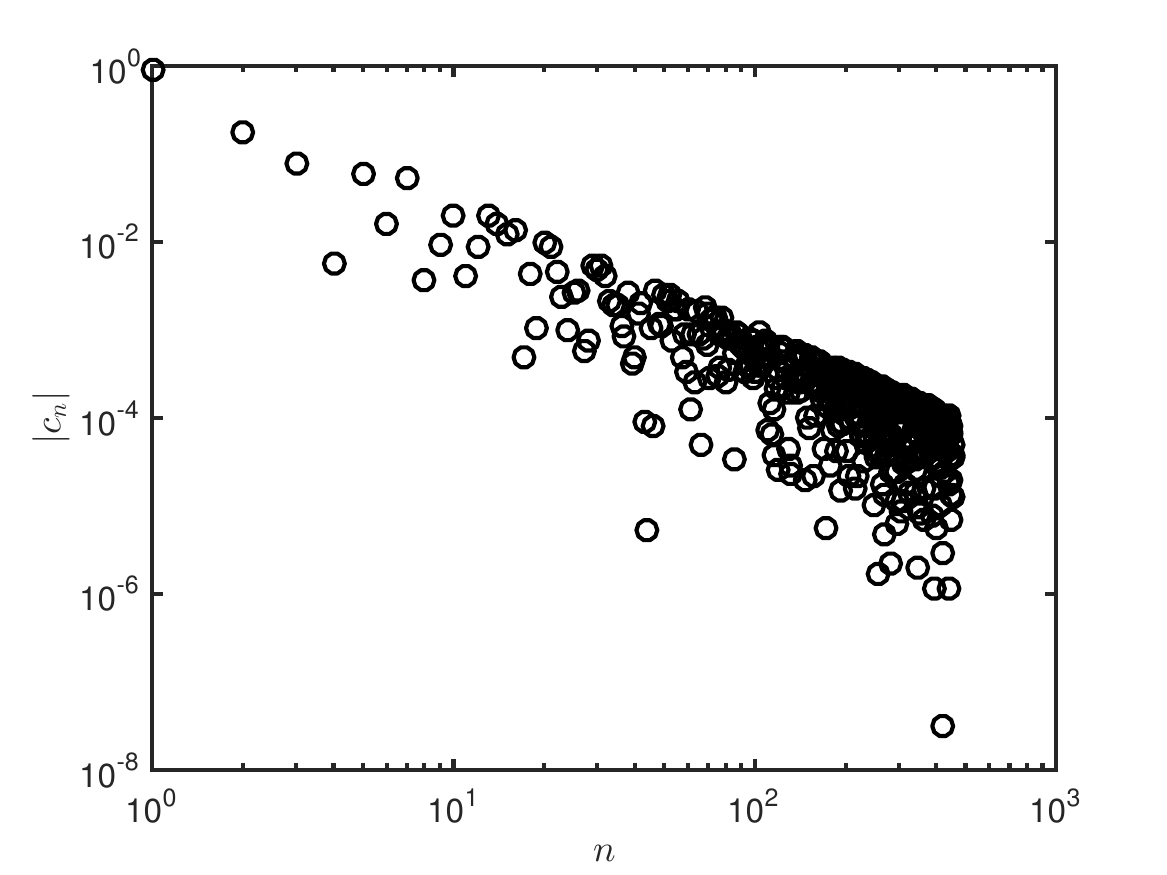}}~
\subfigure[Chebyshev $|\tilde c_n|$]
{\includegraphics[width=0.32\textwidth]{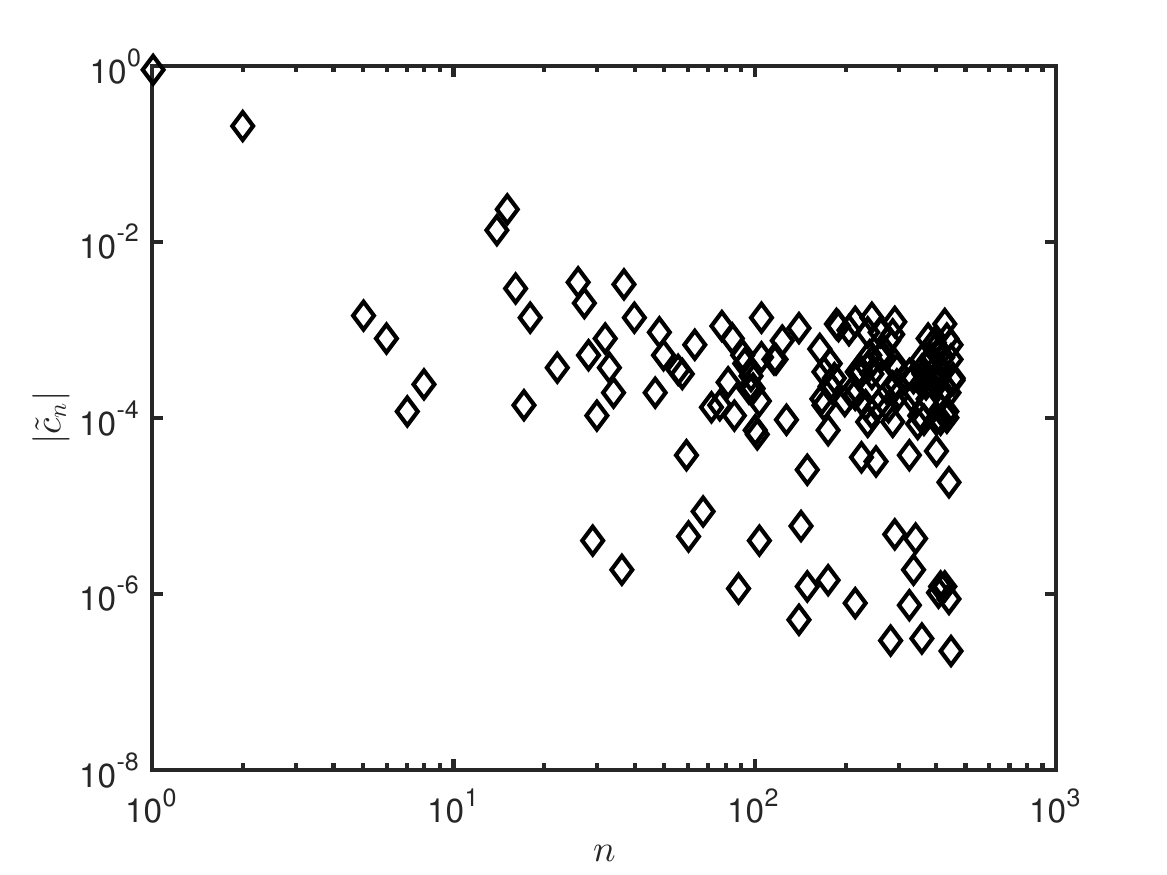}}~
\subfigure[Chebyshev comparison of sparsity]
{\includegraphics[width=0.32\textwidth]{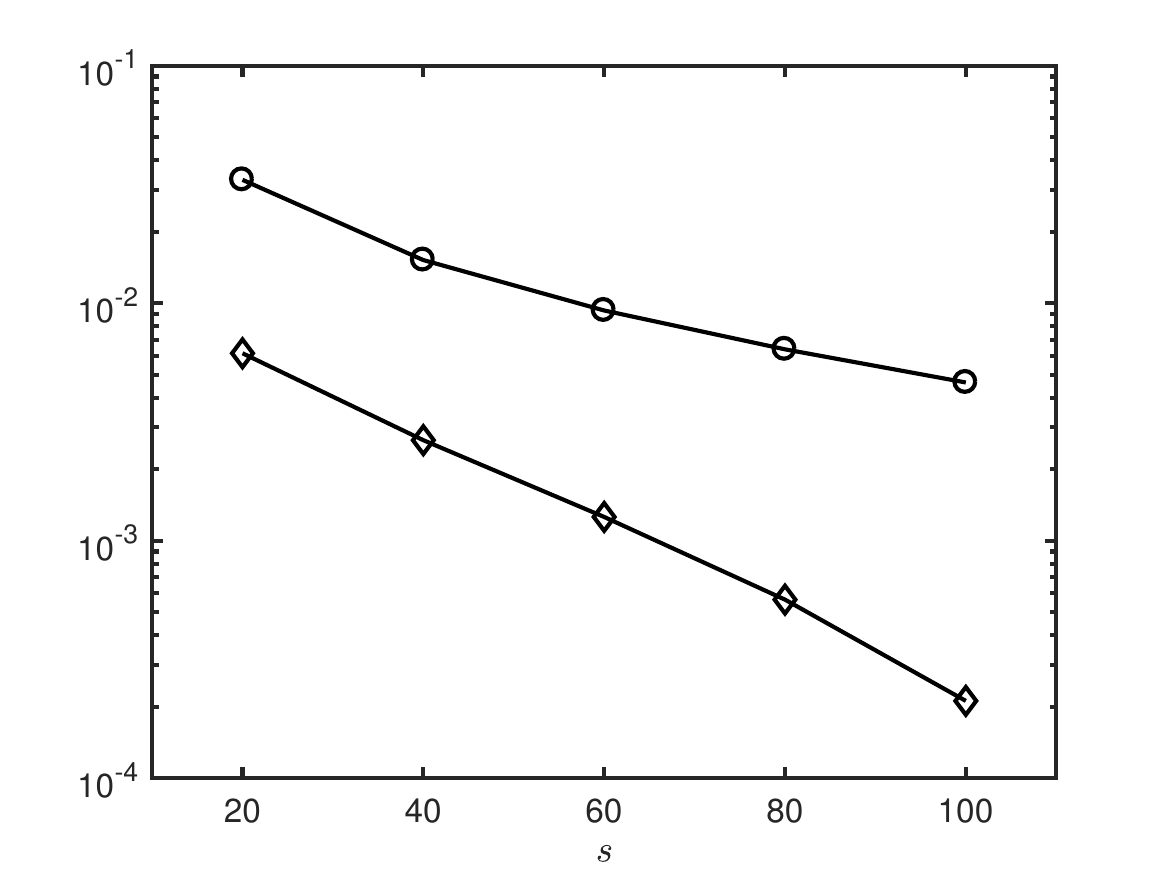}}
\caption{Results for the highly compressible function. Left column: absolute
  values of exact coefficients $c_n$; middle column: absolute values of
  coefficients $\tilde c_n$ after rotations using $160$ samples;
right column: comparison of $\dfrac{\Vert\bm c-\bm c_s\Vert_1}{\sqrt{s}}$
(``$\circ$") and $\dfrac{\Vert\tilde{\bm c}-\tilde{\bm c}_s\Vert_1}{\sqrt{s}}$
(``$\diamond$") with different $s$.}
\label{fig:ex2_coef}
\end{figure}

Consider the following function:
\begin{equation}\label{eq:ex2}
u(\bx) = \sum_{|\ba|=0}^P c_{\ba}\psi_{\ba}(\bx)
       = \sum_{n=1}^Nc_n\psi_n(\bx), \quad
\bx =  (\xi_1,\xi_2,\cdots,\xi_{d})\trans,
\end{equation}
where, $\psi_{\ba}$ are normalized multivariate Legendre or
Chebyshev polynomials, $d=12, P=3, N=455$, and the coefficients $c_n$ are chosen
as uniformly distributed random numbers,
\begin{equation}
c_n = \zeta/n^{1.5}, \quad \zeta\sim \mathcal{U}[0,1].
\end{equation}
For this example, we generate $N$ samples of $\zeta$: 
$\zeta^1,\zeta^2,\cdots,\zeta^N$ then divide them by $n^{1.5}, n=1,2,\cdots,N$
to obtain a random ``compressible signal" $\bm c$. The integrals for the 
relative error are computed by a level-$4$ sparse grid method and, therefore,
are exact. Figure \ref{fig:ex2_rmse} shows the relative $L_2$ errors obtained
by applying our iterative rotation technique to the re-weighted $\ell_1$ approach.
Apparently, introduction of the iterative rotation approach improves the 
accuracy. 
Figure~\ref{fig:ex2_coef} presents a comparison of the absolute values of 
entries of $\bm c$ and $\tilde{\bm c}$ (using $160$ samples). The main 
improvement is that the number of coefficients with magnitudes larger than 
$0.01$ decreased. Also, $c_n$ cluster around the curve $c_n=1/n^{1.5}$ as we set
them in this way, while many $\tilde c_n$ appear below this curve, especially
when $n$ is large. In Figure~\ref{fig:ex2_coef}, we also compare the values of 
$\dfrac{\Vert\bm c-\bm c_s\Vert_1}{\sqrt{s}}$ and
$\dfrac{\Vert\tilde{\bm c}-\tilde{\bm c}_s\Vert_1}{\sqrt{s}}$ to demonstrate
quantitatively the enhancement of the sparsity after rotations.


\subsection{Elliptic equation}
\label{subsec:ex3}
\begin{figure}[h]
\centering
\includegraphics[width=0.45\textwidth]{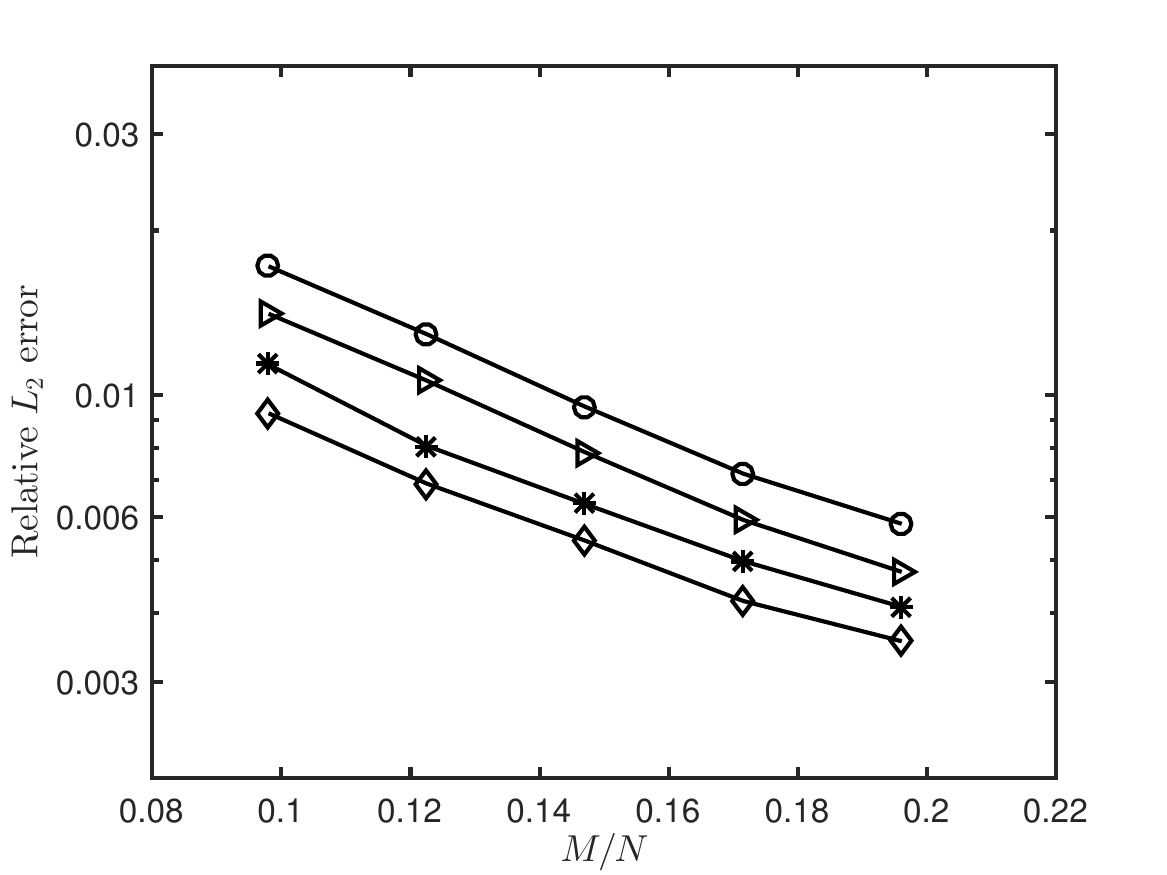}\quad
\includegraphics[width=0.45\textwidth]{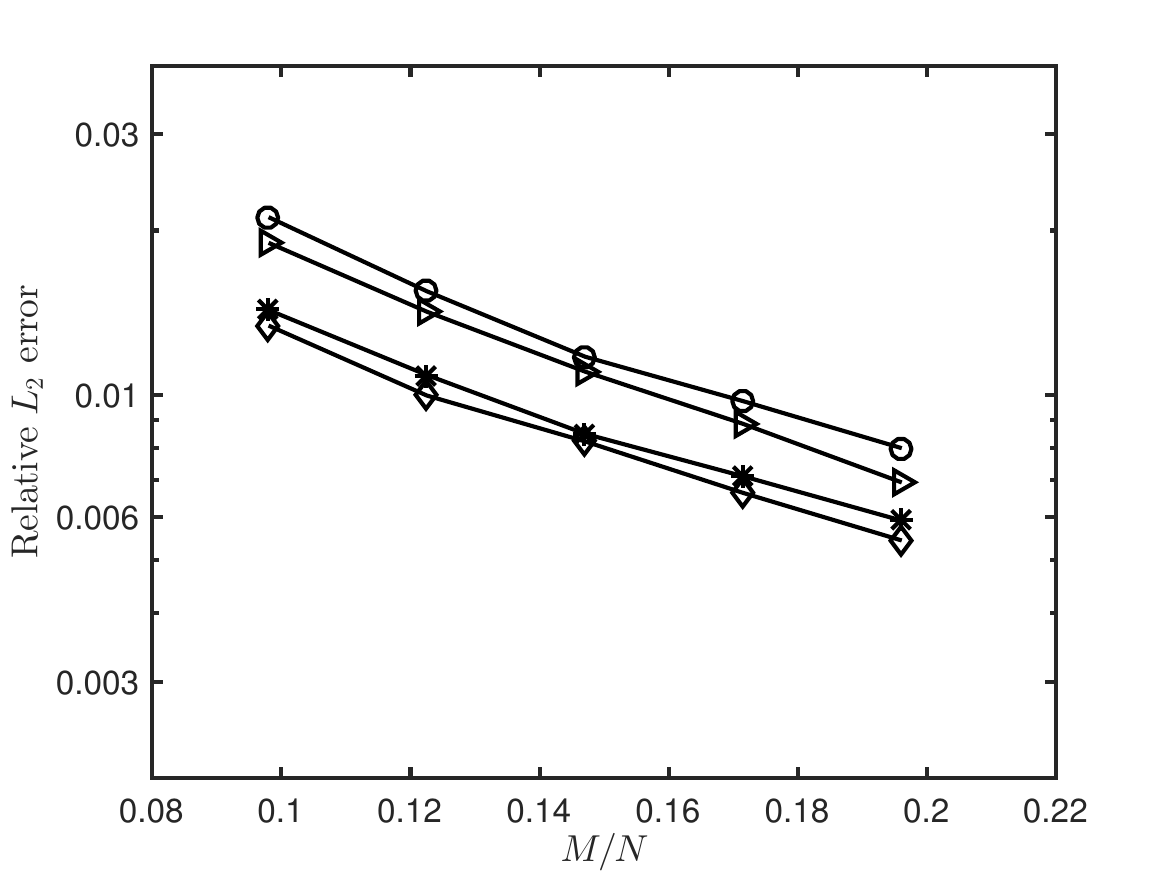}
\caption{Results for the elliptic equation. Left: Legendre polynomial expansion
  (when $\xi_i$ are i.i.d. uniform random variables). Right: Chebyshev 
  polynomial expansion (when $\xi_i$ are i.i.d. Chebyshev random variables).
  ``$\circ$": standard $\ell_1$, ``$\ast$": re-weighted $\ell_1$, 
  ``$\triangleright$": rotated $\ell_1$, ``$\diamond$": re-weighted+rotated $\ell_1$.}
\label{fig:ex3_rmse}
\end{figure}

\begin{figure}[h]
\centering
\subfigure[Legendre $|c_n|$]
{\includegraphics[width=0.32\textwidth]{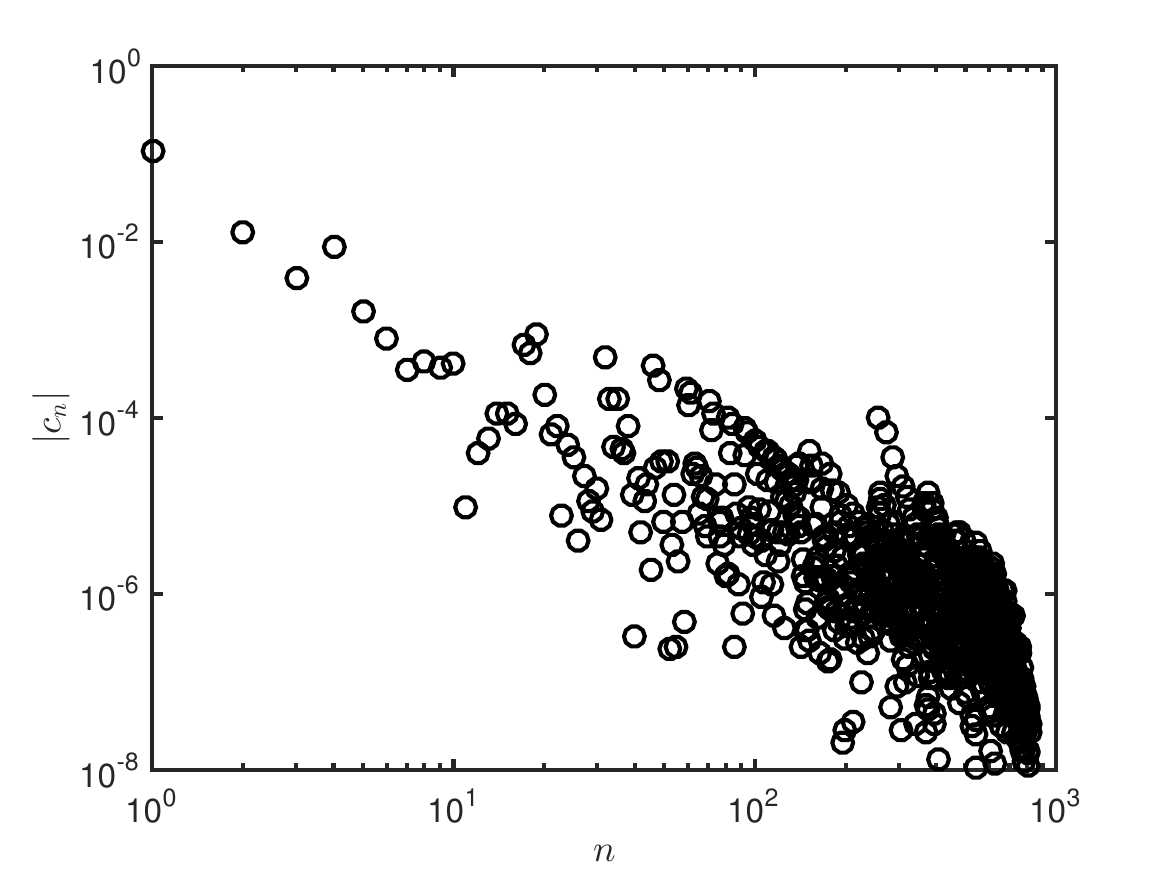}}~
\subfigure[Legendre $|\tilde c_n|$]
{\includegraphics[width=0.32\textwidth]{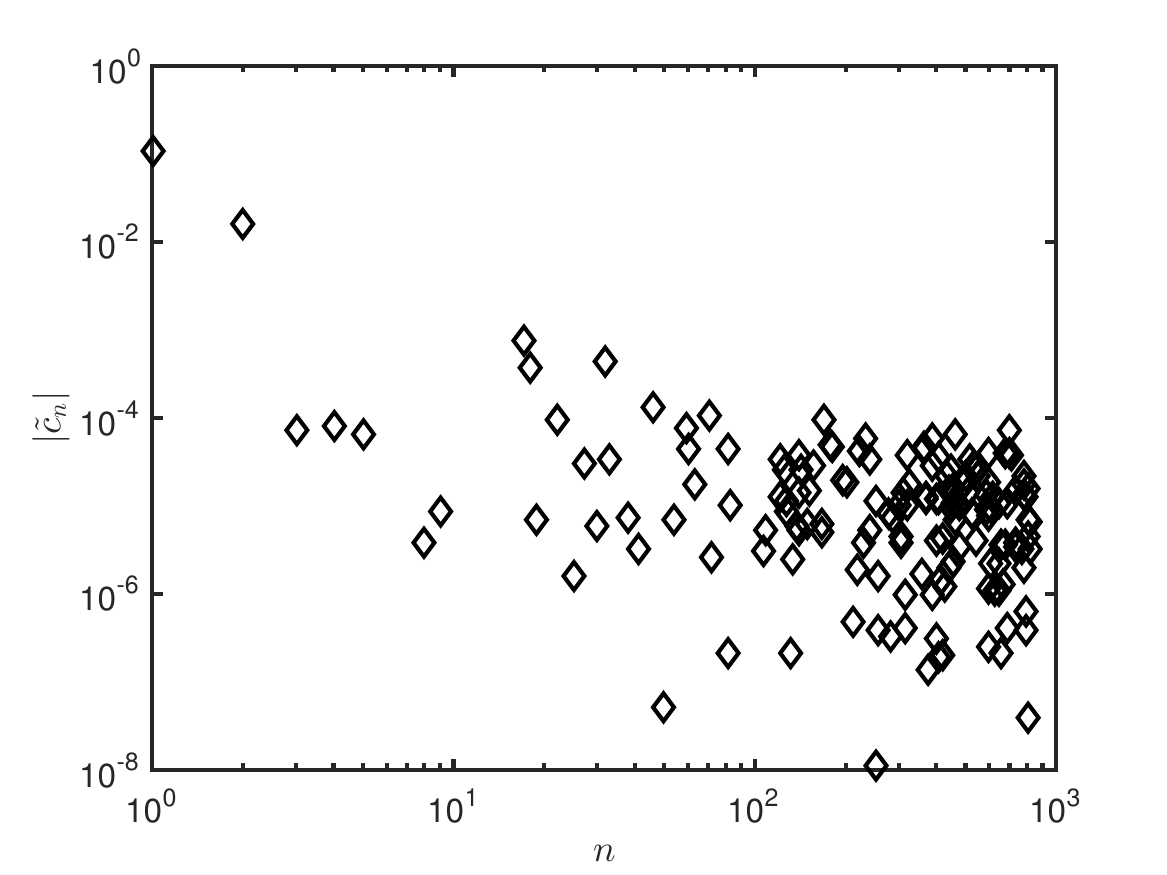}}~
\subfigure[Legendre comparison of sparsity]
{\includegraphics[width=0.32\textwidth]{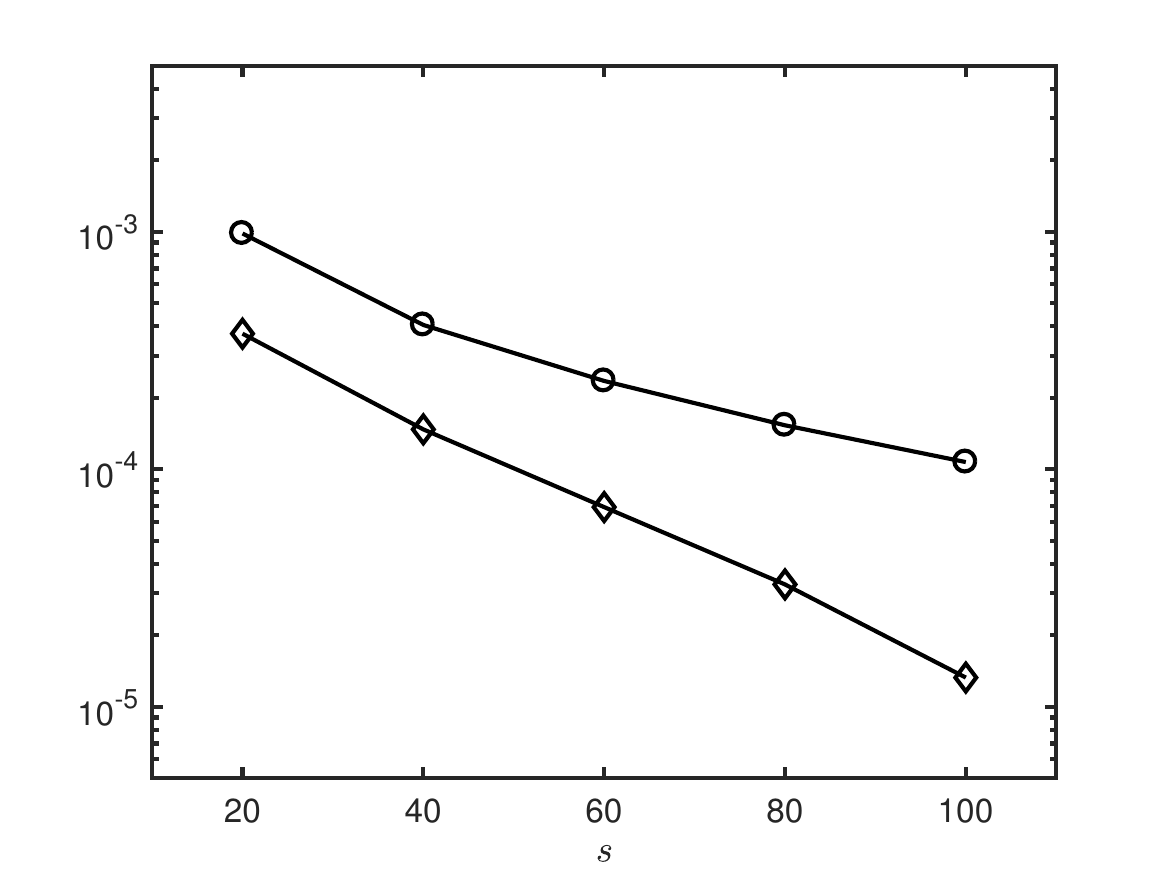}}\\
\subfigure[Chebyshev $|c_n|$]
{\includegraphics[width=0.32\textwidth]{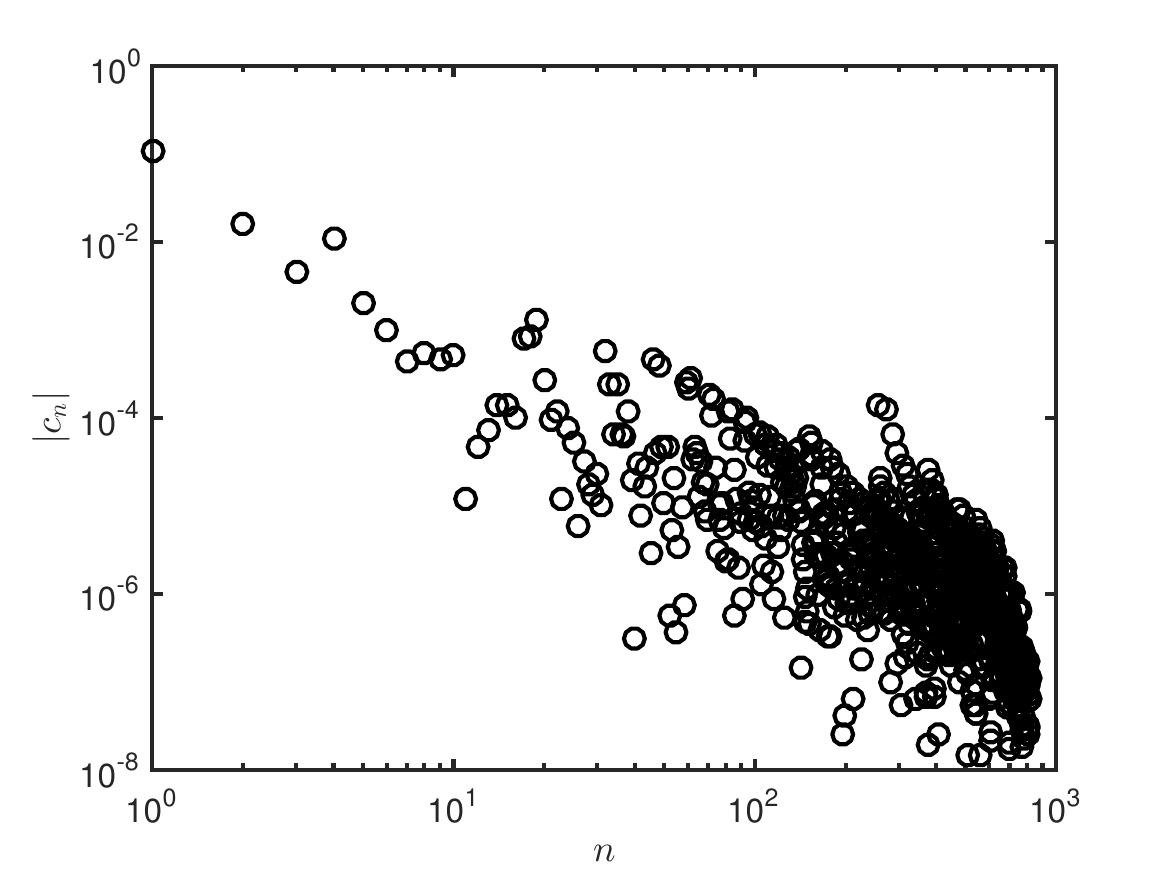}}~
\subfigure[Chebyshev $|\tilde c_n|$]
{\includegraphics[width=0.32\textwidth]{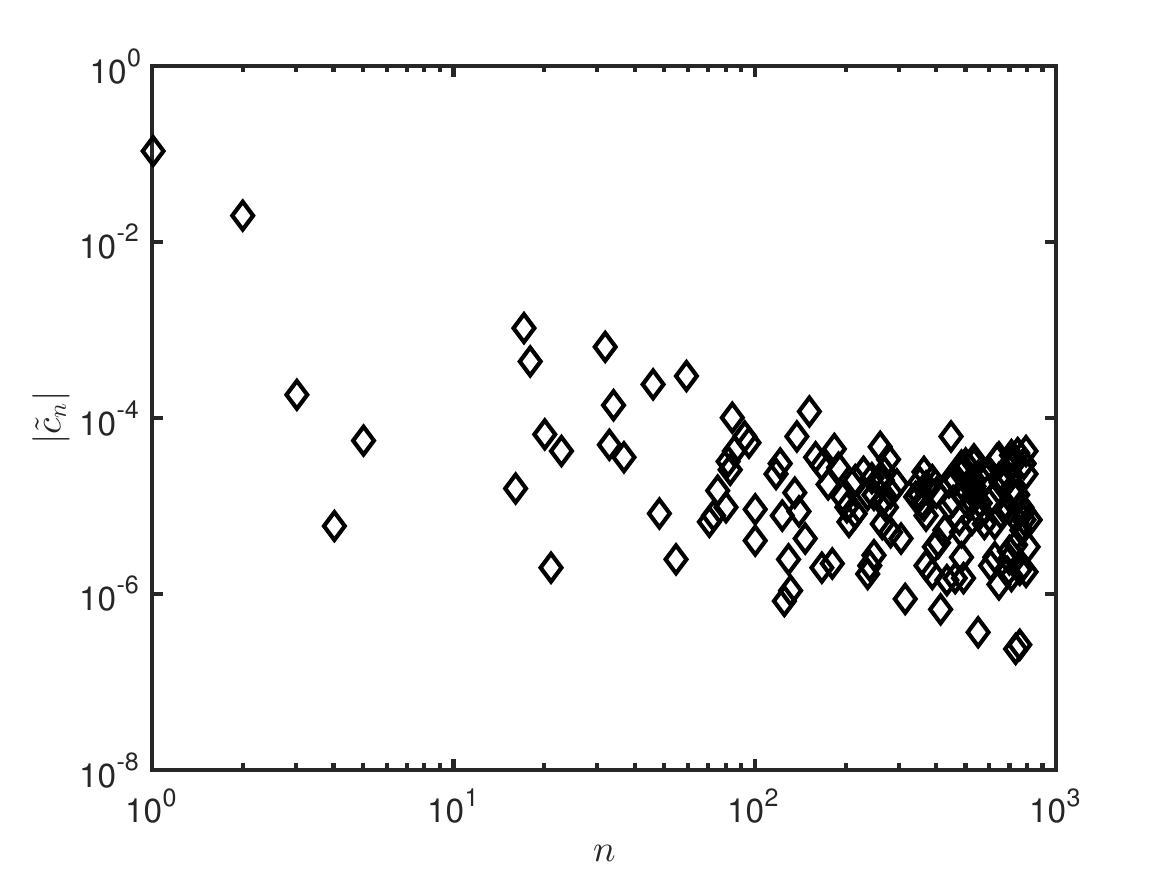}}~
\subfigure[Chebyshev comparison of sparsity]
{\includegraphics[width=0.32\textwidth]{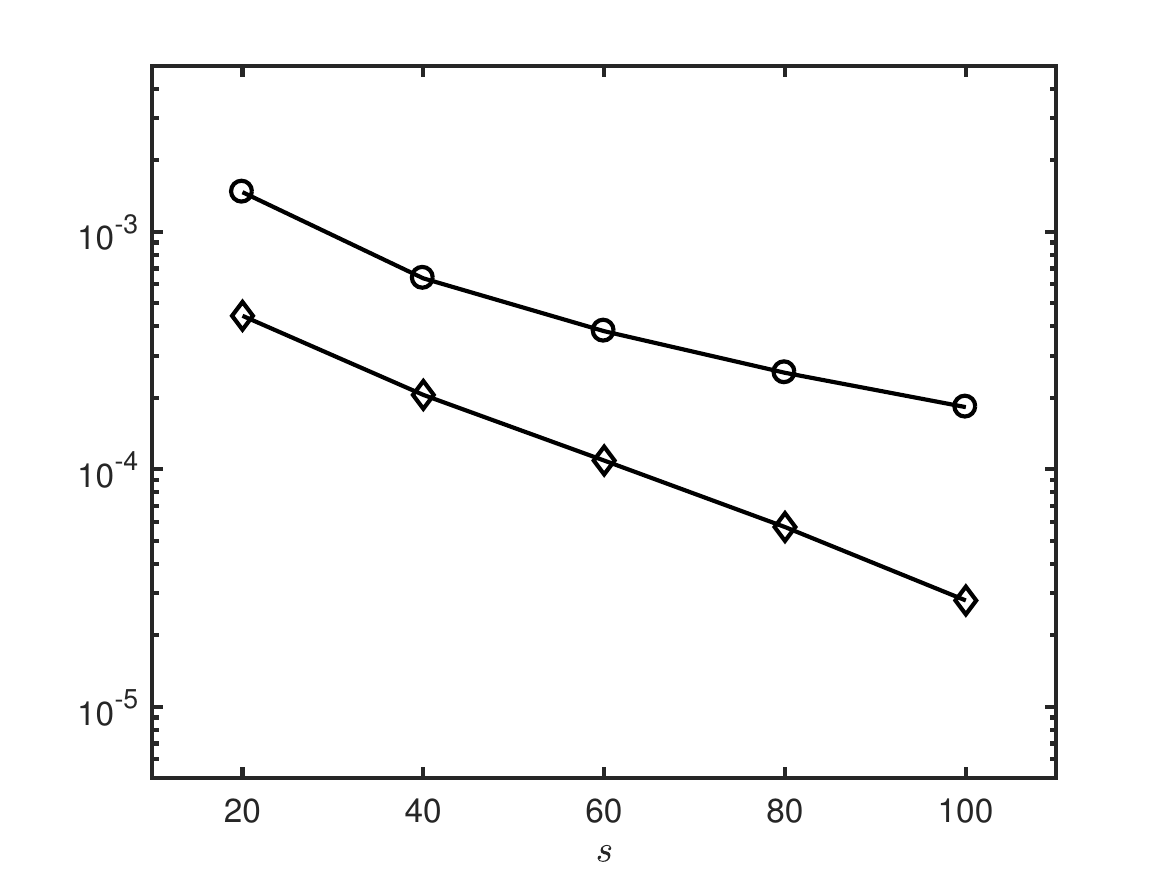}}
\caption{Results for the elliptic equation. Left column: absolute
  value of exact coefficients $|c_n|$; middle column: absolute values of
  coefficients $\tilde c_n$ after rotations using $160$ samples;
right column: comparison of $\dfrac{\Vert\bm c-\bm c_s\Vert_1}{\sqrt{s}}$
(``$\circ$") and $\dfrac{\Vert\tilde{\bm c}-\tilde{\bm c}_s\Vert_1}{\sqrt{s}}$
(``$\diamond$") with different $s$.}
\label{fig:ex3_coef}
\end{figure}

Next, we consider a one-dimensional elliptic differential equation with a random 
coefficient:
\begin{equation}\label{eq:ellip}
\begin{aligned}
-\frac{d}{dx} \left( a(x;\bx)\frac{d u(x;\bx)}{dx} \right) = 1, & \quad x \in (0,1) \\
u(0) = u(1) = 0, &
\end{aligned}
\end{equation}
where $a(x;\bx)$ is a log-normal random field based on Karhunen-Lo\`{e}ve (KL) 
expansion:
\begin{equation}\label{eq:kl}
a(x;\bx) = a_0(x) +
\exp\left(\sigma\sum_{i=1}^d\sqrt{\lambda_i}\phi_i(x)\xi_i\right), 
\end{equation}
where $\{\xi_i\}$ are i.i.d. random variables,
$\{\lambda_i\}_{i=1}^d$, and $\{\phi_i(x)\}_{i=1}^d$ are the largest eigenvalues
and corresponding eigenfunctions of the exponential covariance kernel:
\begin{equation}\label{eq:exp_kernel}
C(x,x') = \exp\left(-\dfrac{|x-x'|}{l_c}\right).
\end{equation}
In the KL expansion, $\lambda_i$ denotes the eigenvalue of the covariance kernel
$C(x,x')$ instead of entries of $\Lambda$ in Eq.~\eqref{eq:grad}. The value of 
$\lambda_i$ and the analytical expressions for $\phi_i$ are available in the 
literature \cite{JardakSK02}. In this example, we set 
$a_0(x)\equiv 0.1, \sigma=0.5, l_c=0.2$, and $d=15$. With this setting, 
$\sum_{i=1}^d\lambda_i > 0.93\sum_{i=1}^{\infty}\lambda_i$. For each input 
sample $\bx^q$, $a$ and $u$ only depend on $x$, and the solution of the 
deterministic elliptic equation can be obtained as \cite{YangK13}:
\begin{equation}\label{eq:ellip_sol}
u(x) = u(0) + \int_0^x \dfrac{a(0)u(0)'-y}{a(y)}\dif y.
\end{equation}
By imposing the boundary condition $u(0)=u(1)=0$, we can compute $a(0)u(0)'$ as
\begin{equation}\label{eq:ellip_const}
a(0)u(0)' = \left(\int_0^1 \dfrac{y}{a(y)}\dif y\right) \Big /
            \left(\int_0^1 \dfrac{1}{a(y)}\dif y\right).
\end{equation}
The integrals in Eqs.~\eqref{eq:ellip_const} and \eqref{eq:ellip_sol} are
obtained by highly accurate numerical integration. For this example, we choose 
the QoI to be $u(x;\bx)$ at $x=0.35$. We aim to build a 
third-order Legendre (or Chebyshev) polynomial expansion that includes $N=816$
basis functions. The relative error is approximated by a level-$6$ sparse grid 
method. Figure \ref{fig:ex3_rmse} shows that accuracy of the re-weighted $\ell_1$
and the iteratively rotated $\ell_1$ method are very close in this case. In the 
Legendre polynomial expansion, the 
incorporation of iterative rotation improves the performance of the other methods.
In the Chebyshev polynomial expansion, the improvement is minimal. This is
related to the compromise of the property of $\tensor\Psi$.
Figure~\ref{fig:ex3_coef} depicts a comparison of $\bm c$ and $\tilde{\bm c}$,
which shows the improvement of the sparsity in the similar manner as in function
with high compressibility in Section~\ref{subsec:ex2}. Here, we plot 
coefficients with absolute values larger than $10^{-8}$ for demonstration
purposes because other entries are negligible in the comparison of sparsity. 
The right column of Figure~\ref{fig:ex3_coef} also quantitatively illustrates
sparsity enhancement.


\subsection{Korteweg-de Vries equation}
\begin{figure}[h]
\centering
\includegraphics[width=0.45\textwidth]{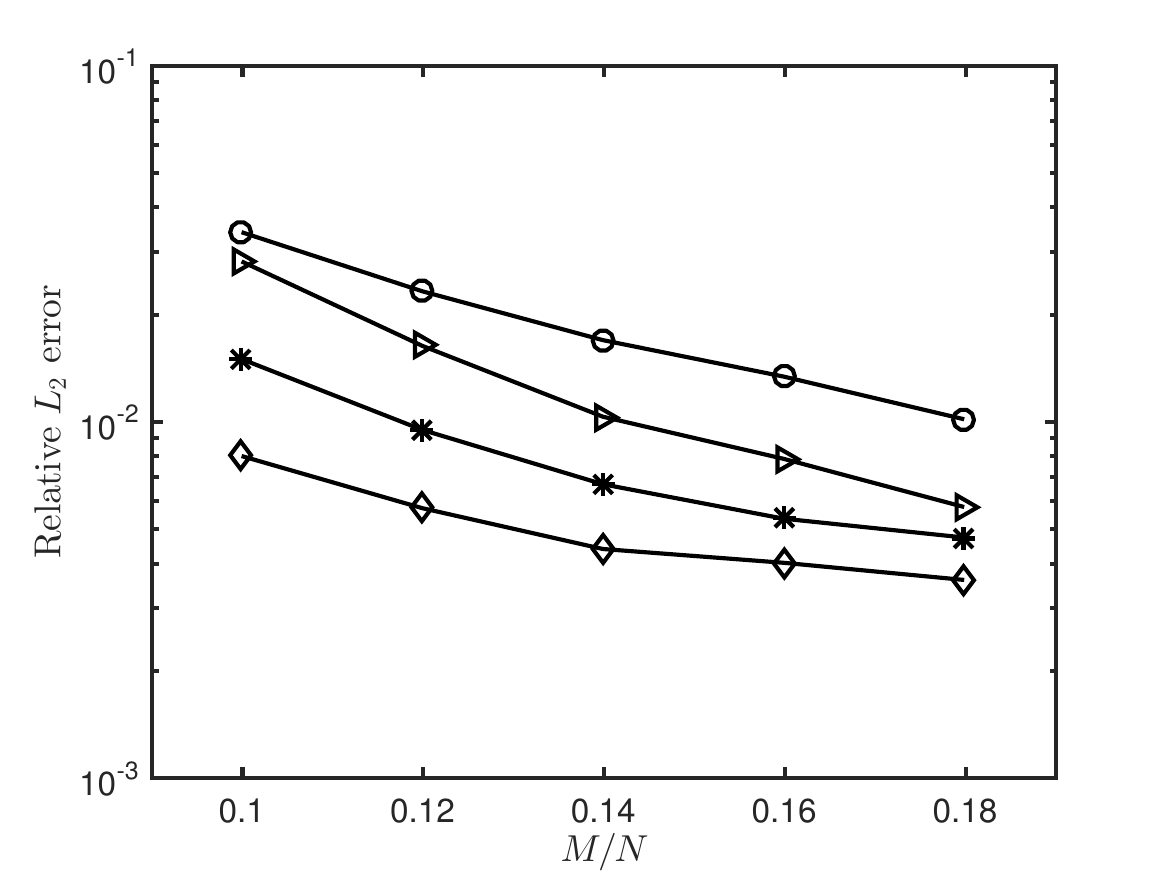}\quad
\includegraphics[width=0.45\textwidth]{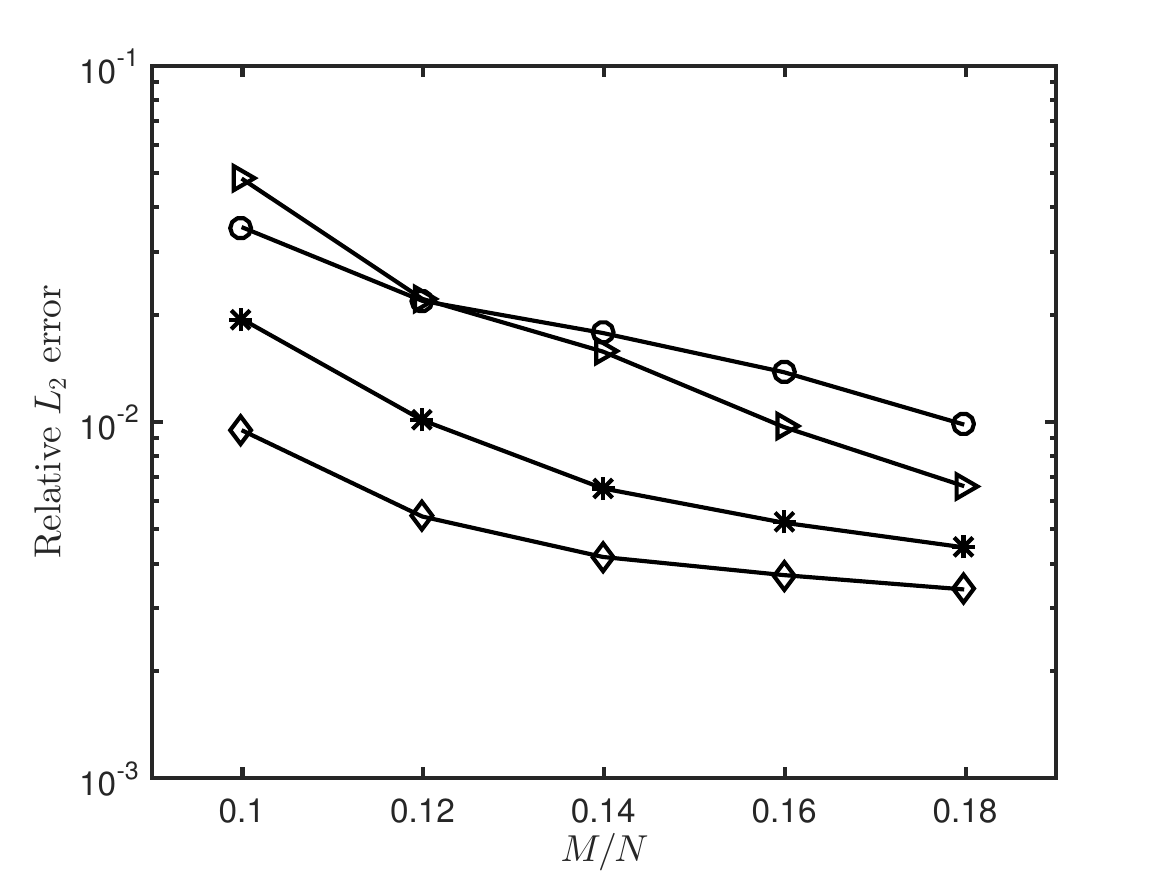}
\caption{Results for the Korteweg-de Vries equation. Left: Legendre polynomial 
  expansion (when $\xi_i$ are i.i.d. uniform random variables). Right: Chebyshev
  polynomial expansion (when $\xi_i$ are i.i.d. Chebyshev random variables).
  ``$\circ$": standard $\ell_1$, ``$\ast$": re-weighted $\ell_1$,
  ``$\triangleright$": rotated $\ell_1$, ``$\diamond$": re-weighted+rotated $\ell_1$.}
\label{fig:ex4_rmse}
\end{figure}

\begin{figure}[h]
\centering
\subfigure[Legendre $|c_n|$]
{\includegraphics[width=0.32\textwidth]{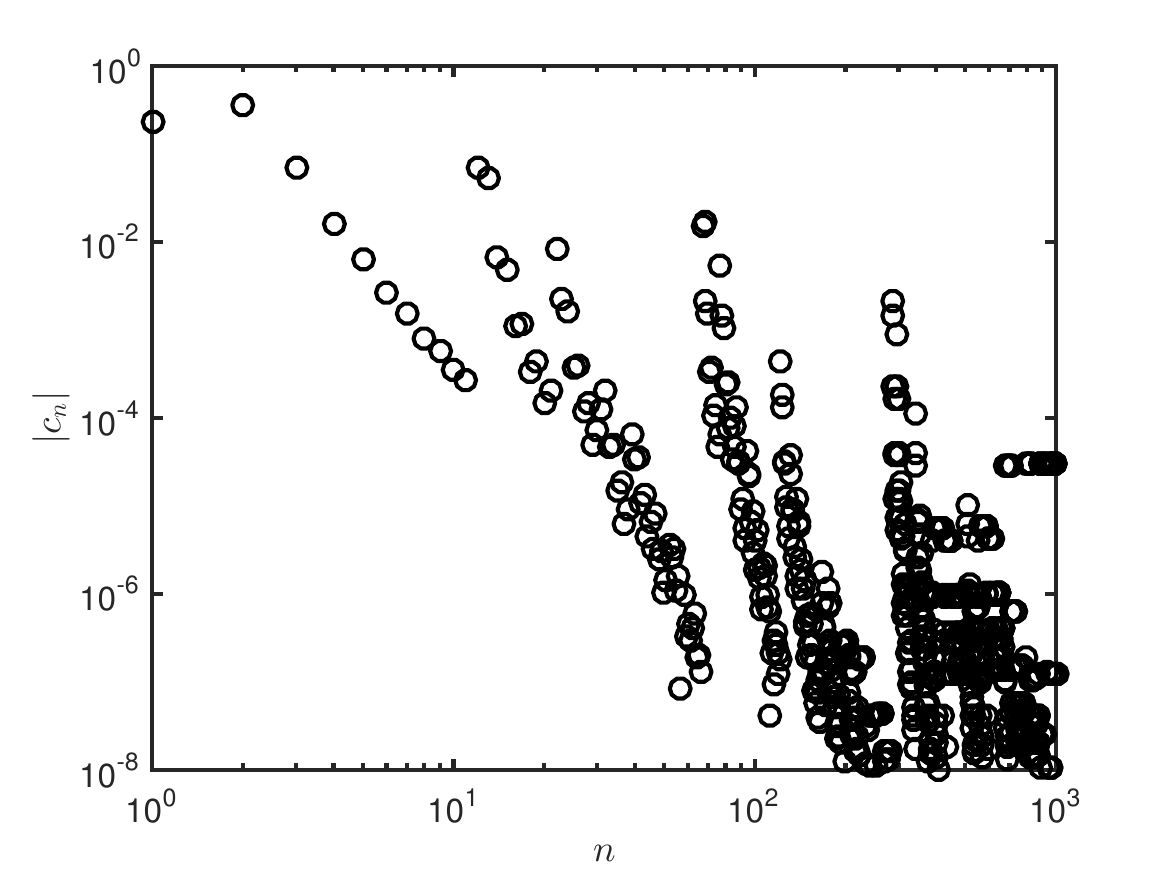}}~
\subfigure[Legendre $|\tilde c_n|$]
{\includegraphics[width=0.32\textwidth]{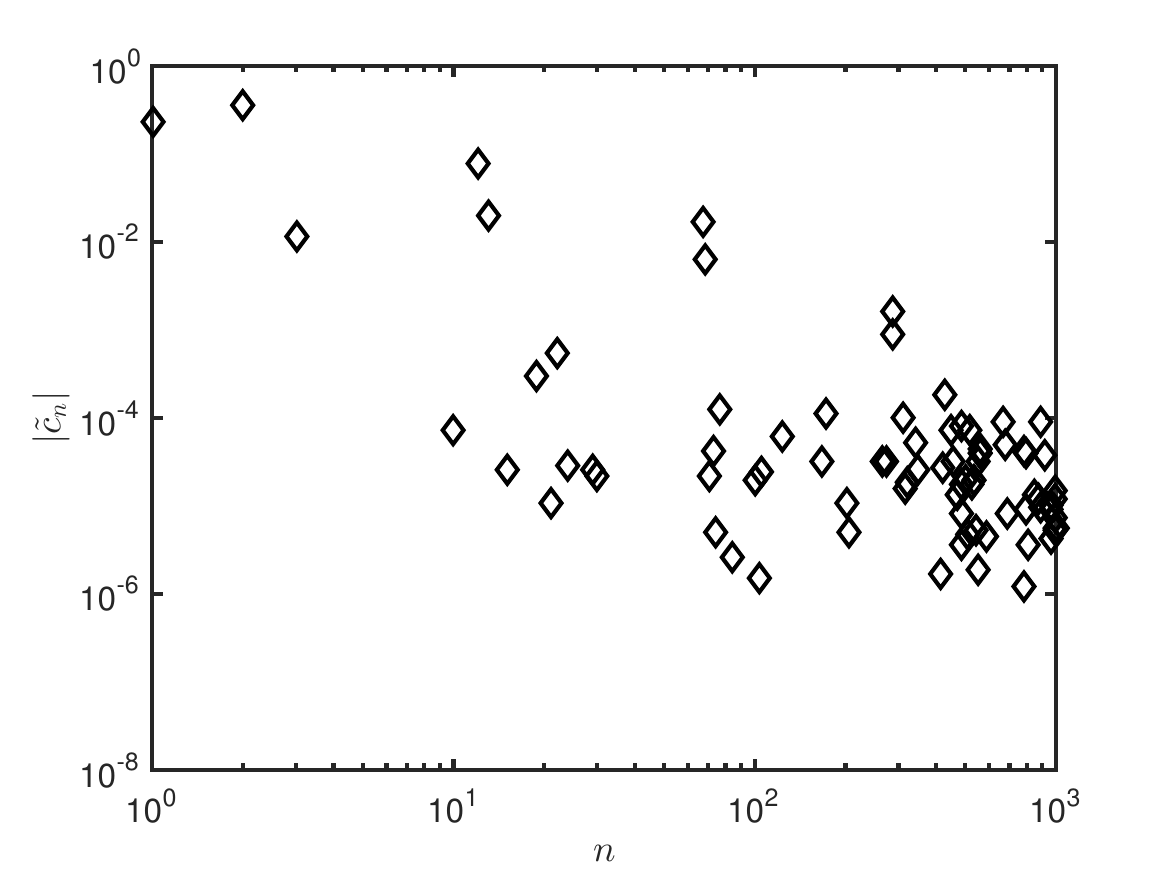}}~
\subfigure[Legendre comparison of sparsity]
{\includegraphics[width=0.32\textwidth]{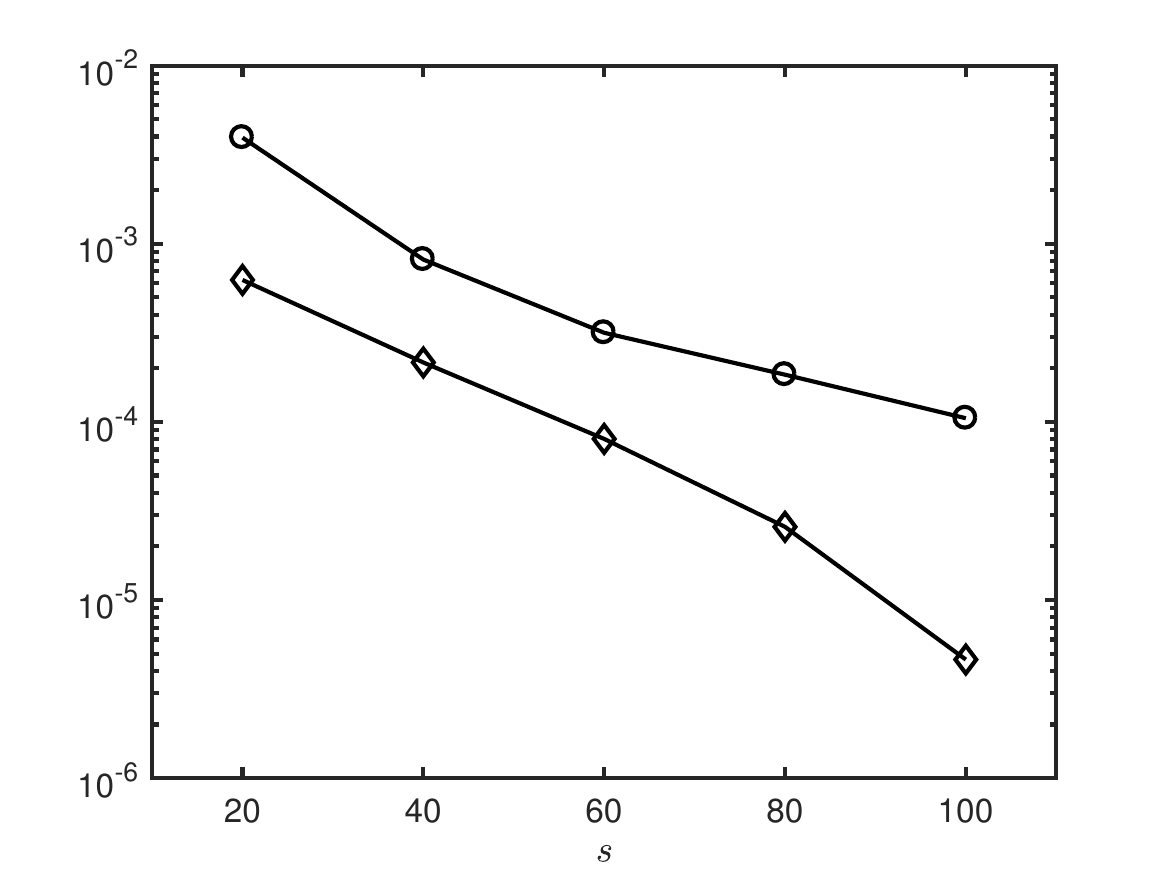}} \\
\subfigure[Chebyshev $|c_n|$]
{\includegraphics[width=0.32\textwidth]{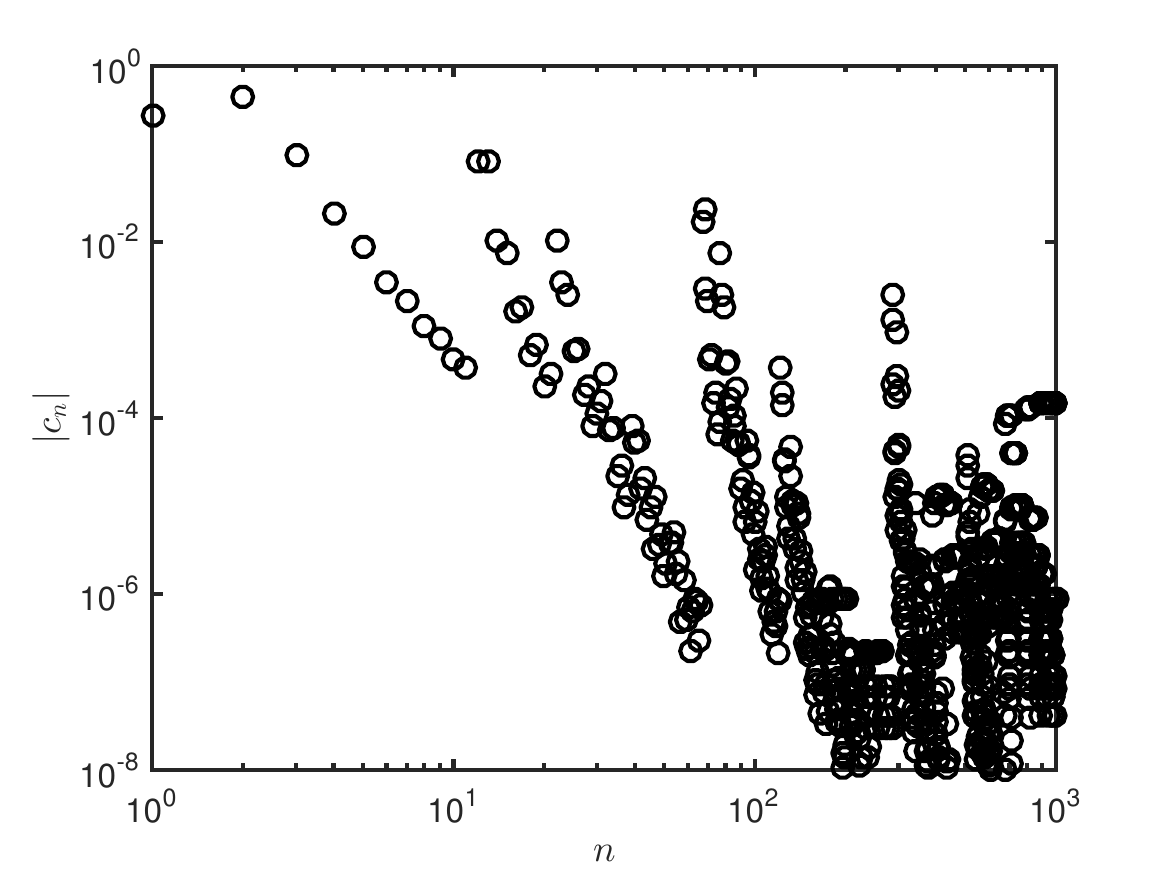}}~
\subfigure[Chebyshev $|\tilde c_n|$]
{\includegraphics[width=0.32\textwidth]{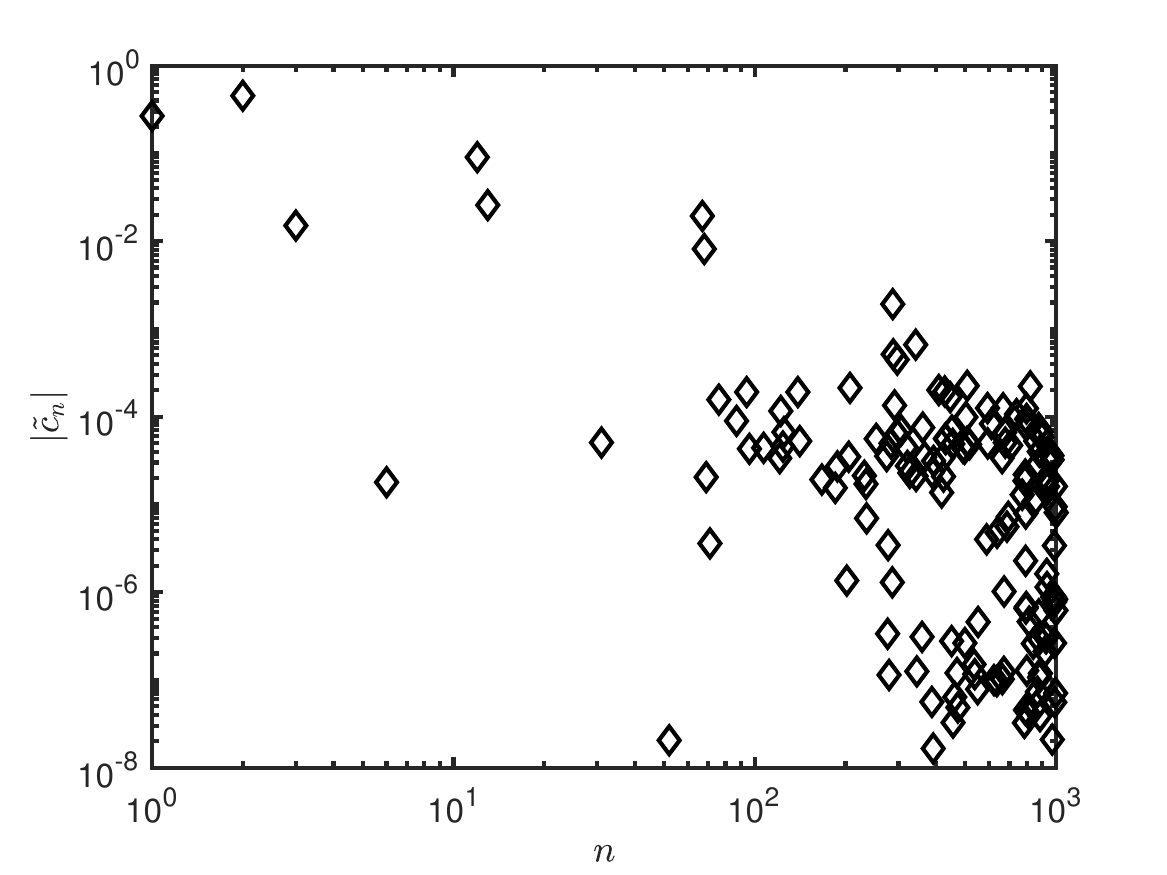}}~
\subfigure[Chebyshev comparison of sparsity]
{\includegraphics[width=0.32\textwidth]{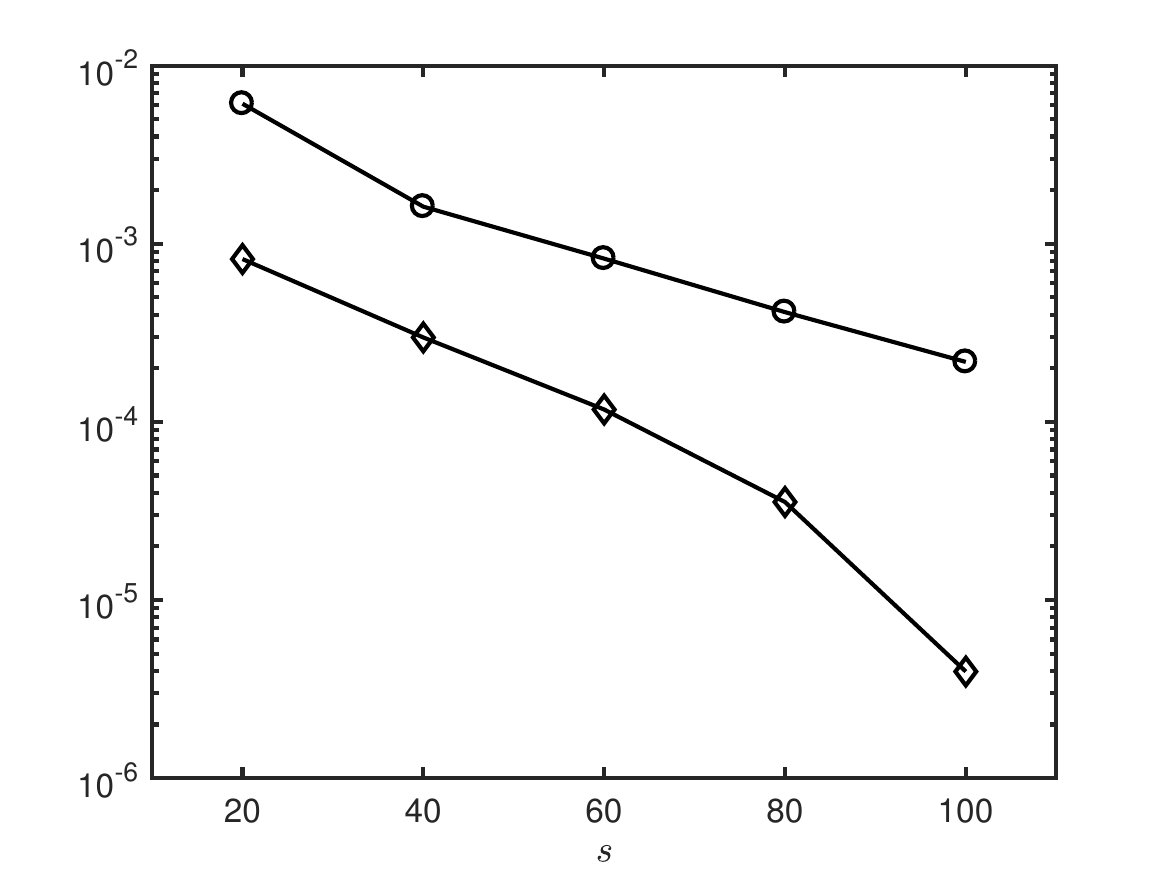}}
\caption{Results for the KdV equation. Left column: absolute
  values of exact coefficients $c_n$; middle column: absolute values coefficients 
  $\tilde c_n$ after rotations using $180$ samples;
right column: comparison of $\dfrac{\Vert\bm c-\bm c_s\Vert_1}{\sqrt{s}}$
(``$\circ$") and $\dfrac{\Vert\tilde{\bm c}-\tilde{\bm c}_s\Vert_1}{\sqrt{s}}$
(``$\diamond$") with different $s$.}

\label{fig:ex4_coef}
\end{figure}

As an example application of our new method to a more complicated and nonlinear
differential equation, we consider the Korteweg-de Vries (KdV) equation with 
time-dependent additive noise \cite{LinGK06}:
\begin{equation}\label{eq:kdv}
\begin{aligned}
& u_t(x,t;\bx)-6u(x,t;\bx)u_x(x,t;\bx)+u_{xxx}(x,t;\bx)=f(t;\bx), 
  \quad x\in (-\infty,\infty), \\
& u(x,0;\bx) = -2 \sech^2(x).
\end{aligned}
\end{equation}
We model $f(t;\bx)$ as a random field represented by the following KL expansion:
\begin{equation}
f(t;\bx) = \sigma\sum_{i=1}^d\sqrt{\lambda_i}\phi_i(t)\xi_i, 
\end{equation}
where $\sigma$ is a constant and $\{\lambda_i,\phi_i(t)\}_{i=1}^d$ are
eigenpairs of the exponential covariance kernel as in Eqs.~\eqref{eq:kl} and
\eqref{eq:exp_kernel}, respectively. In this problem, we set $l_c=0.25$ and 
$d=10$ ($\sum_{i=1}^d\lambda_i > 0.96\sum_{i=1}^{\infty}\lambda_i$). In this 
case, the exact one-soliton solution is
\begin{equation}\label{eq:kdv_sol}
u(x,t;\bx)=\sigma\sum_{i=1}^d\sqrt{\lambda_i}\xi_i\int_0^t\phi_i(y)\dif y 
-2\sech^2\left(x-4t+6\sigma\sum_{i=1}^d\sqrt{\lambda_i}\xi_i
\int_0^t\int_0^{z}\phi_i(y)\dif y\dif z\right).
\end{equation}
The QoI is chosen to be $u(x,t;\bx)$ at $x=6,t=1$ with 
$\sigma=0.4$. Because an analytical expression for $\phi_i$ is available, 
we can compute the integrals in Eq.~\eqref{eq:kdv_sol} with high accuracy. Denoting
\begin{equation}
A_i = \sqrt{\lambda_i}\int_0^1\phi_i(y)\dif y,\quad 
B_i = \sqrt{\lambda_i}\int_0^1\int_0^{z}\phi_i(y)\dif y\dif z, \quad
i=1,2,\cdots,d,
\end{equation}
the analytical solution is
\begin{equation}\label{eq:kdv_sol2}
u(x,t;\bx)\big |_{x=6,t=1}=\sigma\sum_{i=1}^dA_i\xi_i
-2\sech^2\left(2+6\sigma\sum_{i=1}^d B_i\xi_i\right).
\end{equation}
We use a fourth-order gPC expansion to approximate the solution, i.e., $P=4$, 
and the number of gPC basis functions $N=1001$. The $L_2$ error of the Legendre
and Chebyshev polynomial expansions are presented in Figure~\ref{fig:ex4_rmse}.
In this example, the combined iterative rotation and re-weighted $\ell_1$ 
method outperforms all other approaches. However, in the Chebyshev polynomial
expansion, when the sample size is small (i.e., $M/N<0.12$) and if we only use the
rotational method, the result is not as good as that determined by the standard
$\ell_1$ minimization. This phenomenon also is related to the compromise in the
property of $\tensor\Psi$.
Figure~\ref{fig:ex4_coef} presents a comparison of $\bm c$ and $\tilde{\bm c}$ 
(obtained using $180$ samples), showing the sparsity improvement via
the iterative rotation method. Coefficients with absolute values smaller than
$10^{-8}$ are not presented as they are negligible in the comparison of
sparsity. In addition, Figure~\ref{fig:ex4_coef} demonstrates the sparsity 
enhancement of the sparsity quantitatively.


\subsection{High-dimensional function}
\begin{figure}[h]
\centering
\includegraphics[width=0.45\textwidth]{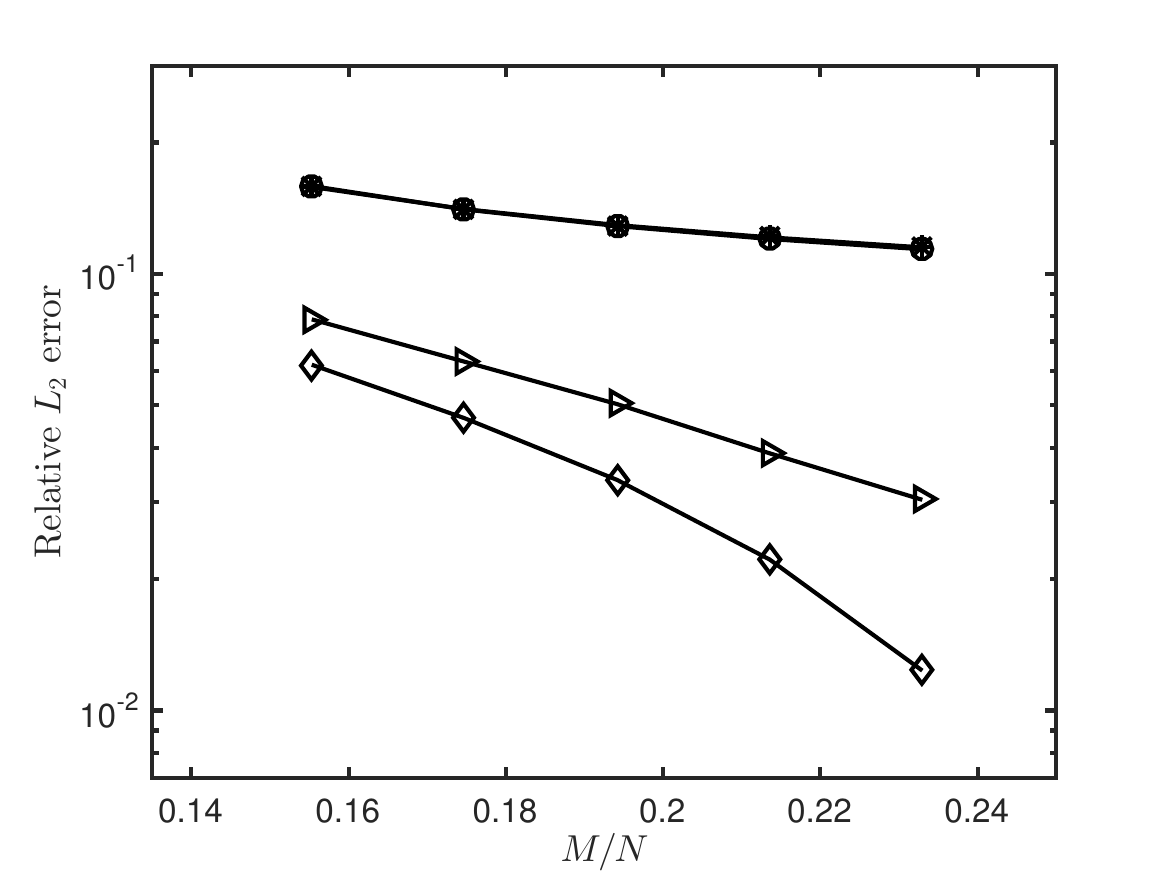}\quad
\includegraphics[width=0.45\textwidth]{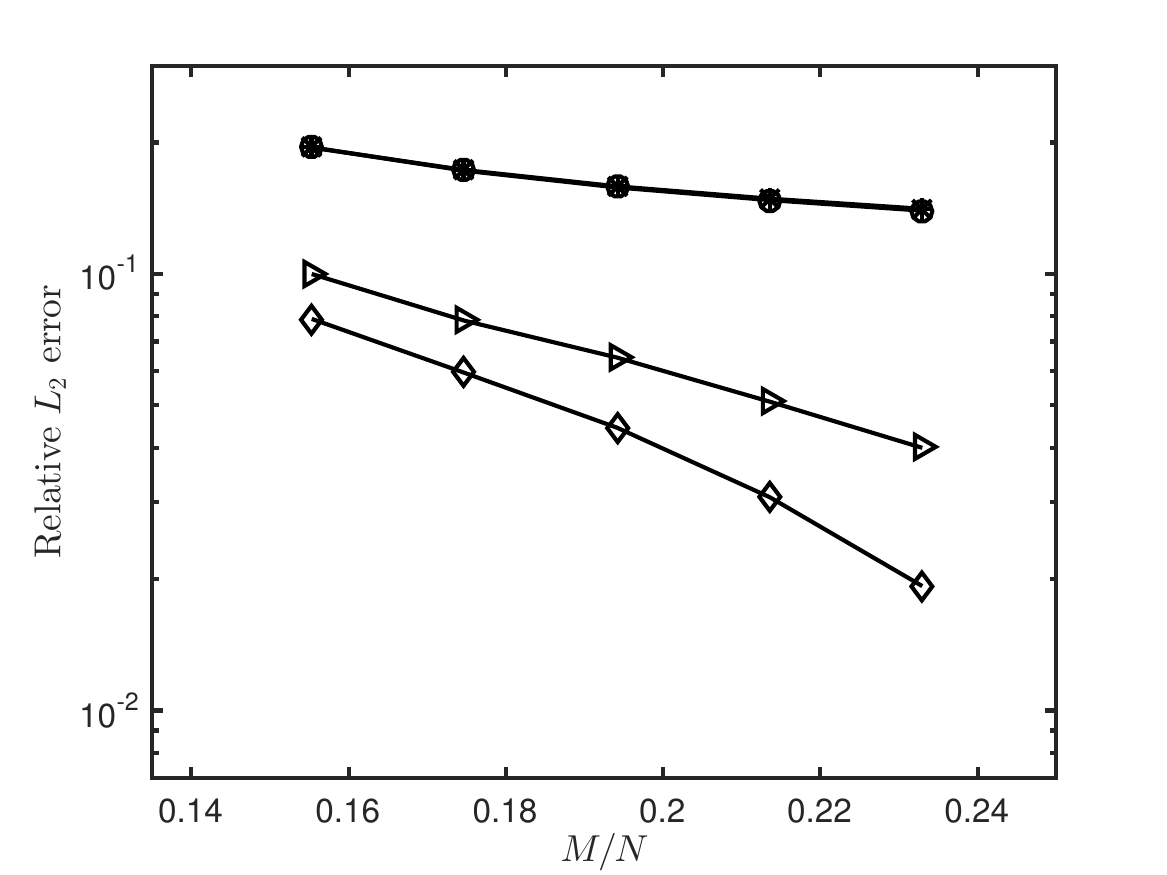}
\caption{Results for the high-dimensional function. Left: Legendre polynomial
expansion (when $\xi_i$ are i.i.d. uniform random variables). Right: Chebyshev
polynomial expansion (when $\xi_i$ are i.i.d. Chebyshev random variables). 
``$\circ$": standard $\ell_1$, ``$\ast$": re-weighted $\ell_1$, 
``$\triangleright$": rotated $\ell_1$, ``$\diamond$": re-weighted+rotated $\ell_1$.}
\label{fig:ex5_rmse}
\end{figure}

\begin{figure}[h]
\centering
\subfigure[Legendre $|c_n|$]
{\includegraphics[width=0.32\textwidth]{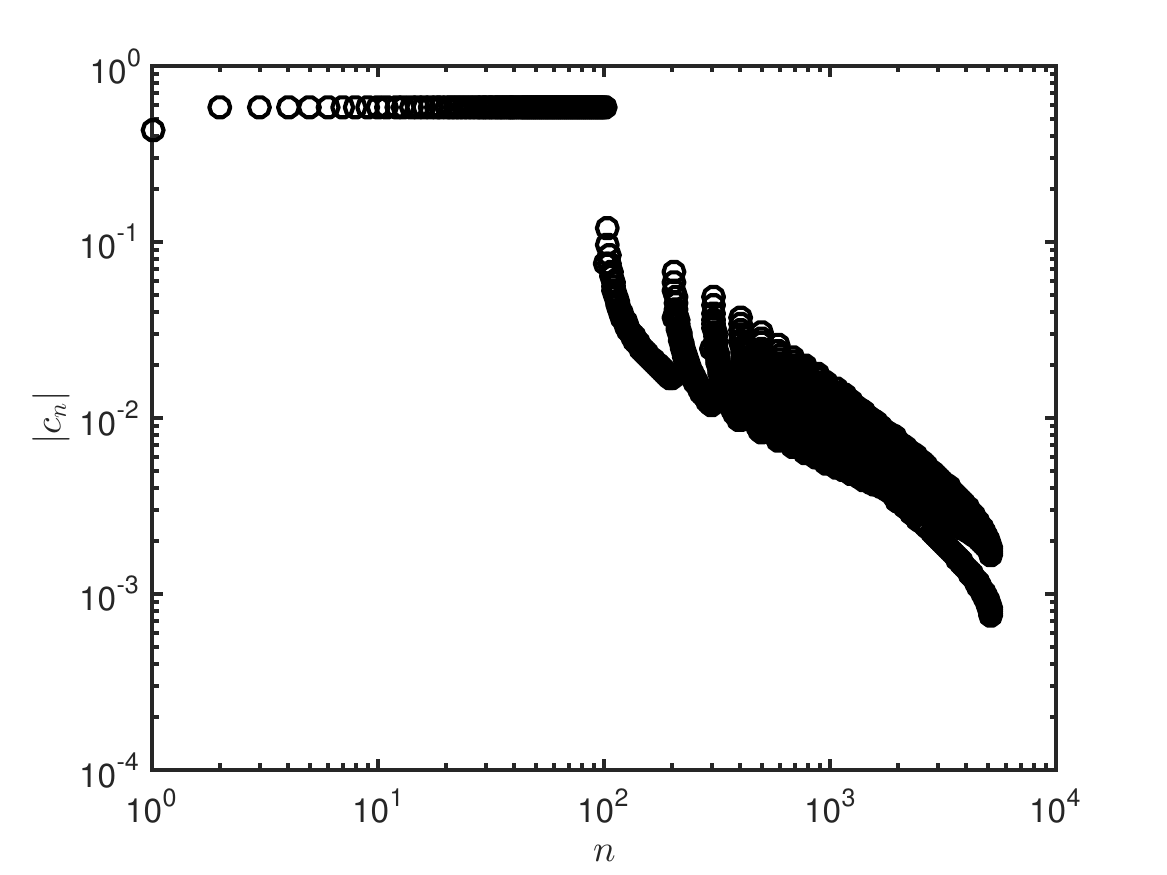}}~
\subfigure[Legendre $|\tilde c_n|$]
{\includegraphics[width=0.32\textwidth]{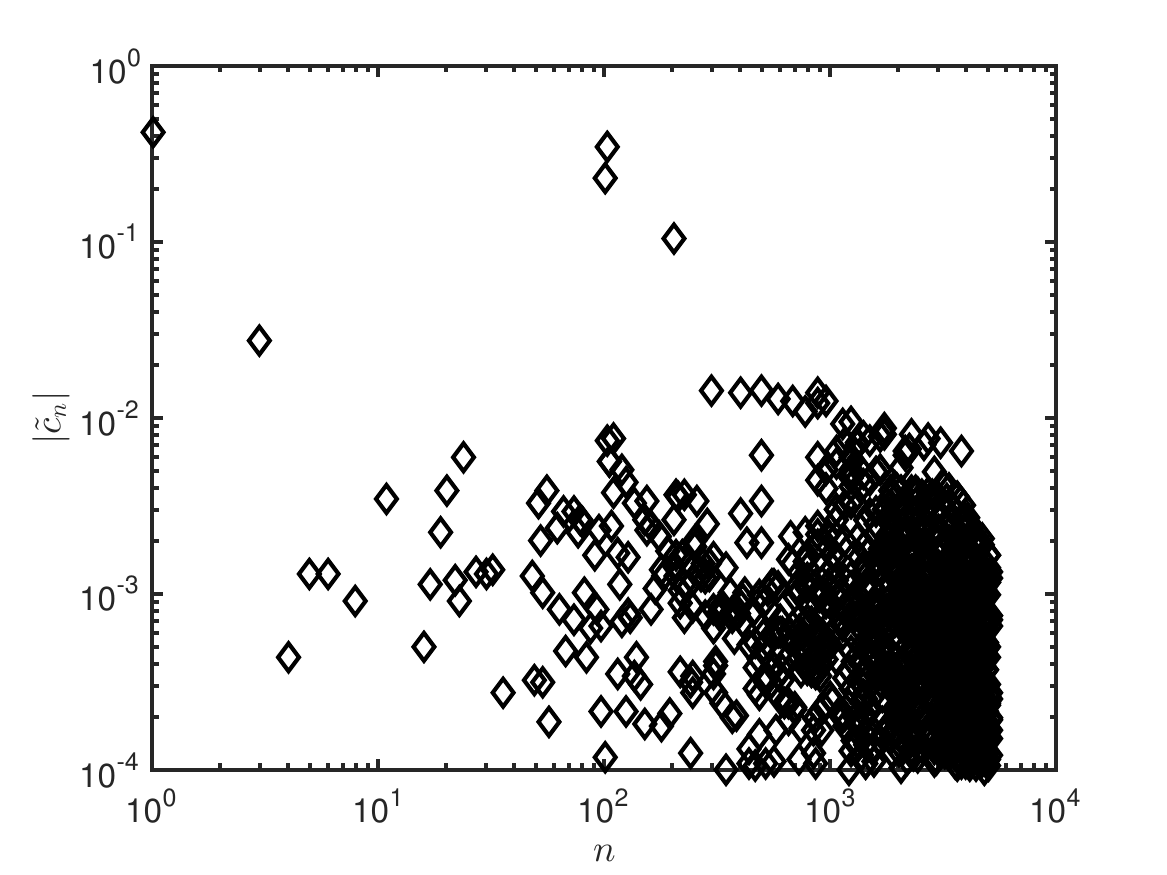}}~
\subfigure[Legendre comparison of sparsity]
{\includegraphics[width=0.32\textwidth]{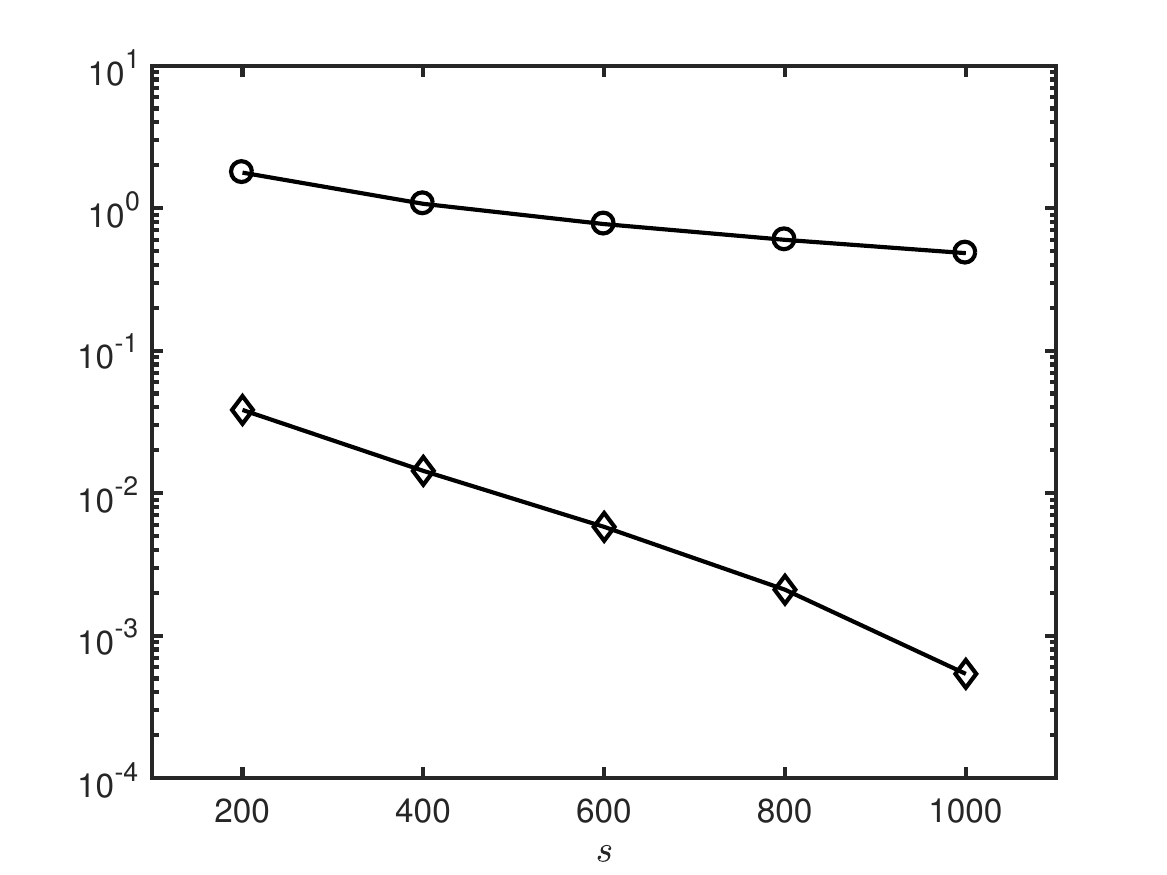}} \\
\subfigure[Chebyshev $|c_n|$]
{\includegraphics[width=0.32\textwidth]{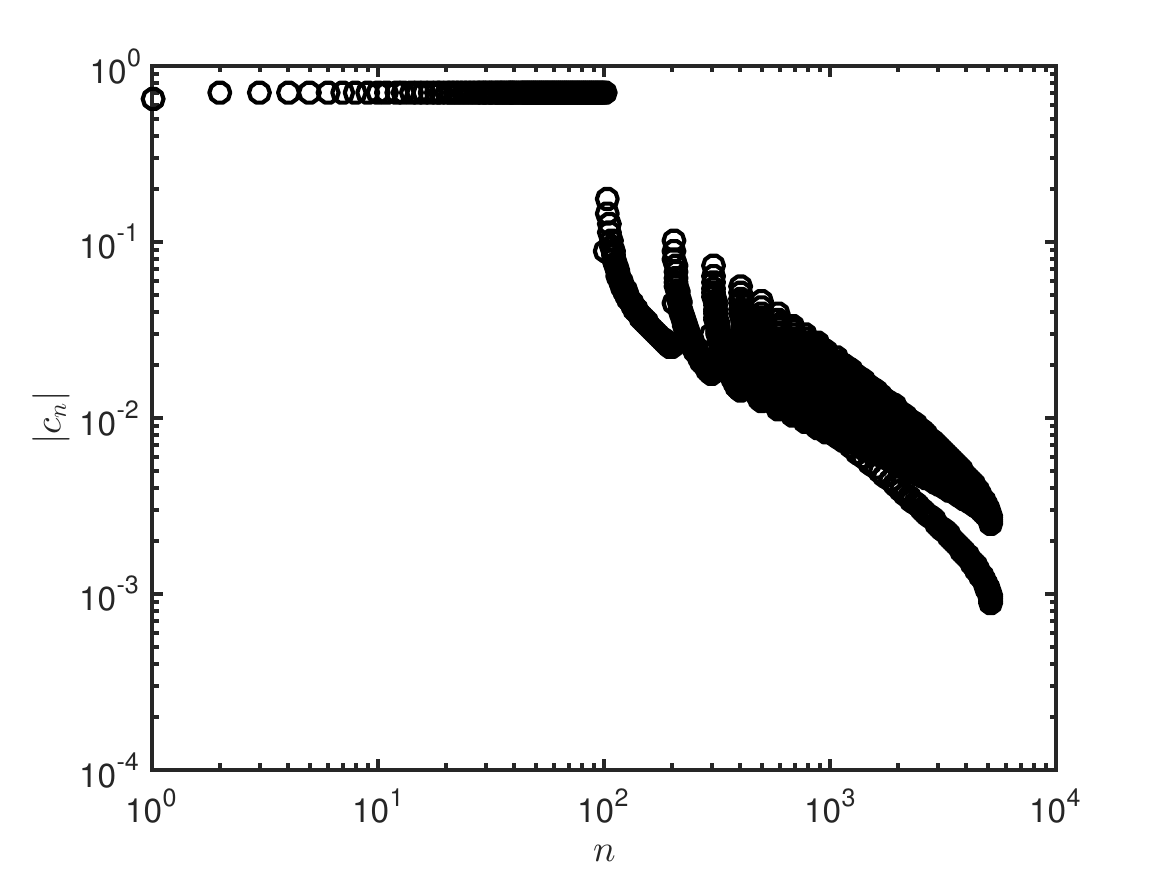}}~
\subfigure[Chebyshev $|\tilde c_n|$]
{\includegraphics[width=0.32\textwidth]{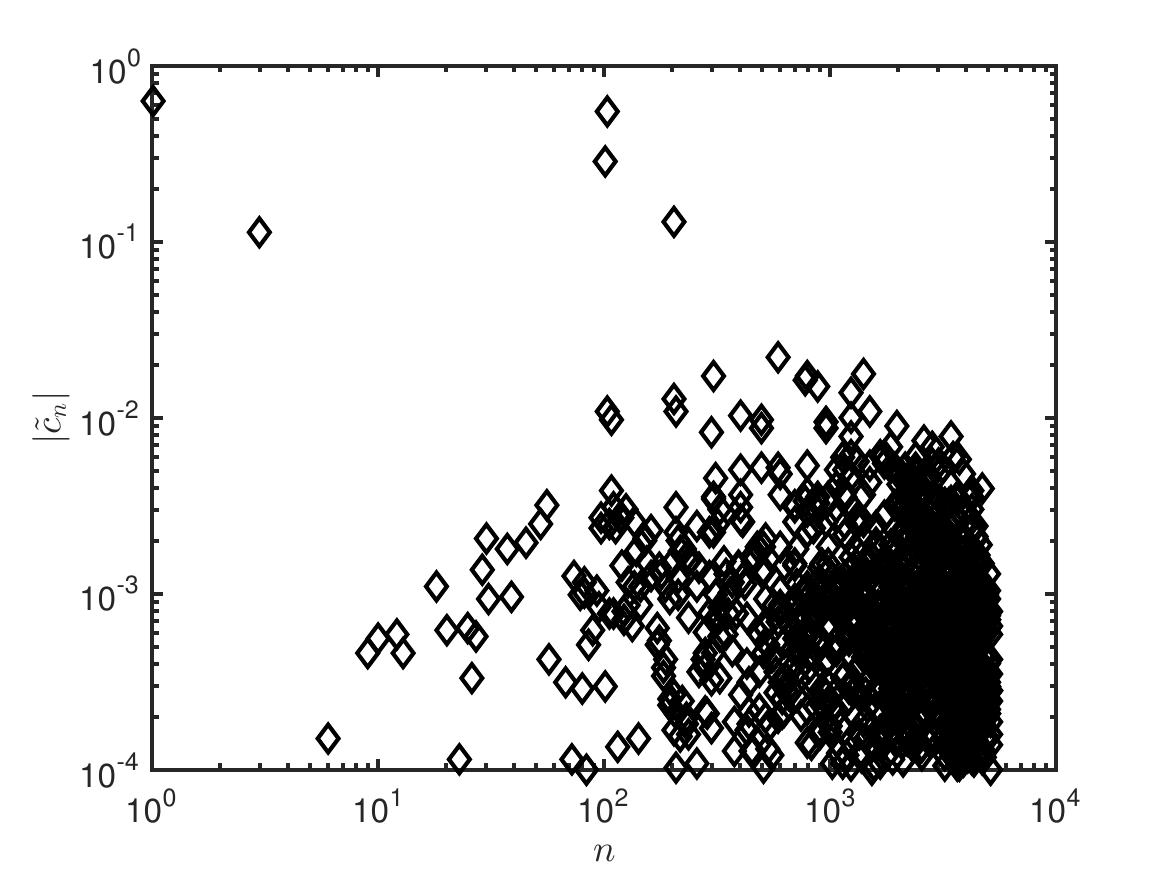}}~
\subfigure[Chebyshev comparison of sparsity]
{\includegraphics[width=0.32\textwidth]{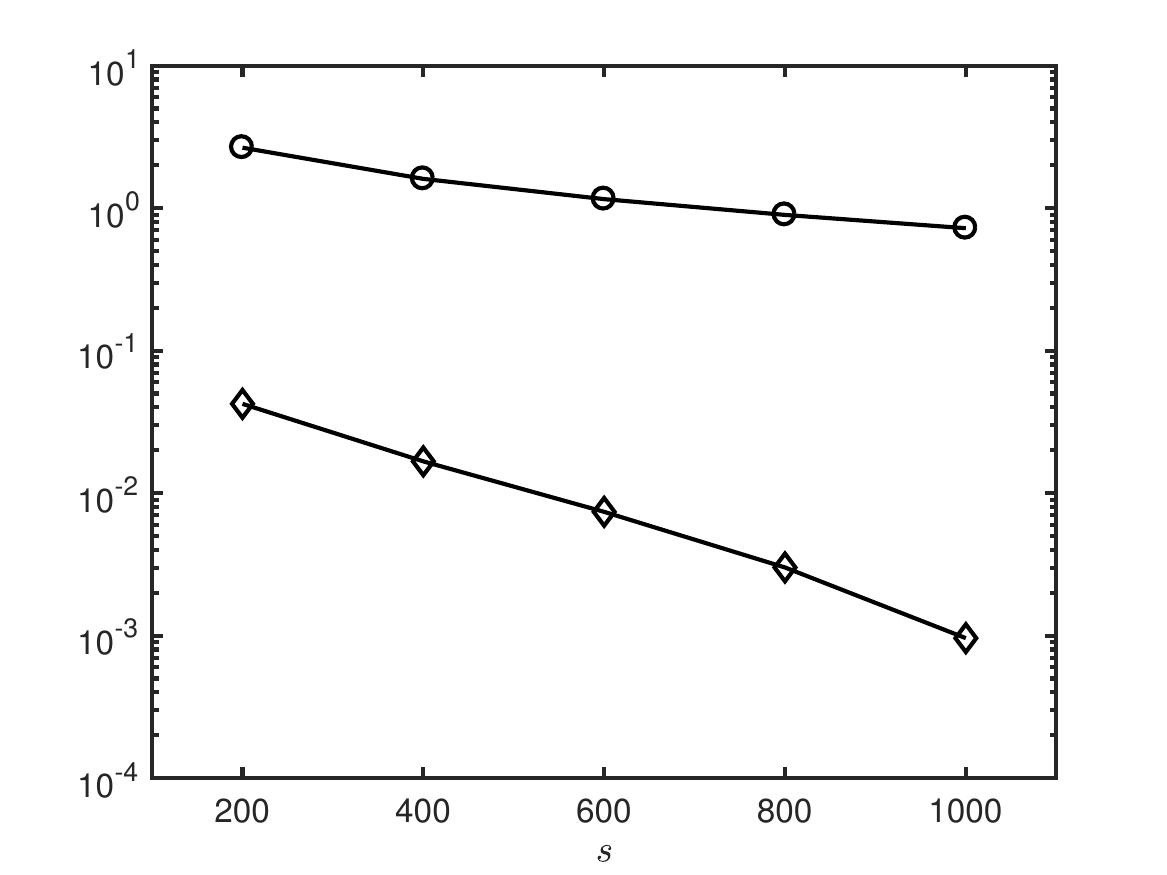}}
\caption{Results for the high-dimensional function. Left column: absolute
  values of exact coefficients $c_n$; middle column: absolute values of 
  coefficients $\tilde c_n$ after rotations using $1200$ samples;
right column: comparison of $\dfrac{\Vert\bm c-\bm c_s\Vert_1}{\sqrt{s}}$
(``$\circ$") and $\dfrac{\Vert\tilde{\bm c}-\tilde{\bm c}_s\Vert_1}{\sqrt{s}}$
(``$\diamond$") with different $s$.}
\label{fig:ex5_coef}
\end{figure}

In this example, we illustrate the potential capability of the rotational
method for dealing with higher-dimensional problems. Specifically, we select a
function similar to the first example (Section~\ref{subsec:ex1}) but with much 
higher dimensionality:
\begin{equation}
u(\bx) = \sum_{i=1}^d \xi_i + 0.25\left(\sum_{i=1}^d \xi_i/\sqrt{i}\right)^2,
  \quad d=100.
\end{equation}
The total number of basis functions for this example is $N=5151$. The relative 
error is computed with a level-$3$ sparse grid method. Hence, the numerical 
integrals are exact. The results are presented in Figure~\ref{fig:ex5_rmse}. As
before, our iterative rotation approach outperforms the existing $\ell_1$ 
methods. Figure~\ref{fig:ex5_coef} features a comparison of $\bm c$ and 
$\tilde{\bm c}$, showing the sparsity enhancement of the sparsity using
the iterative rotation method. Coefficients $\tilde c_n$ (obtained with $1200$
samples) with absolute values smaller than $10^{-4}$ are not presented as they
are two magnitudes smaller than the dominating ones and are negligible in the
comparison of sparsity. The sparsity enhancement is illustrated quantitatively
on the right column of Figure~\ref{fig:ex5_coef}. Notably, for general 
high-dimensional problems, simply truncating the gPC expansion up to a certain
order is not efficient because the number of basis grows exponentially. Hence,
a good approach for high-dimensional problems is to integrate our iterative 
rotation method with a method to reduce $d$ (e.g., ANOVA \cite{YangCLK12}, 
SIR \cite{Li91}) or to reduce $N$ (e.g., adaptive basis selection
\cite{JakemanES14}).

\subsection{Increase in the mutual coherence}
As we point out in Section~\ref{subsec:compromise}, the property of
$\tensor\Psi^{(l)}$ becomes less favorable for $\ell_1$ minimization as we use
Legendre and Chebyshev polynomials in the expansion. Here, we use the mutual
coherence $\mu(\tensor\Psi)$ (see Eq.~\eqref{eq:mutual}) to demonstrate this 
phenomenon. We employ the ridge function setup in example 1 ($d=12, P=3,
N=455$), and we use the exact rotation matrix $\tensor A$ in Eq.~\eqref{eq:ex_rot}
to illustrate how the mutual coherence changes. The matrix $\tensor\Psi$ after 
rotation is computed as $\Psi_{ij}=\psi_j(\tensor A\bx^i)$, where $\tensor A$ is
given in Eq.~\eqref{eq:ex_rot}. We repeat the computing of $\mu(\tensor\Psi)$ 
with $50$ independent sets of $\{\bx^i\}_{i=1}^{180}$ and present the average
value in Table~\ref{tab:coh}. As we expected, for the Hermite polynomial, 
$\mu$ does not change, while, for other types of polynomials, it increases. 
Further, the increase in the Chebyshev polynomial is larger than that of the
Legendre polynomials. This provides a partial explanation as to why the 
rotational method for Chebyshev polynomial is less efficient than for the
Legendre polynomial in some of our test cases. Although theoretical analysis is
not available at this time, Table~\ref{tab:coh} provides an intuitive
understanding of the algorithm.
\begin{table}[h]
  \centering
  \caption{Comparison of mutual coherence before and after rotation for different 
  types of polynomial expansions. $d=12, P=3, N=455, M=180$.}
  \begin{tabular}{C{12em}|C{6em}C{6em}C{6em}}
    \hline \hline
    & Legendre & Chebyshev & Hermite  \\
    \hline
    $\mu(\tensor\Psi)$ before rotation & 0.15 & 0.15 & 0.40 \\
    $\mu(\tensor\Psi)$ after rotation    & 0.45 & 0.50 & 0.40  \\
    \hline \hline
  \end{tabular}
  \label{tab:coh}
\end{table}

\section{Conclusions}
\label{sec:conclusion}

In this work, we extend our previous work on rotation-based \cite{LeiYZLB15} and
iterative-rotation algorithm \cite{YangLBL16} of Hermite polynomial expansion by
providing a general framework for enhancing sparsity of gPC
expansion by using an alternating direction method to identify a rotation 
iteratively. As such, it improves the accuracy of the compressive sensing method
to construct the gPC expansions from a small amount of data. The rotation is 
decided by seeking the directions of maximum variation for the QoI through SVD 
of the gradients at different points in the parameter space. We also demonstrate
that our previously developed iterative method for Hermite polynomial 
expansion \cite{YangLBL16} is a special case of this general framework.  

We combine the iterative rotations with $\ell_1$ minimization. Iterative
rotations also can be integrated with other optimization methods to solve the
compressive sensing problem, e.g., orthogonal matching pursuit
(OMP) \cite{BrucksteinDE09}, $\ell_{1-2}$ minimization \cite{YinLHX15}, and
transformed $\ell_1$ \cite{guo2018stochastic}. {\color{red}Alternatively, it is possible to
alleviate this problem by using the Gram-Schmidt method to reorthogonalize the
measurement matrix \cite{dey2016fuzzy} or using a better sampling and
preconditioning scheme \cite{jakeman2017generalized}.}
In addition, it is possible to further improve the accuracy and reduce the 
number of samples by integrating our method with advanced sampling strategies
(e.g., \cite{HamptonD15}), adaptive basis selection method 
(e.g., \cite{JakemanES14}), Bayesian compressive sensing method 
(e.g.,\cite{KaragiannisBL15}), a better initial guess for the iteration
(e.g.,\cite{yang2017sliced}), etc. This is specifically useful for problems
where the experiments or simulations are costly. The resulting surrogate model
$u_g$ can be used to study parameter sensitivities, and it can be useful in
inverse problems based on the Bayesian framework. 

Along with our previous work, we have demonstrated the effectiveness of the 
rotational method for the Hermite and Legendre polynomial expansions. These two
are the most useful gPC expansions used in the UQ studies as the Gaussian and 
uniform random variables are the most widely used in practice. Our method 
requires fewer samples of QoI to construct surrogate models, which can be a 
great savings of experimental or computational resources. As such, it is useful
for most UQ problems. Of note, the main limitation of applying our method to
other types of gPC expansions is the possible degeneration of the property of 
the measurement matrix. For example, for the Laguerre polynomial expansion, 
numerical tests (not presented in this work) show that our method does not work
for some cases--even if the UQ representation has low dimensional structure. A
systematic numerical analysis on the convergence will be the focus of our future
work. At this time, an intuitive guidance is that our method works well when the
{\color{red}PDF is symmetric (or approximately symmetric) with respect to the origin under 
rotation. }

Finally, as we point out in Section~\ref{subsec:iteration}, removing the
constraint $\tensor A\tensor A\trans=\tensor I$ may result in a more accurate
$u_g$. Also, a nonlinear map from $\bx$ to $\bm\eta$ may work even better.
These improvements allow us to explore low-dimensional structures of the
system more comprehensively, and will be addressed in our future work.

\section*{Acknowledgments}

We would like to thank Dr. Nathan Baker for fruitful discussions.
The research described in this paper was conducted under the Laboratory 
Directed Research and Development Program at Pacific
Northwest National Laboratory, a multiprogram national laboratory operated by
Battelle for the U.S. Department of Energy. PNNL
is operated by Battelle for the DOE under Contract DE-AC05-76RL01830.
A portion of this work was supported by the of Energy, Office of Science, Office 
of Advanced Scientific Computing Research as part of the Collaboratory on 
Mathematics for Mesoscopic Modeling of Materials (CM4).

\bibliographystyle{plain}
\bibliography{uq}

\end{document}